\documentclass[12pt,leqno,twoside]{article} 
\usepackage{bezier,amsmath,amssymb,stmaryrd}

\input xy
\xyoption{all}
\input{rk.sty}
\DeclareMathOperator{\Conc}{Conc}
\newcommand{\CCC}{{\operatorname{CC}}}
\newcommand{\KKK}{{\operatorname{KK}}}
\newcommand{\CCCP}{{\operatorname{CC}^\llcorner}} 
\newcommand{\CCCC}{{\operatorname{CCC}}}
\newcommand{\KKKP}{{\operatorname{KK}^\llcorner}} 
\newcommand{\CCCCE}{{\operatorname{CC}^{\llcorner,\,\operatorname{CE}}}} 
\newcommand{\KKKCE}{{\operatorname{KK}^{\llcorner,\,\operatorname{CE}}}} 
\newcommand{\CCCCP}{{\operatorname{CCC}^{
\linethickness{0.1mm}
\begin{picture}(12,0)(-2,0)
\bezier{10}(0,0)(5,0)(10,0)
\bezier{10}(0,0)(0,5)(0,10)
\bezier{10}(0,0)(3,3)(6,6)
\end{picture}
\thinlines
}}} 
\newcommand{\EEGr}[2]{\EE^{\text{\rm Gr}}_{#1,#2}}
\newcommand{\EEGrp}[2]{\dot\EE^{\text{\rm Gr}}_{#1,#2}}
\newcommand{\EEI}{\EE_\text{\rm I}}
\newcommand{\EEIp}{\dot\EE_\text{\rm I}}
\newcommand{\EEp}{\dot{\EE}}
\newcommand{\tI}{\t_\text{\rm I}}
\newcommand{\Zi}{\Z_\infty}
\newcommand{\bZi}{\b\Z_\infty}
\newcommand{\bZif}{\b\Z_\infty^\#}
\newcommand{\Zif}{\Z_\infty^\#}
\newcommand{\bZiff}{\b\Z_\infty^{\#\#}}
\newcommand{\bZiffp}{\dot{\bar{\Z}}_\infty^{\#\#}}
\newcommand{\dl}{\;\dot{<}\;} 

\begin{document}

\title{Comparison of spectral sequences involving bifunctors}
\author{Matthias K\"unzer}
\maketitle

\begin{small}
\begin{quote}
\begin{center}{\bf Abstract}\end{center}\vspace*{2mm}
Suppose given functors $\Al\ti\Al'\lraa{F}\Bl\lraa{G}\Cl$ between abelian categories, an object $X$ in $\Al$ and an object $X'$ in $\Al'$ such that $F(X,-)$, $F(-,X')$ and $G$ are left exact, and 
such that further conditions hold. We show that, $\EE_1$-terms exempt, the Grothendieck spectral sequence of the composition of $F(X,-)$ and $G$ evaluated at $X'$ is isomorphic to the Grothendieck 
spectral sequence of the composition of $F(-,X')$ and $G$ evaluated at $X$. The respective $\EE_2$-terms are a priori seen to be isomorphic. But instead of trying to compare the differentials and to 
proceed by induction on the pages, we rather compare the double complexes that give rise to these spectral sequences.
\end{quote}
\end{small}

\renewcommand{\thefootnote}{\fnsymbol{footnote}}
\footnotetext[0]{MSC 2000: 18G40.}
\renewcommand{\thefootnote}{\arabic{footnote}}

\begin{footnotesize}
\renewcommand{\baselinestretch}{0.7}%
\parskip0.0ex%
\tableofcontents%
\parskip1.2ex%
\renewcommand{\baselinestretch}{1.0}%
\end{footnotesize}%

\setcounter{section}{-1}

\section{Introduction}

To calculate $\Ext^\ast(X,Y)$, one can either resolve $X$ projectively or $Y$ injectively; the result is, up to isomorphism, the same. To show this, one uses the double complex
arising when one resolves both $X$ and $Y$; cf.\ \bfcite{CE56}{Chap.\ V, Th.\ 8.1}. 

Two problems in this spirit occur in the context of Grothendieck spectral sequences; cf.\ \mb{\S\S\,\ref{SecIntroFirst}, \ref{SecIntroSecond}}.

\subsection{Language}

In \S\ref{SecForm}, we give a brief introduction to 
the Deligne-Verdier spectral sequence language; cf.\ \mb{\bfcite{Ve67}{II.\S 4}}, \bfcite{De94}{App.}; or, on a more basic level, cf.\ \bfcite{Ku06}{Kap.\ 4}. This language amounts to considering 
a diagram $\EE(X)$ containing all the images between the homology groups of the subquotients of a given filtered complex $X$, 
instead of, as is classical, only selected ones. This helps to gain some elbow room in practice\,: to govern the objects of the diagram $\EE(X)$ we can make use of a certain 
short exact sequence; cf.\ \S\ref{SecSES}. 

Dropping the $\EE_1$-terms and similar ones, we obtain the {\it proper} spectral sequence $\dot\EE(X)$ of our
filtered complex $X$. Amongst others, it contains all $\EE_k$-terms for $k\ge 2$ in the classical language; cf.\ \S\S\,\ref{SecCompProp}, \ref{SecClassical}.

\subsection{First comparison}
\label{SecIntroFirst}

Suppose given abelian categories $\Al$, $\Al'$ and $\Bl$ with enough injectives and an abelian category $\Cl$. Suppose given objects $X\in\Ob\Al$ and $X'\in\Ob\Al'$. 
Let $\Al\ti\Al'\lraa{F}\Bl$ be a biadditive functor such that $F(X,-)$ and $F(-,X')$ are left exact. Let $\Bl\lraa{G}\Cl$ be a left exact functor. Suppose further 
conditions to hold; see \S\ref{SecFirstCompIso}.

We have a Grothendieck spectral sequence for the composition $G\0 F(X,-)$ and a Grothendieck spectral sequence for the composition $G\0 F(-,X')$.
We evaluate the former at $X'$ and the latter at $X$. 

In both cases, the $\EE_2$-terms are $(\RR^i G)(\RR^j F)(X,X')$. Moreover, they both converge to 
\mb{$\big(\RR^{i+j} (G\0 F)\big)(X,X')$}. So the \mb{following} assertion is well-motivated.

{\bf Theorem \ref{ThFC1}.} {\it The proper Grothendieck spectral sequences just described are isomorphic; i.e.\
\[
\EEGrp{F(X,-)}{G}(X') \;\iso\; \EEGrp{F(-,X')}{G}(X)\; .
\]
}%
So instead of ``resolving $X'$ twice'', we may just as well ``resolve $X$ twice''.

In fact, the underlying double complexes are connected by a chain of double homotopisms, i.e.\ isomorphisms in the homotopy category as defined in \bfcite{CE56}{IV.\S4}, and rowwise homotopisms
(the proof uses a chain $\bt\llafl{25}{\text{double}} \bt\llafl{25}{\text{roww.}}\bt\lrafl{25}{\text{roww.}}\bt\lrafl{25}{\text{double}}\bt$).
These morphisms then induce isomorphisms on the associated proper first spectral sequences.

\subsection{Second comparison}
\label{SecIntroSecond}

Suppose given abelian categories $\Al$ and $\Bl'$ with enough injectives and abelian categories $\Bl$ and $\Cl$. Suppose given objects $X\in\Ob\Al$ and $Y\in\Ob\Bl$.
Let $\Al\lraa{F} \Bl'$ be a left exact functor. Let $\Bl\ti\Bl'\lraa{G} \Cl\ru{5}$ be a biadditive functor such that $G(Y,-)$ is left exact.

Let $B\in\Ob\CC^{[0}(\Bl)$ be a resolution of $Y$, i.e.\ a complex $B$ admitting a quasiisomorphism $\Conc Y\lra B$. Suppose that
$G(B^k,-)$ is exact for all $k\ge 0$. Let $A\in\Ob\CC^{[0}(\Al)$ be, say, an injective resolution of $X$. Suppose further conditions to hold; see \S\ref{SecSecondCompIso}.

We have a Grothendieck spectral sequence for the composition $G(Y,-)\0 F$, which we evaluate at $X$. 
On the other hand, we can consider the double complex $G(B,FA)$, where the indices of $B$ count rows and the indices of $A$ count columns. To the first filtration of its 
total complex, we can associate the proper spectral sequence $\EEIp\big(G(B,FA)\big)$. 

If $\Bl$ has enough injectives and $B$ is an injective resolution of $Y$, then in both cases the $\EE_2$-terms are a priori seen to be $(\RR^i G)\big(Y,(\RR^j F)(X)\big)$. So also the following assertion
is well-motivated.

{\bf Theorem \ref{ThSC1}.} {\it We have $\;\EEGrp{F}{G(Y,-)}(X) \;\iso\; \EEIp\big(G(B,FA)\big)\;$.}

So instead of ``resolving $X$ twice'', we may just as well ``resolve $X$ once and $Y$ once''.

The left hand side spectral sequence converges to $\big(\RR^{i+j}(G(Y,-)\circ F)\big)(X)$. By this theorem, so does the right hand side one. 

The underlying double complexes are connected by two morphisms of double complexes \mb{(in the directions $\bt\lra\bt\lla\bt$)} that induce isomorphisms on the associated proper spectral sequences.

Of course, Theorems \ref{ThFC1} and \ref{ThSC1} have dual counterparts.

\subsection{Results of Beyl and Barnes}

Let $R$ be a commutative ring. Let $G$ be a group. Let $N\nrm G$ be a normal subgroup. Let $M$ be an $RG$-module.

{\sc Beyl} generalises Grothendieck's setup, allowing for a variant of a Cartan-Eilenberg resolution that consists of acyclic, but no longer necessarily injective objects~\bfcite{Be81}{Th.\ 3.4}. 
We have documented {\sc Beyl}'s Theorem as Theorem \ref{ThAC1} in our framework, without claiming originality.

{\sc Beyl} uses his Theorem to prove that, from the $\EE_2$-term on, the Grothendieck spectral sequence for 
$RG\Modl\lrafl{26}{(-)^N} RN\Modl\lrafl{26}{(-)^{G/N}} R\Modl\ru{5.5}$ at $M$ is isomorphic to the Lyndon-Hochschild-Serre spectral sequence, i.e.\ the spectral sequence associated to 
the double complex $\liu{RG}{(\BBar_{G/N;R}\ts_R\BBar_{G;R}\,,\, M)}$; cf.\ \bfcite{Be81}{Th.\ 3.5}, \bfcite{Ben91}{\S 3.5}. This is now also a consequence 
of Theorems \ref{ThFC1} and \ref{ThSC1}, as explained in \S\S\,\ref{SecApplHSH}, \ref{SecApplLHS}. 

{\sc Barnes} works in a slightly different setup. He supposes given a commutative ring $R$, abelian categories $\Al$, $\Bl$ and $\Cl$ of $R$-modules, and left exact functors 
$F:\Al\lra\Bl$ and $G:\Bl\lra\Cl$, where $F$ is supposed to have an exact left adjoint $J:\Bl\lra\Al$ that satisfies $F\0 J = 1_\Bl$. Moreover, he assumes $\Al$ to have ample injectives 
and $\Cl$ to have enough injectives. In this setup, he obtains a general comparison theorem. See \bfcite{Ba85}{Sec.\ X.5, Def.\ X.2.5, Th.\ X.5.4}.

{\sc Beyl} \bfcit{Be81} and {\sc Barnes} \bfcit{Ba85} also consider cup products; in this article, we do not. 

\subsection{Acknowledgements}

Results of {\sc Beyl} and {\sc Haas} are included for sake of documentation that they work within our framework; cf.\ Theorem \ref{ThAC1} and \S\ref{SecGroth}. 
No originality from my part is claimed.

I thank {\sc B.~Keller} for directing me to \bfcite{ML75}{XII.\S 11}. I thank the referee for helping to considerably improve the presentation, and for suggesting 
Lemma~\ref{LemCrux} and \S\ref{SecApplHSH}. I thank {\sc G.~Carnovale} and {\sc G.~Hi\ss} for help with Hopf algebras.

\begin{footnotesize}
\vspace*{5mm}
{\bf\large Conventions}
\vs

Throughout these conventions, let $\Cl$ and $\Dl$ be categories, let $\Al$ be an additive category, let $\Bl$ and $\Bl'$ be abelian categories, and let $\El$ be an exact category in which all idempotents
split.

\begin{itemize}
\item For $a,\, b\,\in\,\Z$, we write $[a,b] := \{c\in\Z\;:\; a\le c\le b\}$, $[a,b[\; := \{c\in\Z\;:\; a\le c < b\}$, etc.
\item Given $I\tm\Z$ and $i\in\Z$, we write $I_{\ge i} := \{ j\in I\; :\; j\ge i\}$ and $I_{< i} := \{ j\in I\; :\; j < i\}$.
\item The disjoint union of sets $A$ and $B$ is denoted by $A\disj B$.
\item Composition of morphisms is written on the right, i.e.\ $\lraa{a}\lraa{b} = \lraa{ab}$.\vspace*{-2mm} 
\item Functors act on the left. Composition of functors is written on the left, i.e.\ $\lraa{F}\lraa{G} = \lrafl{25}{G\0 F}$
\item Given objects $X$, $Y$ in $\Cl$, we denote the set of morphisms from $X$ to $Y$ by $\liu{\Cl}{(X,Y)}$.
\item The category of functors from $\Cl$ to $\Dl$ and transformations between them is denoted by $\fbo \Cl,\Dl \fbc$.
\item Denote by $\CC(\Al)$ the category of complexes 
\[
X \= (\cdots\;\lraa{d}\; X^{i-1}\;\lraa{d}\; X^i\;\lraa{d}\; X^{i+1}\;\lraa{d}\cdots)
\]
with values in $\Al$. Denote by $\CC^{[0}(\Al)$ the full subcategory of $\CC(\Al)$ consisting of complexes $X$ with $X^i = 0$ for $i < 0$. We have a full embedding
$\Al \lrafl{25}{\Conc} \CC^{[0}(\Al)$, where, given $X\in\Ob\Al$, the complex $\Conc X$ has entry $X$ at position $0$ and zero elsewhere. 
\item Given a complex $X\in\Ob\CC(\Al)$ and $k\in\Z$, we denote by $X^{\bt+k}$ the complex that has differential $X^{i+k}\lrafl{20}{(-1)^k d} X^{i+1+k}$ between positions
$i$ and $i+1$. We also write $X^{\bt-1} := X^{\bt+(-1)}$ etc.
\item Suppose given a full additive subcategory $\Ml\tm\Al$. Then $\Al/\Ml$ denotes the quotient of $\Al$ by $\Ml$, which has the same objects as $\Al$, and which has as morphisms
residue classes of morphisms of $\Al$, where two morphisms are in the same residue class if their difference factors over an object of $\Ml$.
\item A morphism in $\Al$ is {\it split} if it isomorphic, as a diagram on $\bt\lra\bt$, to a morphism of the form \mb{$X\ds Y\lrafl{36}{\smatzz{1}{0}{0}{0}} X\ds Z\ru{6.5}$}.
A complex $X\in\Ob\CC(\Al)$ is {\it split} if all of its differentials are split.
\item An {\it elementary split acyclic} complex in $\CC(\Al)$ is a complex of the form 
\[
\cdots\lra 0\lra T\lraa{1} T\lra 0\lra\cdots\;, 
\]
where the entry $T$ is at positions $k$ and $k+1$ for some $k\in\Z$.
A {\it split acyclic} complex is a complex isomorphic to a direct sum of elementary split acyclic complexes, i.e.\ a complex isomorphic to a complex of the form
\[
\cdots \lrafl{35}{\smatzz{0}{0}{1}{0}} T^i\ds T^{i+1} 
\lrafl{35}{\smatzz{0}{0}{1}{0}} T^{i+1}\ds T^{i+2} 
\lrafl{35}{\smatzz{0}{0}{1}{0}} T^{i+2}\ds T^{i+3} 
\lrafl{35}{\smatzz{0}{0}{1}{0}} \cdots
\]
Let $\CC_\text{sp\,ac}(\Al)\tm\CC(\Al)$ denote the full additive subcategory of split acyclic complexes. Let $\KK(\Al) := \CC(\Al)/\CC_\text{sp\,ac}(\Al)$ denote the homotopy category of complexes
with values in $\Al$. Let $\KK^{[0}(\Al)$ denote the image of $\CC^{[0}(\Al)$ in $\KK(\Al)$. A morphism in $\CC(\Al)$ is a {\it homotopism} if its image in $\KK(\Al)$ is an isomorphism.
\item We denote by $\Inj\Bl\tm\Bl$ the full subcategory of injective objects. 
\item Concerning exact categories, introduced by {\sc Quillen} \bfcite{Qu73}{p.~15}, we use the conventions of \mb{\bfcite{Ku05}{Sec.\ A.2}}. In particular, a commutative quadrangle in $\El$ 
being a pullback is indicated by
\[
\xymatrix{
A\ar[r]\ar[d]\ar@{}[dr]|(0.2)\hookpbo & B\ar[d] \\
C\ar[r]                               & D \; ,\!\!\\
}
\]
a commutative quadrangle being a pushout by
\[
\xymatrix{
A\ar[r]\ar[d]\ar@{}[dr]|(0.8)\hookpou & B\ar[d] \\
C\ar[r]                               & D \; .\!\!\\
}
\]
\item Given $X\in\Ob\CC(\El)$ with pure differentials, and given $k\in\Z$, we denote by $\ZZ^k X$ the kernel of the differential $X^k\lra X^{k+1}$, by $\ZZ'^k X$ the cokernel of the 
differential $X^{k-1}\lra X^k$, and by $\BB^k X$ the image of the differential $X^{k-1}\lra X^k$. Furthermore, we have pure short exact sequences $\BB^k X\lramono \ZZ^k X\lraepi \HH^k X$
and $\HH^k X\lramono\ZZ'^k X\lraepi\BB^{k+1} X$. 
\item A morphism $X\lra Y$ in $\CC(\El)$ between complexes $X$ and $Y$ with pure differentials is a {\it quasiisomorphism} if $\HH^k$ applied to it yields an isomorphism for all $k\in\Z$. A 
complex $X$ with pure differentials is {\it acyclic} if $\HH^k X\iso 0$ for all $k\ge 0$. Such a complex is also called a {\it purely acyclic} complex.
\item Suppose that $\Bl$ has enough injectives. Given a left exact functor $\Bl\lraa{F}\Bl'$, an object $X\in\Ob\Bl$ is {\it $F$-acyclic} if $\RR^i F X\iso 0$ for all $i\ge 1$. In other words,
$X$ is $F$-acyclic if for an injective resolution $I\in\CC^{[0}(\Inj\Bl)$ of $X$ (and then for all such injective resolutions), we have $\HH^i FI\iso 0$ for all $i\ge 1$.
\item By a module, we understand a left module, unless stated otherwise. If $A$ is a ring, we abbreviate $\liu{A}{(-,=)} := \liu{A\Modl}{(-,=)} = \Hom_A(-,=)$.
\end{itemize}
\end{footnotesize}

\section{Double and triple complexes}

\bq
 We fix some notations and sign conventions.

\eq

Let $\Al$ and $\Bl$ be additive categories. Let $\CC(\Al)\lraa{H}\Bl$ be an additive functor. 

\subsection{Double complexes}
\label{SecDoubleComplexes}

\subsubsection{Definition}

A {\it double complex} with entries in $\Al$ is a diagram
\[
\xymatrix{
      & &                & \vdots                             & \vdots                               & \vdots                               &        \\
      & & \cdots\ar[r]^d & X^{i+2,j}\ar[r]^d\ar[u]^(0.4)\dell & X^{i+2,j+1}\ar[r]^d\ar[u]^(0.4)\dell & X^{i+2,j+2}\ar[r]^d\ar[u]^(0.4)\dell & \cdots \\
X\; = & & \cdots\ar[r]^d & X^{i+1,j}\ar[r]^d\ar[u]^\dell      & X^{i+1,j+1}\ar[r]^d\ar[u]^\dell      & X^{i+1,j+2}\ar[r]^d\ar[u]^\dell      & \cdots \\
      & & \cdots\ar[r]^d & X^{i,j}\ar[r]^d\ar[u]^\dell        & X^{i,j+1}\ar[r]^d\ar[u]^\dell        & X^{i,j+2}\ar[r]^d\ar[u]^\dell        & \cdots \\
      & &                & \vdots\ar[u]^\dell                 & \vdots\ar[u]^\dell                   & \vdots\ar[u]^\dell                   &        \\
}
\]
in $\Al$ such that $dd = 0$, $\dell\dell = 0$ and $d\dell = \dell d$ everywhere. As morphisms between double complexes, we take all diagram morphisms.
Let $\CCC(\Al)$ denote the category of double complexes. We may identify $\CCC(\Al) = \CC(\CC(\Al))$. 

The double complexes considered in this \S\ref{SecDoubleComplexes} are stipulated to have entries in $\Al$.

Let $\CCCP(\Al) := \CC^{[0}(\CC^{[0}(\Al))$ be the category of {\it first quadrant double complexes}, 
consisting of double complexes $X$ such that $X^{i,j} = 0$ whenever $i < 0$ or $j < 0$. 

Given a double complex $X$ and $i\in\Z$, we let $X^{i,\ast}\in\Ob\CC(\Al)$ denote the complex that has entry $X^{i,j}$ at position $j\in\Z$, the differentials taken accordingly; $X^{i,\ast}$ is 
called the $i$th {\it row} of $X$.

Similarly, given $j\in\Z$, $X^{\ast,j}\in\Ob\CC(\Al)$ denotes the $j$th {\it column} of $X$.

\subsubsection{Applying $H$ in different directions}

Given $X\in\Ob\CCC(\Al)$, we let $H(X^{\ast,-})\in\Ob\CC(\Al)$ denote the complex that has $H(X^{\ast,j})$ at position $j\in\Z$, and as differential $H(X^{\ast,j})\lra H(X^{\ast,j+1})$ the image
of the morphism $X^{\ast,j}\lra X^{\ast,j+1}$ of complexes under $H$. Similarly, $H(X^{-,\ast})\in\Ob\CC(\Al)$ has $H(X^{j,\ast})$ at position $j\in\Z$.

\bq
 In other words, a ``$\ast$'' denotes the index direction to which $H$ is applied, a ``$-$'' denotes the surviving index direction. For short, ``$\ast$'' before ``$-$''.

\eq

\subsubsection{Concentrated double complexes}

Given a complex $U\in\Ob\CC^{[0}(\Al)$, we denote by $\Conc_2 U\in\Ob\CCCP(\Al)$ the double complex whose $0$th row is given by $U$, and whose other rows are zero; i.e.\ given $j\in\Z$, then
$(\Conc_2 U)^{i,j}$ equals $U^j$ if $i = 0$, and $0$ otherwise, the differentials taken accordingly. Similarly, $\Conc_1 U\in\Ob\CCCP(\Bl)$ denotes the double complex whose $0$th column is given 
by $U$, and whose other columns are zero.

\subsubsection{Row- and columnwise notions}

A morphism $X\lraa{f} Y$ of double complexes is called a {\it rowwise homotopism} if $X^{i,\ast}\lraa{f^{i,\ast}} Y^{i,\ast}$ is a homotopism for all $i\in\Z$.
Provided $\Al$ is abelian, it is called a {\it rowwise quasiisomorphism} if $X^{i,\ast}\lraa{f^{i,\ast}} Y^{i,\ast}$ is a quasiisomorphism for all $i\in\Z$.

A morphism $X\lraa{f} Y$ of double complexes is called a {\it columnwise homotopism} if $X^{\ast,j}\lraa{f^{\ast,j}} Y^{\ast,j}$ is a homotopism for all $j\in\Z$.
Provided $\Al$ is abelian, it is called a {\it columnwise quasiisomorphism} if $X^{\ast,j}\lraa{f^{\ast,j}} Y^{\ast,j}\ru{5}$ is a quasiisomorphism for all $j\in\Z$.

Provided $\Al$ is abelian, a double complex $X$ is called {\it rowwise split} if $X^{i,\ast}$ is split for all $i\in\Z$; a short exact sequence $X'\lra X\lra X''$ of double complexes
is called {\it rowwise split short exact} if $X'^{i,\ast}\lra X^{i,\ast}\lra X''^{i,\ast}$ is split short exact for all $i\in\Z$.

A double complex $X$ is called {\it rowwise split acyclic} if $X^{i,\ast}$ is a split acyclic complex for all $i\in\Z$. It is called {\it columnwise split acyclic} 
if $X^{\ast,j}$ is a split acyclic complex for all $j\in\Z$. 

\subsubsection{Horizontally and vertically split acyclic double complexes}

An {\it elementary horizontally split acyclic} double complex is a double complex of the form
\[
\xymatrix@R=12mm@C=12mm{
             & \vdots         & \vdots                       & \vdots                   & \vdots         &        \\
\cdots\ar[r] & 0\ar[u]\ar[r]  & T^{i+1}\ar@{=}[r]\ar[u]      & T^{i+1}\ar[r]\ar[u]      & 0\ar[u]\ar[r]  & \cdots \\
\cdots\ar[r] & 0\ar[u]\ar[r]  & T^i \ar@{=}[r]\ar[u]^{\dell} & T^i \ar[r]\ar[u]^{\dell} & 0\ar[u]\ar[r]  & \cdots \\
             & \vdots\ar[u]   & \vdots\ar[u]                 & \vdots\ar[u]             & \vdots\ar[u]   & .      \\ 
}
\]
A {\it horizontally split acyclic} double complex is a double complex isomorphic to a direct sum of elementary horizontally split acyclic double complexes, i.e.\ to one of the form
\[
\xymatrix@R=12mm@C=12mm{
             & \vdots                                                                                      & \vdots                                                           &        \\
\cdots\ar[r] & T^{i+1,j}\dk T^{i+1,j+1} \ar[r]^(0.47){\smatzz{0}{0}{1}{0}}\ar[u]                           & T^{i+1,j+1}\dk T^{i+1,j+2}\ar[r]\ar[u]                           & \cdots \\
\cdots\ar[r] & T^{i,j}\dk T^{i,j+1} \ar[r]^(0.48){\smatzz{0}{0}{1}{0}}\ar[u]^{\smatzz{\dell}{0}{0}{\dell}} & T^{i,j+1}\dk T^{i,j+2}\ar[r]\ar[u]^{\smatzz{\dell}{0}{0}{\dell}} & \cdots \\
             & \vdots\ar[u]                                                                                & \vdots\ar[u]                                                     & .      \\ 
}
\]
An {\it elementary vertically split acyclic} double complex is a double complex of the form
\[
\xymatrix@R=12mm@C=12mm{
             & \vdots                 & \vdots                   &        \\
\cdots\ar[r] & 0\ar[r]\ar[u]          & 0\ar[r]\ar[u]            & \cdots \\
\cdots\ar[r] & T^i\ar[r]^d\ar[u]      & T^{i+1}\ar[r]\ar[u]      & \cdots \\
\cdots\ar[r] & T^i \ar[r]^d\ar@{=}[u] & T^{i+1}\ar[r]\ar@{=}[u]  & \cdots \\
\cdots\ar[r] & 0\ar[r]\ar[u]          & 0\ar[r]\ar[u]            & \cdots \\
             & \vdots\ar[u]           & \vdots\ar[u]             & .      \\ 
}
\]
A {\it vertically split acyclic} double complex is a double complex isomorphic to a direct sum of elementary vertically split acyclic double complexes, i.e.\ to one of the form
\[
\xymatrix@R=12mm@C=12mm{
             & \vdots                                                                                  & \vdots                                                         &        \\
\cdots\ar[r] & T^{i+1,j}\dk T^{i+2,j} \ar[r]^(0.47){\smatzz{d}{0}{0}{d}}\ar[u]                         & T^{i+1,j+1}\dk T^{i+2,j+1}\ar[r]\ar[u]                         & \cdots \\
\cdots\ar[r] & T^{i,j}\dk T^{i+1,j} \ar[r]^(0.48){\smatzz{d}{0}{0}{d}}\ar[u]^{\smatzz{0}{0}{1}{0}}     & T^{i,j+1}\dk T^{i+1,j+1}\ar[r]\ar[u]^{\smatzz{0}{0}{1}{0}}     & \cdots \\
             & \vdots\ar[u]                                                                            & \vdots\ar[u]                                                   & .      \\ 
}
\]

A horizontally split acyclic double complex is in particular rowwise split acyclic. A vertically split acyclic double complex is in particular columnwise split acyclic. 

A double complex is called {\it split acyclic} if it is isomorphic to the direct sum of a horizontally and a vertically split acyclic double complex. Let $\CCC_\text{sp\,ac}(\Al)$ denote the
full additive subcategory of split acyclic double complexes. Let 
\[
\KKK(\Al) \; :=\; \CCC(\Al)/\CCC_\text{sp\,ac}(\Al)\; ; 
\]
cf.\ \bfcite{CE56}{IV.\S4}. A morphism in $\CCC(\Al)$ that is mapped to an isomorphism in $\KKK(\Al)$ is called a {\it double homotopism}.

\bq
 A speculative aside. The category $\KK(\Al)$ is Heller triangulated; cf.\ \mb{\bfcite{Ku05}{Def.\ 1.5.(i), Th.\ 4.6}}. Such a Heller triangulation hinges on two induced shift functors, one of 
 them induced by the shift functor on $\KK(\Al)$. Now $\KKK(\Al)$ carries two shift functors, and so there might be more isomorphisms between induced shift functors one can fix. How can the formal 
 structure of $\KKK(\Al)$ be described?

\eq

\subsubsection{Total complex}
\label{SecTotal}

Let $\KKKP(\Al)$ be the full image of $\CCCP(\Al)$ in $\KKK(\Al)$.

The {\it total complex} $\t X$ of a double complex $X\in\Ob\CCCP(\Al)$ is given by the complex
\[
\hspace*{-5mm}
\t X \= 
\Big(X^{0,0} 
\lrafl{25}{\smatez{d}{\dell}} 
X^{0,1}\ds X^{1,0} 
\mrafl{35}{\rsmatzd{d}{\dell}{0}{0}{-d}{-\dell}} 
X^{0,2}\ds X^{1,1}\ds X^{2,0} 
\mrafl{55}{
\enger{\left(\ba{rrrr}
\scm d  &\scm \dell & \scm  0     & \scm 0       \\
\scm 0  &\scm -d    & \scm -\dell & \scm 0       \\
\scm 0  &\scm  0    & \scm  d     & \scm\;\dell  \\
\ea\right)
}}
X^{0,3}\ds X^{1,2}\ds X^{2,1}\ds X^{3,0} 
\lra\cdots
\Big)\ru{10}
\]
in $\Ob\CC^{[0}(\Al)$. Using the induced morphisms, we obtain a total complex functor $\CCCP(\Al)\lraa{\t} \CC^{[0}(\Al)$. Since $\t$ maps elementary horizontally or vertically split acyclic double 
complexes to split acyclic complexes, it induces a functor $\KKKP(\Al)\lraa{\t} \KK^{[0}(\Al)$. If, in addition, $\Al$ is abelian, the total complex functor maps rowwise 
quasiisomorphisms and column\-wise quasiisomorphisms to quasiisomorphisms, as one sees using the long exact homology sequence and induction on a suitable filtration.

Note that we have an isomorphism $U\lraiso \t\Conc_1 U$, natural in $U\in\Ob\CC^{[0}(\Al)$, having entries $1_{U_0}$, $1_{U_1}$, $-1_{U_2}$, $-1_{U_3}$, $1_{U_4}$, etc. Moreover,
$U = \t\Conc_2 U$, natural in $U\in\Ob\CC^{[0}(\Al)$.

\subsubsection{The homotopy category of first quadrant double complexes as a quotient}

\begin{Lemma}
\label{LemDTC1}
The residue class functor $\CCC(\Al)\lra\KKK(\Al)$, restricted to $\CCCP(\Al)\lra\KKKP(\Al)$, induces an equivalence
\[
\CCCP(\Al)/\big(\CCC_\text{\rm sp\,ac}(\Al)\cap\CCCP(\Al)\big)\;\;\lraiso\;\; \KKKP(\Al)\; .
\]
\end{Lemma}

{\it Proof.} We have to show faithfulness; i.e.\ that if a morphism $X\lra Y$ in $\CCCP(\Al)$ factors over a split acyclic double complex, then it factors over a split acyclic double complex that 
lies in $\Ob\CCCP(\Al)$. By symmetry and additivity, it suffices to show that if a morphism $X\lra Y$ in $\CCCP(\Al)$ factors over a horizontally split acyclic double complex, then it factors over a 
horizontally split acyclic double complex that lies in $\Ob\CCCP(\Al)$. Furthermore, we may assume $X\lra Y$ to factor over an elementary horizontally split acyclic double complex $S$ concentrated
in the columns $k$ and $k+1$ for some $k\in\Z$. We may assume that $S^{i,j} = 0$ for $i < 0$ and $j\in\Z$. If $k < 0$, and in particular, if $k = -1$, then $X\lra Y$ is zero because $S\lra Y$ is
zero, so that in this case we may assume $S = 0$. On the other hand, if $k\ge 0$, then $S\in\Ob\CCCP(\Al)$.\qed

\bq
 Cf.\ also the similar Remark \ref{RemCEqis4_5}.

\eq

\subsection{Triple complexes}

\subsubsection{Definition}

Let $\CCCC(\Al) := \CC(\CC(\CC(\Al)))$ be the category of {\it triple complexes.} A triple complex $Y$ has entries $Y^{k,\ell,m}$ for $k,\,\ell,\, m\,\in\,\Z$. 

We denote the differentials in the three directions by $Y^{k,\ell,m}\lraa{d_1} Y^{k+1,\ell,m}$, $Y^{k,\ell,m}\lraa{d_2} Y^{k,\ell+1,m}$ and
$Y^{k,\ell,m}\lraa{d_3} Y^{k,\ell,m+1}$, respectively.

Let $k,\,\ell,\, m\,\in\,\Z$. We shall use the notation $Y^{-,\ell,=}$ for the double complex having at position $(k,m)$ the entry $Y^{k,\ell,m}$, differentials taken accordingly. Similarly the
complex $Y^{k,\ell,\ast}$ etc. 

Given a triple complex $Y\in\Ob\CCCC(\Al)$, we write $H Y^{-,=,\ast}\in\Ob\CCC(\Al)$ for the double complex having at position $(k,\ell)$ the entry $H(Y^{k,\ell,\ast})$, differentials
taken accordingly. 

Denote by $\CCCCP(\Al)\tm\CCCC(\Al)$ the full subcategory of {\it first octant triple complexes;} i.e.\ triple complexes $Y$ having $Y^{k,\ell,m} = 0$ whenever $k < 0$ or $\ell < 0$ or $m < 0$.

\subsubsection{Planewise total complex}
\label{SecPlanewise}

For $Y\in\Ob\CCCCP(\Al)$ we denote by $\t_{1,2} Y\in\Ob\CCCP(\Al)$ the {\it planewise total complex} of $Y$, defined for $m\in\Z$ as
\[
(\t_{1,2} Y)^{\ast,m} \; :=\; \t (Y^{-,=,m})\; ,
\]
with the differentials of $\t_{1,2} Y$ in the horizontal direction being induced by the differentials in the third index direction of $Y$, and with the differentials of $\t_{1,2} Y$ in the 
vertical direction being given by the total complex differentials. Explicitly, given $k,\,\ell\,\ge\, 0$, we have
\[
(\t_{1,2} Y)^{k,\ell} = \Ds_{i,\, j\,\ge\, 0,\; i+j\, =\, k} Y^{i,j,\ell}\; .
\]
By means of induced morphisms, this furnishes a functor
\[
\barcl
\CCCCP(\Al) & \lraa{\t_{1,2}} & \CCCP(\Al) \\
Y           & \lramaps        & \t_{1,2} Y \; .\\
\ea
\]

\section{Cartan-Eilenberg resolutions}
\label{SecCE}

\bq
 We shall use {\sc Quillen}'s language of exact categories~\bfcite{Qu73}{p.~15} to deal with Cartan-Eilenberg resolutions~\bfcite{CE56}{XVII.\S 1}, as it has been done by {\sc Mac Lane} already before 
 this language was available; cf.\ \mb{\bfcite{ML75}{XII.\S 11}}. The assertions in this section are for the most part wellknown.

\eq

\subsection{A remark}

\begin{Remark}
\label{RemCEqis4_5}
Let $\Al$ be an additive category. Then $\CC^{[0}(\Al)/\big(\CC^{[0}(\Al)\cap\CC_\text{\rm sp\,ac}(\Al)\big)\lra \KK^{[0}(\Al)$ is an equivalence.
\end{Remark}

{\it Proof.} Faithfulness is to be shown. A morphism $X\lra Y$ in $\CC^{[0}(\Al)$ that factors over an elementary split acyclic complex of the form $(\cdots\lra 0\lra T\ident T\lra0\lra\cdots$) with $T$ 
in positions $k$ and $k+1$ is zero, provided $k < 0$.\qed

\subsection{Exact categories}
\label{SecEx}

Concerning the terminology of exact categories, introduced by {\sc Quillen}~\bfcite{Qu73}{p.~15}, we refer to \mb{\bfcite{Ku05}{Sec.\ A.2}}. 

Let $\El$ be an exact category in which all idempotents split. An object $I\in\Ob\El$ is called {\it relatively injective,} or a {\it relative injective} (relative to the set of pure short 
exact sequences, that is), if $\liu{\El}{(-,I)}$ maps pure short exact sequences of $\El$ to short exact sequences. We say that $\El$ has {\it enough relative injectives,} if for all $X\in\Ob\El$, 
there exists a relative injective $I$ and a pure monomorphism $X\lramono I$.

In case $\El$ is an abelian category, with all short exact sequences stipulated to be pure, then we omit ``relative'' and speak of ``injectives'' etc.

\begin{Definition}
\label{DefEx2}\rm
Suppose given a complex $X\in\Ob\CC^{[0}(\El)$ with pure differentials. A {\it relatively injective complex resolution} of $X$ is a complex $I\in\Ob\CC^{[0}(\El)$, together with a quasiisomorphism
$X\lra I$, such that the following properties are satisfied.
\begin{itemize}
\item[(1)] The object entries of $I$ are relatively injective.
\item[(2)] The differentials of $I$ are pure.
\item[(3)] The quasiisomorphism $X\lra I$ consists of pure monomorphisms.
\end{itemize}
We often refer to such a relatively injective complex resolution just by $I$.

A {\it relatively injective object resolution}, or just a {\it relatively injective resolution,} of an object $Y\in\Ob\El$ is a relatively injective complex resolution of $\Conc Y$.

A {\it relatively injective resolution} is the complex of a relatively injective object resolution of some object in $\El$.
\end{Definition}

\begin{Remark}
\label{RemEx1}
Suppose that $\El$ has enough relative injectives. Every complex $X\in\Ob\CC^{[0}(\El)$ with pure differentials has a relatively injective complex resolution $I\in\Ob\CC^{[0}(\El)$.

In particular, every object $Y\in\Ob\El$ has a relatively injective resolution $J\in\Ob\CC^{[0}(\El)$.
\end{Remark}

{\it Proof.} Let $X^0\lramono I^0$ be a pure monomorphism into a relatively injective object $I^0$. Forming a pushout along $X^0\lramono I^0$, we obtain a pointwise purely monomorphic morphism 
of complexes $X\lra X'$ with $X'^0 = I^0$ and $X'^k = X^k$ for $k\ge 2$. By considering its cokernel, we see that it is a quasiisomorphism. So we may assume $X^0$ to be relatively injective. 

Let $X^1\lramono I^1$ be a pure monomorphism into a relatively injective object $I^1$. Form a pushout along $X^1\lramono I^1$ etc.
\qed

\begin{Remark}
\label{RemEx3}
Suppose given $X\in\Ob\CC^{[0}(\El)$ with pure differentials such that $\HH^k X \iso 0$ for $k\ge 1$. Suppose given $I\in\Ob\CC^{[0}(\El)$ such that $I^k$ is purely injective for $k\ge 0$, and such that
the differential $I^0\lraa{d} I^1$ has a kernel in $\El$. Then the map
\[
\liu{\KK^{[0}(\El)}{(X,I)} \;\lra\; \liu{\El}{\big(\Kern(X^0\lraa{d} X^1),\, \Kern(I^0\lraa{d} I^1)\big)} \\
\]
that sends a representing morphism of complexes to the morphism induced on the mentioned kernels, is bijective.
\end{Remark}

Suppose $\El$ to have enough relative injectives. Let $\Il\tm\El$ denote the full subcategory of relative injectives.
Let $\CC^{[0,\,\text{\rm res}}(\Il)$ denote the full subcategory of $\CC^{[0}(\Il)$ consisting of complexes $X$ with pure differentials such that $\HH^k X \iso 0$ for $k\ge 1$.
Let $\KK^{[0,\,\text{\rm res}}(\Il)$ denote the image of $\CC^{[0,\,\text{\rm res}}(\Il)$ in $\KK(\El)$.

\begin{Remark}
\label{RemEx4}
The functor $\CC^{[0,\,\text{\rm res}}(\Il)\lra \El$, $X\lramaps\HH^0(X)$, induces an equivalence 
\[
\KK^{[0,\,\text{\rm res}}(\Il)\;\lraiso\; \El \; .
\]
\end{Remark}

{\it Proof.} This functor is dense by Remark~\ref{RemEx1}, and full and faithful by Remark~\ref{RemEx3}.\qed

\begin{Remark}[exact Horseshoe Lemma]
\label{RemEx5}
\ \\
Given a pure short exact sequence $X'\lra X\lra X''$ and relatively injective resolutions $I'$ of $X'$ and $I''$ of $X''$, there exists a relatively injective resolution $I$ of $X$ and
a pointwise split short exact sequence $I'\lra I\lra I''$ that maps under $\HH^0$ to $X'\lra X\lra X''$.
\end{Remark}

{\it Proof.} Choose pure monomorphisms $X'\lramono I'^0$ and $X''\lramono I''^0$ into relative injectives $I'^0$ and $I''^0$. Embed them into a morphism from the pure short exact sequence
$X'\lramono X\lraepi X''$ to the split short exact sequence $I'\lrafl{25}{\smatez{1}{0}} I'\ds I''\lrafl{36}{\smatze{0}{1}} I''\ru{8}$. Insert the pushout $T$ of $X'\lramono X$ along $X'\lramono I'^0$ and
the pullback of $I'^0\ds I''^0\lraepi I''^0$ along $X''\lramono I''^0$ to see that $X\lra I'^0\ds I''^0$ is purely monomorphic. So we can take the cokernel 
$\BB^1 I'\lra \BB^1 I \lra \BB^1 I''$ of this morphism of pure short exact sequences. Considering the cokernels on the commutative triangle \mb{$(X,T,I'^0\ds I''^0)$} of pure monomorphisms, we obtain
a bicartesian square $(T,\,I'^0\ds I''^0,\, \BB^1 I',\,\BB^1 I)$ and conclude that the sequence of cokernels is itself purely short exact. So we can iterate.\qed

\subsection{An exact category structure on $\CC(\Al)$}
\label{SecCompl}

Let $\Al$ be an abelian category with enough injectives.

\begin{Remark}
\label{RemCEqis1}
The following conditions on a short exact sequence $X'\lra X\lra X''$ in $\CC(\Al)$ are equivalent.
\begin{itemize}
\item[{\rm (1)}\phantom{$'$}] All connectors in its long exact homology sequence are equal to zero.
\item[{\rm (2)}\phantom{$'$}] The sequence $\BB^k X'\lra\BB^k X\lra\BB^k X''$ is short exact for all $k\in\Z$.
\item[{\rm (3)}\phantom{$'$}] The morphism $\ZZ^k X\lra \ZZ^k X''$ is epimorphic for all $k\in\Z$.
\item[{\rm (3$'$)}] The morphism $\ZZ'^k X'\lra \ZZ'^k X$ is monomorphic for all $k\in\Z$.
\item[{\rm (4)}\phantom{$'$}] The diagram 
\[
\xymatrix{
\BB^k X'\ar[r]\ar[d]  & \ZZ^k X'\ar[r]\ar[d]  & \HH^k X'\ar[d]  \\
\BB^k X\ar[r]\ar[d]   & \ZZ^k X\ar[r]\ar[d]   & \HH^k X\ar[d]   \\
\BB^k X''\ar[r]       & \ZZ^k X''\ar[r]       & \HH^k X'' \\
}
\]
has short exact rows and short exact columns for all $k\in\Z$.
\end{itemize}
\end{Remark}

{\it Proof.} We consider the diagram in (4) as a (horizontal) short exact sequence of (vertical) complexes and regard its long exact homology sequence. 
Taking into account that all assertions are supposed to hold for all $k\in\Z$, we can employ the long exact homology sequence on $X'\lra X\lra X''$ to prove the equivalence of (1), (2), (3) and (4).

Now the assertion (1) $\equ$ (3) is dual to the assertion (1) $\equ$ (3$'$).\qed

\begin{Remark}
\label{RemCEqis2_5}
The category $\CC(\Al)$, equipped with the set of short exact sequences that have zero connectors on homology as pure short exact sequences, is an exact category with enough relatively injective 
objects in which all idempotents split. With respect to this exact category structure on $\CC(\Al)$, a complex is relatively injective if and only if it is split and has injective object entries.
\end{Remark}

Cf.\ \bfcite{ML75}{XII.\S 11}, where pure short exact sequences are called {\it proper}. A relatively injective object in $\CC(\Al)$ is also referred to as an {\it injectively split complex.} 
To a relatively injective resolution of a complex $X\in\Ob\CC(\Al)$, we also refer as a {\it Cartan-Eilenberg-resolution,} or, for short, as a {\it CE-resolution} of $X$; cf.\ \bfcite{CE56}{XVII.\S 1}. 
A {\it CE-resolution} is a CE-resolution of some complex. Considered as a double complex, it is in particular rowwise split and has injective object entries. 

Given a morphism $X\lraa{f} X'$ in $\CC(\Al)$, CE-resolutions $J$ of $X$ and $J'$ of $X'$, a morphism $J\lraa{\h f} J'\ru{5}$ in $\CCC(\Al)$ such that
$(J^{i,j}\lraa{\h f^{i,j}} J'^{i,j}) = (0\lra 0)\ru{5}$ for $i < 0$ and such that
\[
\HH^0(J^{\ast,-}\lraa{\h f^{\ast,-}} J'^{\ast,-})\;\=\; (X\lraa{f} X')
\]
is called a {\it CE-resolution} of $X\lraa{f} X'$. By Remarks \ref{RemCEqis2_5} and \ref{RemEx4}, each morphism in $\CC(\Al)$ has a CE\nobreakdash-resolution.

{\it Proof of Remark~{\rm\ref{RemCEqis2_5}}.} We {\it claim} that $\CC(\Al)$, equipped with the said set of short exact sequences, is an exact category. We verify the conditions (Ex 1,\,2,\,3) listed 
in \bfcite{Ku05}{Sec.\ A.2}. The conditions \mb{(Ex 1$^\0$,\,2$^\0$,\,3$^\0$)} then follow by duality.

Note that by Remark~\ref{RemCEqis1}.(3$'$), a monomorphism $X\lra Y$ in $\CC(\Al)$ is pure if and only if $\ZZ'^k(X\lra Y)$ is monomorphic in $\Al$ for all $k\in\Z$.

Ad (Ex 1). To see that a split monomorphism is pure, we may use additivity of the functor $\ZZ'^k$ for $k\in\Z$.

Ad (Ex 2). To see that the composition of two pure monomorphisms is pure, we may use $\ZZ'^k$ being a functor for $k\in\Z$.

Ad (Ex 3). Suppose given a commutative triangle
$$
\xymatrix{
                              & Y\ar~+{|*\dir{|}}[dr] & \\
X\ar[ur] \ar~+{|*\dir{*}}[rr] &                       & Z\; , \!\! \\
}
$$
in $\CC(\Al)$. Applying the functor $\ZZ'^k$ to it, for $k\in\Z$, we conclude that $\ZZ'^k(X\lra Y)$ is monomorphic, whence $X\lra Y$ is purely monomorphic. So we may complete to
$$
\xymatrix{
A\ar~+{|*\dir{*}}[dr] \ar[rr]              &                                           & B \\
                                           & Y\ar~+{|*\dir{|}}[dr]\ar~+{|*\dir{|}}[ur] & \\
X\ar~+{|*\dir{*}}[ur] \ar~+{|*\dir{*}}[rr] &                                           & Z \\
}
$$
in $\CC(\Al)$ with $(X,Y,B)$ and $(A,Y,Z)$ pure short exact sequences. Applying $\ZZ'^k$ to this diagram, we conclude that $\ZZ'^k(A\lra B)$ is a monomorphism for $k\in\Z$, whence $A\lra B$ is 
a pure monomorphism. 

This proves the {\it claim.}

Note that idempotents in $\CC(\Al)$ are split since $\CC(\Al)$ is also an abelian category.

We {\it claim\,} relative injectivity of complexes with split differentials and injective object entries. By a direct sum decomposition, and using the fact that any monomorphism from an 
elementary split acyclic complex with injective entries to an arbitrary complex is split, we are reduced to showing that a pure monomorphism from a complex with a single nonzero injective entry, 
at position $0$, say, to an arbitrary complex is split. So suppose given $I\in\Ob\Inj\Al$, $X\in\Ob\CC(\Al)$ and a pure monomorphism $\Conc I\lramono X$. Using 
Remark~\ref{RemCEqis1}.(3$'$), we may choose a retraction to the composite $(I\lra X^0\lra \ZZ'^0 X)$. This yields a retraction to $I\lra X^0$ that composes to $0$ with $X^{-1}\lra X^0$, which can be 
employed for the sought retraction $X\lra\Conc I$. This proves the {\it claim.}

Let $X\in\Ob\CC(\Al)$. We {\it claim\,} that there exists a pure monomorphism from $X$ to a relatively injective complex. Since $\Al$ has enough injectives, by a direct sum decomposition
we are reduced to finding a pure monomorphism from $X$ to a split complex. Consider the following morphism $\phi_k$ of complexes for $k\in\Z$,
\[
\xymatrix{
\cdots \ar[r] & 0 \ar[r]              & X^k\ar[r]^(0.32){\smatez{1}{0}} & X^k\ds \ZZ'^k X \ar[r]            & 0 \ar[r]             & \cdots        \\ 
\cdots \ar[r] & X^{k-2}\ar[r]^d\ar[u] & X^{k-1}\ar[r]^d\ar[u]^d         & X^k\ar[r]^d\ar[u]_{\smatez{1}{p}} & X^{k+1} \ar[r]\ar[u] & \cdots \zw{,} \\ 
}
\]
where $X^k\lraepia{p} \ZZ'^k X$ is taken from $X$. The functor $\ZZ'^k$ maps it to the identity. We take the direct sum of the upper complexes over $k\in\Z$ and let the morphisms $\phi_k$ be the 
components of a morphism $\phi$ from $X$ to this direct sum. At position $k$, this morphism $\phi$ is monomorphic because $\phi_k$ is. Moreover, $\ZZ'^k(\phi)$ is a monomorphism because $\ZZ'^k(\phi_k)$ 
is. Hence $\phi$ is purely monomorphic by condition (3$'$) of Remark \ref{RemCEqis1}. This proves the {\it claim.}\qed

\begin{Remark}
\label{RemCEqis3}
Write $\El := \CC(\Al)$. Given $\ell\ge 0$, we have a homology functor $\El\lraa{\HH^\ell}\Al$, which induces a functor $\CC(\El)\lrafl{27}{\CC(\HH^\ell)}\CC(\Al)$.
Suppose given a purely acyclic complex $X\in\Ob\CC(\El)$. Then $\CC(\HH^\ell) X\in\Ob\CC(\Al)$ is acyclic.
\end{Remark}

{\it Proof.} This follows using the definition of pure short exact sequences, i.e.\ Remark~\ref{RemCEqis1}.(1).\qed

\subsection{An exact category structure on $\CC^{[0}(\Al)$}
\label{SecCompl2}

Write $\CCCCE(\Inj\Al)$ for the full subcategory of $\CCCP(\Al)$ whose objects are CE-resolutions. Write $\KKKCE(\Inj\Al)$ for the full subcategory of $\KKKP(\Al)$ whose objects are CE-resolutions. 

\begin{Remark}
\label{RemCEqis3_5}
The category $\CC^{[0}(\Al)$, equipped with the short exact sequences that lie in $\CC^{[0}(\Al)$ and that are pure in $\CC(\Al)$ in the sense of Remark~{\rm\ref{RemCEqis2_5}} as pure short exact
sequences, is an exact category wherein idempotents are split. It has enough relative injectives, viz.\ injectively split complexes that lie in $\CC^{[0}(\Al)$.
\end{Remark}

{\it Proof.} To show that it has enough relative injectives, we replace $\phi_0$ in the proof of Remark~\ref{RemCEqis2_5} by $X\lraa{\phi'_0}\Conc X^0$, defined by $X_0\lrafl{30}{1_{X_0}} X_0$ at
position $0$.\qed

\subsection{The Cartan-Eilenberg resolution of a quasiisomorphism}
\label{SecCEqis}
 
Abbreviate $\El := \CC(\Al)$, which is an exact category as in Remark~\ref{RemCEqis2_5}. Consider $\CCCP(\Al) \tm \CC^{[0}(\El)$, where the second index of $X\in\Ob\CCCP(\Al)$ counts the positions 
in $\El = \CC(\Al)$; i.e.\ when $X$ is viewed as a complex with values in $\El$, its entry at position $k$ is given by $X^{k,\ast}\in\El = \CC(\Al)$.

\begin{Remark}
\label{RemCeqis7_5}
Suppose given a split acyclic complex $X\in\Ob\CC^{[0}(\Al)$. There exists a horizontally split acyclic CE-resolution $J\in\Ob\CCCCE(\Inj\Al)$ of $X$. 
\end{Remark}

{\it Proof.} This holds for an elementary split acyclic complex, and thus also in the general case by taking a direct sum.\qed

\begin{Lemma}
\label{LemCEnew1}
Suppose given $X\in\Ob\CCCP(\Al)$ with pure differentials when considered as an object of $\CC^{[0}(\El)$, and with $\HH^k\big( X^{\ast,-}\big) \iso 0$ in $\CC^{[0}(\Al)$ for $k \ge 1$.

Suppose given $J\in\Ob\CCCP(\Inj\Al)$ with split rows $J^{k,\ast}$ for $k\ge 1$. In other words, $J$ is supposed to consist of relative injective object entries when considered 
as an object of $\CC^{[0}(\El)$.

Then the map
$$
\liu{\KKKP(\Al)}{(X,J)}\;\;\mrafl{30}{\HH^0\left((-)^{\ast,-}\right)}\;\; \liu{\KK^{[0}(\Al)}{\left(\HH^0\big(X^{\ast,-}\big),\,\HH^0\big(J^{\ast,-}\big)\right)}
\leqno (\ast)
$$
is bijective.
\end{Lemma}

{\it Proof.} First, we observe that by Remark \ref{RemEx3}, we have
$$
\liu{\KK^{[0}(\El)}{(X,J)}\;\;\mraisofl{30}{\HH^0\left((-)^{\ast,-}\right)}\;\; \liu{\El}{\left(\HH^0\big(X^{\ast,-}\big),\,\HH^0\big(J^{\ast,-}\big)\right)}\; . \ru{5}
\leqno (\ast\ast)
$$
So it remains to show that $(\ast)$ is injective.
Let $X\lraa{f} J$ be a morphism that vanishes under $(\ast)$. Then $\HH^0\big(X^{\ast,-}\big)\lra \HH^0\big(J^{\ast,-}\big)$ factors over a split acyclic complex $S\in\Ob\CC^{[0}(\Al)$; 
cf.\ Remark~\ref{RemCEqis4_5}. Let $K$ be a horizontally split acyclic CE-resolution of $S$; cf.\ Remark~\ref{RemCeqis7_5}. By Remark~\ref{RemEx3}, we obtain a morphism $X\lra K$ that lifts 
$\HH^0\big(X^{\ast,-}\big)\lra S$ and a morphism $K\lra J$ that lifts $S\lra \HH^0\big(J^{\ast,-}\big)$. The composite $X\lra K\lra J$ vanishes in $\KKKP(\Al)$. The difference
\[
(X\lraa{f} J) - (X\lra K\lra J)
\]
lifts $\HH^0\big(X^{\ast,-}\big)\lraa{0} \HH^0\big(J^{\ast,-}\big)$. Hence by $(\ast\ast)$, it vanishes in $\KK^{[0}(\El)$ and so a fortiori in $\KKKP(\Al)$. Altogether,
$X\lraa{f} J$ vanishes in $\KKKP(\Al)$.\qed

\begin{Proposition}
\label{PropCEqis6} 
The functor $\CCCCE(\Inj\Al)\mrafl{28}{\HH^0\left((-)^{\ast,-}\right)}\CC^{[0}(\Al)$ induces an equivalence
\[
\KKKCE(\Inj\Al)\;\;\mraisofl{28}{\HH^0\left((-)^{\ast,-}\right)}\;\;\KK^{[0}(\Al)\; .
\]
\end{Proposition}

{\it Proof.} By Lemma \ref{LemCEnew1}, this functor is full and faithful. By Remark \ref{RemEx1}, it is dense.\qed

\begin{Corollary}
\label{CorCeqis7}
Suppose given $X,\, X'\,\in\,\Ob\Cl^{[0}(\Al)$. Let $J$ be a CE-resolution of $X$. Let $J'$ be a CE-resolution of $X'$. If $X$ and $X'$ are isomorphic in $\KK^{[0}(\Al)$, then
$J$ and $J'$ are isomorphic in $\KKKP(\Al)$.
\end{Corollary}

\bq
 The following lemma is to be compared to Remark \ref{RemCeqis7_5}.

\eq

\begin{Lemma}
\label{LemCeqis8}
Suppose given an acyclic complex $X\in\Ob\CC^{[0}(\Al)$. There exists a rowwise split acyclic CE-resolution $J$ of $X$. Each CE-resolution of $X$ is isomorphic to $J$ in $\KKKP(\Al)$.
\end{Lemma}

{\it Proof.} By Corollary \ref{CorCeqis7}, it suffices to show that there exists a rowwise split acyclic \mb{CE-resolution} of $X$. Recall that a CE-resolution of an arbitrary complex 
$Y\in\Ob\CC^{[0}(\Al)$ can be constructed by a choice of injective resolutions of $\HH^k Y$ and $\BB^k Y$ for $k\in\Z$, followed by an application of the abelian Horseshoe Lemma to the 
short exact sequences $\BB^k Y\lra \ZZ^k Y\lra \HH^k Y$ for $k\in\Z$ and then to $\ZZ^k Y\lra Y^k\lra\BB^{k+1} Y$ for $k\in\Z$; cf.\ \bfcite{CE56}{Chap.~XVII, Prop.~1.2}. Since 
$\HH^k X = 0$ for $k\in\Z$, we may choose the zero resolution for it. Applying this construction, we obtain a rowwise split acyclic CE-resolution.
\qed 

Given $X\lraa{f} X'$ in $\CC^{[0}(\Al)$, a morphism $J\lraa{\h f} J'$ in $\CCCP(\Al)$ is called a {\it CE-resolution} of $X\lraa{f} X'$ if $\HH^0(\h f^{\ast,-}) \iso f$, as diagrams of the form
$\bt\lra\bt$. By Remark~\ref{RemEx3}, given CE-resolutions $J$ of $X$ and $J'$ of $X'$, there exists a CE-resolution $J\lraa{\h f} J'$ of $X\lraa{f} X'$.

\begin{Proposition} \hspace*{1mm} 
\label{PropCEqis9}
Let $X\lraa{f} X'$ be a quasiisomorphism in $\CC^{[0}(\Al)$. Let $J\lraa{\h f} J'$ be a \mb{CE-resolution} of $X\lraa{f} X'$. Then $\h f$ can be written as a composite in $\CCCCE(\Inj\Al)$ of a 
rowwise homotopism, followed by a double homotopism.
\end{Proposition}

{\it Proof.} Choose a pointwise split monomorphism $X\lraa{a} A$ into a split acyclic complex $X$. We can factor 
\[
(X\lraa{f} X') \; \=\; \left(X\lrafl{28}{\smatez{f}{a}} X'\ds A\lrafl{35}{\smatze{1}{0}} X'\right)\; ,
\]
so that $\smatez{f}{a}$ is a pointwise split monomorphism. Let $B$ be a CE-resolution of $A$. Choosing a CE-resolution $b$ of $a$, we obtain the factorisation
\[
(J\lraa{\h f} J') \; \=\; \left(J\lrafl{28}{\smatez{\h f}{b}} J'\ds B\lrafl{35}{\smatze{1}{0}} J'\right)\; .
\]
Since $X'\ds A\lrafl{35}{\smatze{1}{0}} X'$ is a homotopism, $J'\ds B\lrafl{35}{\smatze{1}{0}} J$ is a double homotopism; cf.\ Corollary~\ref{CorCeqis7}. Hence $\h f$ is a composite of a
rowwise homotopism and a double homotopism if and only if this holds for $\smatez{\h f}{b}$. So we may assume that $f$ is pointwise split monomorphic, so in particular, monomorphic.

By Proposition~\ref{PropCEqis6}, we may replace the given CE-resolution $\h f$ by an arbitrary CE-resolution of $f$ between $J$ and an arbitrarily chosen CE-resolution of $X'$ without changing
the property of being a composite of a rowwise homotopism and a double homotopism for this newly chosen CE-resolution of $f$.

Let $X\lraa{f} X'\lra \b X$ be a short exact sequence in $\CC^{[0}(\Al)$. Since $f$ is a quasiisomorphism, $\b X\in\Ob\CC^{[0}(\Al)$ is acyclic. Let $\b J$ be a rowwise split acyclic 
CE\nobreakdash-resolution of $\b X$; cf.\ Lemma~\ref{LemCeqis8}. The short exact sequence $X\lraa{f} X'\lra\b X$ is pure by acyclicity of $\b X$; cf.\ Remark~\ref{RemCEqis1}.(1). Hence by the exact 
Horseshoe Lemma, there exists a rowwise split short exact sequence $J\lra \w J'\lra \b J$ of CE-resolutions that maps to $X\lraa{f} X'\lra \b X$ under $\HH^0\big((-)^{\ast,-}\big)$; cf.\ 
Remark~\ref{RemEx5}. Since $\b J$ is rowwise split acyclic and since the sequence $J\lra \w J'\lra \b J$ is rowwise split short exact, $J\lra \w J'$ is a rowwise homotopism. Since $J\lra \w J'$ is a 
CE-resolution of $X\lraa{f} X'$, this proves the proposition.\qed 

\section{Formalism of spectral sequences}
\label{SecForm}

\bq
 We follow essentially {\sc Verdier} \bfcite{Ve67}{II.4}; cf.\ \bfcite{De94}{App.}; on a more basic level, cf.\ \bfcite{Ku06}{Kap.\ 4}. 

\eq

Let $\Al$ be an abelian category.

\subsection{Pointwise split and pointwise finitely filtered complexes}

Let $\Zi := \{-\infty\}\disj\Z\disj\{\infty\}$, considered as a linearly ordered set, and thus as a category. Write $]\al,\be] := \{\sa\in\Zi\; :\; \al < \sa \le \be\}$ for $\al,\,\be\,\in\,\Zi$
such that $\al\le\be$; etc.

Given $X\in\Ob\bo\Zi,\CC(\Al)\bc$, the morphism of $X$ on $\al\le\be$ in $\Zi$ shall be denoted by $X(\al)\lraa{x} X(\be)$.

An object $X\in\Ob\bo\Zi,\CC(\Al)\bc$ is called a {\it pointwise split and pointwise finitely filtered complex (with values in $\Al$),} provided (SFF 1,\,2,\,3) hold.

\begin{itemize}
\item[(SFF 1)] We have $X(-\infty) = 0$. 
\item[(SFF 2)] The morphism $X(\al)^i\lraa{x^i} X(\be)^i$ is split monomorphic for all $i\in\Z$ and all $\al\le\be$ in $\Zi$.
\item[(SFF 3)] For all $i\in\Z$, there exist $\be_0,\,\al_0\,\in\,\Z$ such that $X(\al)^i\lraa{x^i} X(\be)^i$ is an identity whenever $\al\le\be\le\be_0$ or $\al_0\le\al\le\be$ in $\Zi$.
\end{itemize}

The pointwise split and pointwise finitely filtered complexes with values in $\Al$ form a full subcategory $\text{SFFC}(\Al)\tm\bo\Zi,\CC(\Al)\bc$.

Suppose given a pointwise split and pointwise finitely filtered complex $X$ with values in $\Al$ for the rest of the present \S\ref{SecForm}. 

Let $\al\in\Zi$. Write $\b X(\al) := \mb{$\Cokern\big(X(\al - 1)\lra X(\al)\big)$}$ for $\al\in\Z$. Given $i\in\Z$, we obtain $X(\al)^i \iso \Ds_{\sa\in ]-\infty,\al]}\b X(\sa)^i\ru{-2.5}$, which is a 
finite direct sum. We identify along this isomorphism. In particular, we get as a matrix representation for the differential
\[
\big(X(\al)^i\lraa{d} X(\al)^{i+1}\big) \; \=\; \left(\Ds_{\sa\in ]-\infty,\al]}\b X(\sa)^i \;\;\mrafl{30}{(d^i_{\sa,\ta})_{\sa,\ta}}\;\; \Ds_{\ta\in ]-\infty,\al]}\b X(\ta)^{i+1} \right)\; ,
\]
where $d^i_{\sa,\ta} = 0$ whenever $\sa < \ta$; a kind of lower triangular matrix.

\subsection{Spectral objects}

Let $\bZi := \Zi\ti\Z$. Write $\al^{+k} := (\al,k)$, where $\al\in\Zi$ and $k\in\Z$. Let $\al^{+k} \le \be^{+\ell}$ in $\bZi$ if $k < \ell$ or ($k = \ell$ and $\al\le\be$), i.e.\ let $\bZi$
be linearly ordered via a lexicographical ordering. We have an automorphism $\al^{+k}\lramaps\al^{+k+1}$ of the poset $\bZi$, to which we refer as {\it shift.} Note that $-\infty^{+k} = (-\infty)^{+k}$.

We have an order preserving injection $\Zi\lra\bZi\,$, $\al\lramaps\al^{+0}$. We use this injection as an identification of $\Zi$ with its image in $\bZi\,$, i.e.\ we sometimes write $\al := \al^{+0}$
by abuse of notation.

Let $\bZif := \{ (\al,\be) \in \bZi\ti\bZi\; :\; \be^{-1}\le\al\le\be\le\al^{+1}\}$. We usually write $\be/\al := (\al,\be)\in\bZif$; reminiscent of a quotient. The set $\bZif$ is partially ordered
by $\be/\al\le\be'/\al'$ $:\equ$ ($\be\le\be'$ and $\al\le\al'$). We have an automorphism $\be/\al\lramaps (\be/\al)^{+1} := \al^{+1}/\be$ of the poset $\bZif$, to which, again, we refer as {\it shift}. 

We write $\Zif := \{ \be/\al \in\bZif \; :\; -\infty\le\al\le\be\le\infty\}$. Note that any element of $\bZif$ can uniquely be written as $(\be/\al)^{+k}$ for some $\be/\al\in\Zif$ and
some $k\in\Z$.

We shall construct the {\it spectral object} $\text{Sp}(X)\in\Ob\bo\bZif,\KK(\Al)\bc$. The morphism of 
$\text{Sp}(X)$ on $\be/\al\le\be'/\al'$ in $\bZif$ shall be denoted by $X(\be/\al)\lraa{x} X(\be'/\al')$.

We require that 
\[
\left(X\big((\be/\al)^{+k}\big)\;\lraa{x}\; X\big((\be'/\al')^{+k}\big)\right) \;\; =\;\; \left(X(\be/\al)\;\lraa{x}\; X(\be'/\al')\right)^{\bt + k}
\]
for $\be/\al\le\be'/\al'$ in $\bZif$; i.e., roughly put, that $\text{Sp}(X)$ be periodic up to shift of complexes.

Define 
\[
X\big(\be/\al\big) \;:=\; \Cokern\big(X(\al)\lraa{x} X(\be)\big) 
\]
for $\be/\al\in\Zif$. By periodicity, we conclude that $X\big(\al/\al\big) = 0$ and $X\big(\al^{+1}/\al\big) = 0$ for all $\al\in\bZi$.

Write 
\[
D^i_{\be/\al,\,\be'/\al'} \; := \; (d^i_{\sa,\ta})_{\sa\in ]\al,\be],\; \ta\in ]\al',\be']} \; :\; X(\be/\al)^i\;\lra\; X(\be'/\al')^{i+1}
\]
for $i\in\Z$ and \mb{$\be/\al,\;\be'/\al'\,\in\,\Zif$}. 

Given $-\infty\le\al\le\be\le\ga\le\infty$ and $i\in\Z$, we let
\[
\barcl
\bigg(X(\be/\al)^i \;\lraa{x^i}\; X(\ga/\al)^i\bigg)      & := & \bigg(X(\be/\al)^i \;\lrafl{25}{\smatez{1}{0}}\;  X(\be/\al)^i\ds X(\ga/\be)^i\bigg) \vspace*{1mm}\\
\bigg(X(\ga/\al)^i \;\lraa{x^i}\; X(\ga/\be)^i\bigg)      & := & \bigg(X(\be/\al)^i\ds X(\ga/\be)^i \;\lrafl{35}{\smatze{0}{1}}\; X(\ga/\be)^i \bigg) \vspace*{1mm}\\
\bigg(X(\ga/\be)^i \;\lraa{x^i}\; X(\al^{+1}/\be)^i\bigg) & := & \bigg(X(\ga/\be)^i\;\mrafl{35}{D_{\ga/\be,\, \be/\al}^i}\; X(\be/\al)^{i+1}\bigg)\; . \\
\ea
\]
By periodicity up to shift of complexes, this defines $\text{Sp}(X)$. The construction is functorial in $X\in\Ob\text{SFFC}(\Al)$.

\subsection{Spectral sequences}

Let $\bZiff := \{(\ga/\al,\de/\be)\in\bZif\ti\bZif\; :\; \de^{-1}\le\al\le\be\le\ga\le\de\le\al^{+1}\}$. Given $(\ga/\al,\de/\be)\in\bZiff$, we usually write $\de/\be\bby\ga/\al := (\ga/\al,\de/\be)$.
The set $\bZiff$ is partially ordered by 
\[
\de/\be\bby\ga/\al \le \de'/\be'\bby\ga'/\al' \;\; :\equ \;\; \text{($\ga/\al\le\ga'/\al'$ and $\de/\be\le\de'/\be'$)}\; . 
\]
Define the {\it spectral sequence} $\EE(X)\in\Ob\bo\bZiff,\Al\bc$ of $X$ by letting its value on 
\[
\de/\be\bby\ga/\al \;\le\; \de'/\be'\bby\ga'/\al' 
\]
in $\bZiff$ be the morphism that appears in the middle column of the diagram
\[
\xymatrix{
\HH^0\big(X(\ga/\al)\big)\ar~+{|(0.44)*\dir{|}}[r]\ar[d]_{\HH^0(x)} & \EE(\de/\be\bby\ga/\al)(X)\ar~+{|(0.55)*\dir{*}}[r]\ar[d]_e & \HH^0\big(X(\de/\be)\big)\ar[d]_{\HH^0(x)} \\
\HH^0\big(X(\ga'/\al')\big)\ar~+{|(0.44)*\dir{|}}[r]                & \EE(\de'/\be'\bby\ga'/\al')(X)\ar~+{|(0.55)*\dir{*}}[r]     & \HH^0\big(X(\de'/\be')\big)\; . \!\!       \\
}
\]
Given $\de/\be\bby\ga/\al\in\bZiff$ and $k\in\Z$, we also write 
\[
\EE(\de/\be\bby\ga/\al)^{+k}(X) \;:=\; \EE\big((\de/\be)^{+k}\bby(\ga/\al)^{+k}\big)(X) \; .
\]

Altogether,
\[
\ba{rcrclcl}
\bo\Zi,\CC(\Al)\bc & \om & \text{SFFC}(\Al) & \lra     & \bo\bZif,\KK(\Al)\bc & \lra     & \bo\bZiff,\Al\bc \\
                   &     & X                & \lramaps & \text{Sp}(X)         & \lramaps & \EE(X)\; . \\
\ea
\]

\subsection{A short exact sequence}
\label{SecSES}

\begin{Lemma}
\label{LemSES1}
Given $\eps^{-1}\le\al\le\be\le\ga\le\de\le\eps\le\al^{+1}$ in $\bZi$, we have a short exact sequence
\[
\EE(\eps/\be\bby\ga/\al)(X) \;\;\lramonoa{e}\;\; \EE(\eps/\be\bby\de/\al)(X) \;\;\lraepia{e}\;\; \EE(\eps/\ga\bby\de/\al)(X)\; .
\]
\end{Lemma}

{\it Proof.} See \bfcite{Ku05}{Lem.\ 3.9}.\qed

\begin{Lemma}
\label{CorSES2}
Given $\eps^{-1}\le\al\le\be\le\ga\le\de\le\eps\le\al^{+1}$ in $\bZi$, we have a short exact sequence
\[
\EE(\eps/\ga\bby\de/\al)(X) \;\;\lramonoa{e}\;\; \EE(\eps/\ga\bby\de/\be)(X) \;\;\lraepia{e}\;\; \EE(\al^{+1}/\ga\bby\de/\be)(X)\; .
\]
\end{Lemma}

{\it Proof.} Apply the functor induced by $\be/\al\lramaps \al^{+1}/\be$ to $\text{Sp}(X)$. Then apply \bfcite{Ku05}{Lem.\ 3.9}.\qed

The short exact sequence in Lemma \ref{LemSES1} is called a {\it fundamental short exact sequence (in first notation)}, the short exact sequence in Lemma \ref{CorSES2} is called 
a {\it fundamental short exact sequence (in second notation)}. They will be used without further comment.

\subsection{Classical indexing}
\label{SecClassical}

Let $1\le r\le \infty$ and let $p,\,q\,\in\,\Z$. Denote 
\[
\EE_r^{p,q} \= \EE_r^{p,q}(X) \; := \; \EE(-p-1+r/\!-\!p-1\bby\!-\!p/\!-\!p-r)^{+p+q}(X) \; ,
\]
where $i + \infty := \infty$ and $i - \infty := -\infty$ for all $i\in\Z$.

\bq
 \begin{Example}
 \label{ExSES2_5}\rm
 The short exact sequences in Lemmata \ref{LemSES1}, \ref{CorSES2} allow to derive the exact couples of Massey. Write 
 $\DD_r^{i,j} = \DD_r^{i,j}(X) := \EE(-i/\!-\!\infty\fbby\!-\!i\!-\! r\!+\! 1/\!-\!\infty)^{+i+j}(X)$ for $i,\, j\,\in\,\Z$ and $r\ge 1$. We obtain an exact sequence
 \[
 \DD_r^{i,\,j}\;\lraa{e}\;\DD_r^{i-1,\,j+1}\;\lraa{e}\;\EE_r^{i+r-2,\,j-r+2}\;\lraa{e}\;\DD_r^{i+r-1,\,j-r+2}\;\lraa{e}\;\DD_r^{i+r-2,\,j-r+3}
 \]
 by Lemmata \ref{LemSES1}, \ref{CorSES2}. 
 \end{Example}

\eq

\subsection{Comparing proper spectral sequences} 
\label{SecCompProp}

Let $X\lraa{f} Y$ be a morphism in $\text{SFFC}(\Al)$, i.e.\ a morphism of pointwise split and pointwise finitely filtered complexes with values in $\Al$. Write $\EE(X)\lrafl{27}{\EE(f)} \EE(Y)$ 
for the induced morphism on the spectral sequences.

For $\al,\,\be\,\in\,\bZi$, we write $\al\dl\be$ if 
\[
\big(\text{$\al < \be$}\big) \hspace*{5mm} \text{or} \hspace*{5mm} \big(\text{$\al = \be\;\;$ and $\;\;\al\in\{\infty^{+k}\; :\; k\in\Z\}\cup \{-\infty^{+k}\; :\; k\in\Z\}$}\big) \; .
\]
We write
\[
\bZiffp \; :=\; \{ \de/\be\bby\ga/\al\in\bZiff\; :\; \de^{-1}\le \al \dl \be \le \ga \dl \de \le \al^{+1}\} \; .
\]
We write
\[
\EEp \= \EEp(X) \; :=\; \EE(X)|_{\bZiffp} \;\in\; \Ob\bo\bZiffp,\Al\bc
\]
for the {\it proper spectral sequence} of $X$; analogously for the morphisms.

\begin{Lemma}
\label{LemProper1}
If $\EE(\al+1/\al-1\bby\al/\al-2)^{+k}(f)$ is an isomorphism for all $\al\in\Z$ and all $k\in\Z$, then $\EEp(f)$ is an isomorphism.
\end{Lemma}

{\it Proof.} {\it Claim 1.} We have an isomorphism $\EE(\ga/\be-1\bby\be/\be-2)^{+k}(f)$ for all $k\in\Z$, all $\be\in\Z$ and all $\ga\in\Z$ such that $\ga > \be$. We have an isomorphism
$\EE(\be+1/\be-1\bby\be/\al-1)^{+k}(f)$ for all $k\in\Z$, all $\be\in\Z$ and all $\al\in\Z$ such that $\al < \be$.

The assertions follow by induction using the exact sequences
\[
\EE(\ga+2/\ga\bby\ga+1/\be)^{+k-1}\;\lraa{e}\; \EE(\ga/\be-1\bby\be/\be-2)^{+k}\;\lraa{e}\; \EE(\ga+1/\be-1\bby\be/\be-2)^{+k} \;\lra\; 0 
\]
and
\[
0\;\lra\;\EE(\be+1/\be-1\bby\be/\al-2)^{+k}\;\lraa{e}\;\EE(\be+1/\be-1\bby\be/\al-1)^{+k}\;\lraa{e}\;\EE(\be-1/\al-2\bby\al-1/\al-3)^{+k+1}\; .
\]
 
{\it Claim 2.} We have an isomorphism $\EE(\ga/\be-1\bby\be/\al-1)^{+k}(f)$ for all $k\in\Z$ and all $\al,\,\be,\,\ga\,\in\,\Z$ such that $\al < \be < \ga$. 

We proceed by induction on $\ga - \al$. By Claim 1, we may assume that $\al < \be - 1 < \be + 1 < \ga$. Consider the image diagram
\[
\EE(\ga-1/\be-1\bby\be/\al-1)^{+k}\;\lraepia{e}\;\EE(\ga/\be-1\bby\be/\al-1)^{+k}\;\lramonoa{e}\;\EE(\ga/\be-1\bby\be/\al)^{+k}\; .
\] 

{\it Claim 3.} We have an isomorphism $\EE(\de/\be\bby\ga/\al)^{+k}(f)$ for all $k\in\Z$ and all $\al,\,\be,\,\ga,\,\de\,\in\,\Z$ such that $\al < \be \le \ga < \de$.

We may assume that $\ga - \be \ge 1$, for $\EE(\de/\be\bby\be/\al)^{+k} = 0$. We proceed by induction on $\ga - \be$. By Claim 2, we may assume that $\ga - \be \ge 2$. Consider the short exact sequence
\[
\EE(\de/\be\bby\ga-1/\al)^{+k}\;\lramonoa{e}\;\EE(\de/\be\bby\ga/\al)^{+k}\;\lraepia{e}\;\EE(\de/\ga-1\bby\ga/\al)^{+k} \; .
\]

{\it Claim 4.} We have an isomorphism $\EE(\de/\be\bby\ga/\al)^{+k}(f)$ for all $k\in\Z$ and all $\al,\,\be,\,\ga,\,\de\,\in\,\Zi$ such that $\al < \be \le \ga < \de$.

In view of Claim 3, it suffices to choose $\w\al\in\Z$ small enough such that $\EE(\de/\be\bby\ga/\w\al)^{+k}(f) = \EE(\de/\be\bby\ga/\!-\!\infty)^{+k}(f)$; etc.

{\it Claim 5.} We have an isomorphism $\EE(\de/\be\bby\ga/\al)^{+k}(f)$ for all $k\in\Z$ and all $\al,\,\be,\,\ga,\,\de\,\in\,\Zi$ such that $\al \dl \be \le \ga \dl \de$.

In view of Claim 4, it suffices to choose $\w\be\in\Z$ small enough such that $\EE(\de/\w\be\bby\ga/\!-\!\infty)^{+k}(f) = \EE(\de/\!-\!\infty\bby\ga/\!-\!\infty)^{+k}(f)$; etc.

{\it Claim 6.} We have an isomorphism $\EE(\de/\be\bby\ga/\al)^{+k}(f)$ for all $k\in\Z$ and all $\al,\,\be,\,\ga,\,\de\,\in\,\bZi$ such that 
$-\infty \le \de^{-1} \le \al \dl \be \le \ga \le \infty < -\infty^{+1} \le \de\le\al^{+1}$.

In view of Claim 5, it suffices to consider the short exact sequence
\[
\EE(\infty/\be\bby\ga/\de^{-1})^{+k}\;\lramonoa{e}\;\EE(\infty/\be\bby\ga/\al)^{+k}\;\lraepia{e}\;\EE(\de/\be\bby\ga/\al)^{+k} \; .
\]

{\it Claim 7.} The morphism $\EEp(f)$ is an isomorphism.

Suppose given $\al,\,\be,\,\ga,\,\de\,\in\,\bZi$ such that $\de^{-1}\le \al \dl \be \le \ga \dl \de \le \al^{+1}$. Via a shift, we may assume that we are in the situation of Claim 5 or of Claim 6.
\qed

\subsection{The first spectral sequence of a double complex}

Let $\Al$ be an abelian category. Let $X\in\Ob\CCCP(\Al)$. Given $n\in\Zi$, we write $X^{[n,\ast}$ for the double complex arising from $X$ by replacing $X^{i,j}$ by $0$ for all $i\in [0,n[$.
We define a pointwise split and pointwise finitely filtered complex $\tI X$, called the {\it first filtration of $\t X$,} by letting $\tI X(\al) := \t X^{[-\al,\ast}$ for $\al\in\Zi$; and by letting 
$\tI X(\al)\lra\tI X(\be)$ be the pointwise split inclusion $\t X^{[-\al,\ast}\lra \t X^{[-\be,\ast}$ for $\al,\,\be\,\in\,\Zi$ such that $\al\le\be$. Let $\EEI = \EEI(X) := \EE(\tI X)$. This 
construction is functorial in $X\in\Ob\CCCP(\Al)$. Note that $\ol{\tI X}(\al) = X^{-\al,k+\al}$.

\bq
 We record the following wellknown lemma in the language we use here.

\eq

\begin{Lemma}
\label{LemSSDC1}
Let $\al\in\, ]\!-\!\infty,0]$. Let $k\in\Z$ such that $k\ge -\al$. We have
\[
\barcl
\EEI(\al/\al-1\bby\al/\al-1)^{+k}(X)   & = & \HH^{k+\al}( X^{-\al,\ast} ) \\
\EEI(\al+1/\al-1\bby\al/\al-2)^{+k}(X) & = & \HH^{-\al}\big(\HH^{k+\al}( X^{-,\ast})\big)\; , \\
\ea
\]
naturally in $X\in\Ob\CCCP(\Al)$.
\end{Lemma}

{\it Proof.} The first equality follows by $\EEI(\al/\al-1\bby\al/\al-1)^{+k} = \HH^k \tI X(\al/\al-1) = \HH^{k+\al} (X^{-\al,\ast})$.

The morphism $\tI X(\al/\al-1) \lra \tI X\big((\al - 2)^{+1}/\al-1\big) = \tI X\big(\al-1/\al-2\big)^{\bt+1}$ from $\text{Sp}(\tI X)$ is at position $k\ge 0$ given by
\[
\ol{\tI X}(\al)^k = X^{-\al,k+\al}\;\mrafl{27}{(-1)^\al\,\dell}\; X^{-\al+1,k+\al} = \ol{\tI X}(\al-1)^{k+1}\; ;
\]
cf.\ \S\ref{SecTotal}. In particular, the morphisms
\[
\EEI(\al+1/\al\bby\al+1/\al)^{+k-1}\;\lraa{e}\;\EEI(\al/\al-1\bby\al/\al-1)^{+k}\;\lraa{e}\;\EEI(\al-1/\al-2\bby\al-1/\al-2)^{+k+1}
\]
are given by
\[
\HH^{k+\al}(X^{-\al-1,\ast}) \;\;\mrafl{30}{(-1)^{\al+1}\HH^{k+\al}(\dell)}\;\;\HH^{k+\al}(X^{-\al,\ast}) \;\;\mrafl{30}{(-1)^\al\HH^{k+\al}(\dell)}\;\; \HH^{k+\al}(X^{-\al+1,\ast}) \; .
\]
Now the second equality follows by the diagram
\[
\xymatrix@C=5mm{
                                                                              & \EEI(\al+1/\al-1\bby\al/\al-2)^{+k}\ar~+{|*\dir{*}}[rd]^e        &                                               \\
\EEI(\al/\al-1\bby\al/\al-2)^{+k}\ar~+{|*\dir{*}}[rd]^e\ar~+{|*\dir{|}}[ru]^e &                                                                  & \EEI(\al+1/\al-1\bby\al/\al-1)^{+k}           \\
\EEI(\al+1/\al\bby\al+1/\al)^{+k-1} \ar[r]^-e                                 & \EEI(\al/\al-1\bby\al/\al-1)^{+k}\ar[r]^-e\ar~+{|*\dir{|}}[ru]^e & \EEI(\al-1/\al-2\bby\al-1/\al-2)^{+k+1}\zw{.} \\
}
\]
\qed

\begin{Remark} 
\label{RemSSDC2}\hspace*{3mm}
Let $X\lraa{f} Y$ be a rowwise quasiisomorphism in $\CCCP(\Al)$. Then $\EEI(\de/\be\bby\ga/\al)^{+k}(f)$ is an isomorphism for $\de^{-1}\le\al\le\be\le\ga\le\de\le\al^{+1}$ in $\bZi$ and $k\in\Z$.
\end{Remark}

{\it Proof.} It suffices to show that the morphism $\text{Sp}(\tI f)$ in $\bo\bZif,\,\KK(\Al)\bc$ is pointwise a quasiisomorphism. To have this, it suffices to show that $\t f^{[k,\ast}$ is a 
quasiisomorphism for $k\ge 0$. But $f^{[k,\ast}$ is a rowwise quasiisomorphism for $k\ge 0$; cf.\ \S\ref{SecTotal}.\qed

\begin{Lemma} 
\label{LemSSDC3}
The functor $\;\CCCP(\Al)\;\lraa{\EEIp}\; \bo\bZiffp,\,\Al\bc\;$ factors over
\[
\KKKP(\Al)\;\lraa{\EEIp}\; \bo\bZiffp,\,\Al\bc \; .
\]
\end{Lemma}

{\it Proof.} By Lemma \ref{LemDTC1}, we have to show that $\EEIp$ annihilates all elementary horizontally split acyclic double complexes in $\Ob\CCCP(\Al)$ and all elementary vertically split acyclic 
double complexes in $\Ob\CCCP(\Al)$.

Let $U\in\Ob\CCCP(\Al)$ be an elementary vertically split acyclic double complex concentrated in rows $i$ and $i+1$, where $i\ge 0$. Let $V\in\Ob\CCCP(\Al)$ be an elementary horizontally split acyclic 
double complex concentrated in columns $j$ and $j+1$, where $j\ge 0$.

Since $V$ is rowwise acyclic, $\EEI$ annihilates $V$ by Remark~\ref{RemSSDC2}, whence so does $\EEIp$.

Suppose given 
\[
-\infty\le\al\dl\be\le\ga\dl\de\le\infty
\leqno (\ast)
\]
in $\bZi$ and $k\in\Z$. We {\it claim} that the functor $\EEI(\de/\be\bby\ga/\al)^{+k}$ annihilates $U$. We may assume that $\be < \ga$. Note that $\EEI(\de/\be\bby\ga/\al)^{+k}(U)$ is the image
of 
\[
\HH^k\big(\tI U(\ga/\al)\big)\;\lra\;\HH^k\big(\tI U(\de/\be)\big)\; .
\]

The double complex $U^{[-\de,\ast}/U^{[-\be,\ast}$ is columnwise acyclic except possibly if $-\be = i+1$ or if $-\de = i+1$.
The double complex $U^{[-\ga,\ast}/U^{[-\al,\ast}$ is columnwise acyclic except possibly if $-\al = i+1$ or if $-\ga = i+1$. All three remaining combinations of these exceptional cases are 
excluded by $(\ast)$, however. Hence $\EEI(\de/\be\bby\ga/\al)^{+k}(U) = 0$. This proves the {\it claim.}

Suppose given 
\[
\de^{-1}\le\al\dl\be\le\ga\le\infty\le-\infty^{+1}\le\de\le\al^{+1} \; . 
\leqno (\ast\ast)
\]
in $\bZi$ and $k\in\Z$. We {\it claim} that the functor $\EEI(\de/\be\bby\ga/\al)^{+k}$ annihilates $U$. We may assume that 
$\be < \ga$ and that $\de^{-1} < \al$. Note that $\EEI(\de/\be\bby\ga/\al)^{+k}(U)$ is the image of 
\[
\HH^k\big(\tI U(\ga/\al)\big) \;\lra\; \HH^{k+1}\big(\tI U(\be/\de^{-1})\big)\; .
\]
 
The double complex $U^{[-\be,\ast}/U^{[-(\de^{-1}),\ast}$ is columnwise acyclic except possibly if $-(\de^{-1}) = i+1$ or if $-\be = i+1$. The double complex $U^{[-\ga,\ast}/U^{[-\al,\ast}$ is columnwise 
acyclic except possibly if $-\ga = i+1$ or if $-\al = i+1$. Both remaining combinations of these exceptional cases are excluded by $(\ast\ast)$, however.
Hence $\EEI(\de/\be\bby\ga/\al)^{+k}(U) = 0$. This proves the {\it claim.}

Both claims taken together show that $\EEIp$ annihilates $U$.\qed

\section{Grothendieck spectral sequences}
\label{SecGroth}

\subsection{Certain quasiisomorphisms are preserved by a left exact functor}
\label{SecALemma}

Suppose given abelian categories $\Al$, $\Bl$, and suppose that $\Al$ has enough injectives. Let $\Al\lraa{F}\Bl$ be a left exact functor.

\begin{Remark} 
\label{RemGr1}
Suppose given an $F$-acyclic object \mb{$X\in\Ob\Al$} and an injective resolution \mb{$I\in\Ob\CC^{[0}(\Inj\Al)$} of $X$. Let $\Conc X\lraa{f} I$ be its quasiisomorphism. Then $\Conc FX\lraa{Ff} FI$ is a 
quasiisomorphism.
\end{Remark}

{\it Proof.} This follows since $F$ is left exact and since $\HH^i (FI) \iso (\RR^i F) X \iso 0$ for $i \ge 1$.\qed

\begin{Remark} 
\label{RemGr1_5}
Suppose given a complex $U\in\Ob\CC^{[0}(\Al)$ consisting of $F$-acyclic objects. There exists an injective complex resolution $I\in\Ob\CC^{[0}(\Inj\Al)$ of $U$ such that its quasiisomorphism 
$U\lraa{f} I$ maps to a quasiisomorphism $FU\lraa{Ff} FI$.
\end{Remark}

{\it Proof.} Let $J\in\Ob\CCCCE(\Inj\Al)$ be a CE-resolution of $U$; cf.\ Remark \ref{RemCEqis2_5}. Since the morphism of double complexes $\Conc_2 U\lra J$ is a columnwise quasiisomorphism consisting
of monomorphisms, taking the total complex, we obtain a quasiisomorphism $U\lra \t J$ consisting of monomorphisms. By $F$-acyclicity of the entries of $U$, the image $\Conc_2 FU\lra FJ$ under $F$ 
is a columnwise quasiisomorphism, too; cf.\ Remark~\ref{RemGr1}. Hence $F$ maps the quasiisomorphism $U\lra\t J$ to the quasiisomorphism $FU\lra F\t J$. So we may take $I := \t J$. \qed

\begin{Lemma}
\label{LemCEqis4}
Suppose given a complex $U\in\Ob\CC^{[0}(\Al)$ consisting of $F$-acyclic objects and an injective complex resolution $I\in\Ob\CC^{[0}(\Inj\Al)$ of $U$. Let $U\lraa{f} I$ be its quasiisomorphism. 
Then $FU\lraa{Ff} FI$ is a quasiisomorphism.
\end{Lemma}

{\it Proof.} Let $U\lra I'$ be a quasiisomorphism to an injective complex resolution $I'$ that is mapped to a quasiisomorphism by $F$; cf.\ Remark \ref{RemGr1_5}. Since $U\lra I'$ is a quasiisomorphism, 
the induced map $\liu{\KK(\Al)}{(U,I)} \lla \liu{\KK(\Al)}{(I',I)}$ is surjective, so that there exists a morphism $I'\lra I$ such that $(U\lra I'\lra I) \= (U\lraa{f} I)$ in $\KK(\Al)$. Since, moreover, 
$U\lraa{f} I$ is a quasiisomorphism, $I'\lra I$ is a homotopism. Since $FU\lra FI'$ is a quasiisomorphism and $FI'\lra FI$ is a homotopism, we conclude that $FU\lra FI$ is a quasiisomorphism.\qed

\subsection{Definition of the Grothendieck spectral sequence functor}
\label{SecGrothDef}

Suppose given abelian categories $\Al$, $\Bl$ and $\Cl$, and suppose that $\Al$ and $\Bl$ have enough injectives. Let $\Al\lraa{F}\Bl$ and $\Bl\lraa{G}\Cl$ be left exact functors. 

A {\it $(F,G)$-acyclic resolution} of $X\in\Ob\Al$ is a complex $A\in\Ob\CC^{[0}(\Al)$, together with a quasiisomorphism $\Conc X\lra A$, such that the following hold.

\begin{itemize}
\item[(A\,1)] The object $A^i$ is $F$-acyclic for $i\ge 0$.
\item[(A\,2)] The object $A^i$ is $(G\0 F)$-acyclic for $i\ge 0$.
\item[(A\,3)] The object $F A^i$ is $G$-acyclic for $i\ge 0$.
\end{itemize}

An object $X\in\Ob\Al$ that possesses an $(F,G)$-acyclic resolution is called {\it $(F,G)$-acyclicly resolvable.} The full subcategory of $(F,G)$-acyclicly resolvable objects in $\Al$ 
is denoted by $\Al_{(F,G)}$.

A complex $A\in\Ob\CC^{[0}(\Al)$, together with a quasiisomorphism $\Conc X\lra A$, is called an {\it $F$\!\nobreakdash-acyclic resolution} of $X\in\Ob\Al$ if (A\,2) holds.

\begin{Remark} 
\label{RemGr0}
If $F$ carries injective objects to $G$-acyclic objects, then {\rm (A\,1)} and {\rm (A\,3)} imply {\rm (A\,2)}.
\end{Remark}

{\it Proof.} Given $i\ge 0$, we let $I$ be an injective resolution of $A^i$, and $\w I$ the acyclic complex obtained by appending $A^i$ to $I$ in position $-1$. Since $A^i$ is $F$-acyclic,
the complex $F\w I$ is acyclic; cf.\ Remark \ref{RemGr1}. Note that $F\BB^0\w I\iso F A^i$ is $G$-acyclic by assumption. Since 
\[
(\RR^k G) F\w I^j\;\lra\; (\RR^k G) F\BB^{j+1}\w I\;\lra\;(\RR^{k+1}G) F\BB^j\w I
\]
is exact in the middle for $j\ge 0$ and $k\ge 1$, we may conclude by induction on $j$ and by $G$\nobreakdash-acyclicity assumption on $F\w I^j$ that $F\BB^j\w I$ is $G$-acyclic for $j\ge 0$. 
In particular, we have $(\RR^1 G)(F\BB^j\w I) \iso 0$ for $j\ge 0$, whence 
\[
GF\BB^j\w I\;\lra\; GF\w I^j\;\lra\; GF\BB^{j+1}\w I
\]
is short exact for $j\ge 0$. We conclude that $(G\0 F)\w I$ is acyclic. Hence $A^i$ is $(G\0 F)$-acyclic.\qed

\bq
 To see Remark \ref{RemGr0}, one could also use a Grothendieck spectral sequence, once established.

\eq

\begin{Remark} 
\label{RemGr2} 
Suppose given $X\in\Ob\Al$, an injective resolution $I$ of $X$ and an $F$-acyclic re\-solution $A$ of $X$. Then there exists a quasiisomorphism $A\lra I$ that is mapped to $1_X$ by $\HH^0$. Moreover,
any morphism $A\lraa{u} I$ that is mapped to $1_X$ by $\HH^0$ is a quasiisomorphism and is mapped to a quasiisomorphism $FA\lrafl{26}{Fu} FI$ by $F$.
\end{Remark}

{\it Proof.} Let $I'$ be an injective complex resolution of $A$ such that its quasiisomorphism $A\lra I'$ is mapped to a quasiisomorphism by $F$; cf.\ Remark \ref{RemGr1_5}. We use the 
composite quasiisomorphism $\Conc X\lra A\lra I'$ to resolve $X$ by $I'$.

To prove the first assertion, note that there is a homotopism $I'\lra I$ resolving $1_X$; whence the composite $(A\lra I'\lra I)$ is a quasiisomorphism resolving $1_X$. 

To prove the second assertion, note that the induced map $\liu{\KK(\Al)}{(A,I)} \lla \liu{\KK(\Al)}{(I',I)}$ is surjective, whence there is a factorisation $(A\lra I'\lra I) = (A\lraa{u} I)$ 
in $\KK(\Al)$ for some morphism $I'\lra I$, which, since resolving $1_X$ as well, is a homotopism. In particular, $A\lraa{u} I$ is a quasiisomorphism. Finally, since $FI'\lra FI$ is a homotopism, 
also $FA\lraa{Fu} FI$ is a quasiisomorphism.\qed

\bq
 Alternatively, in the last step of the preceding proof we could have invoked Lemma \ref{LemCEqis4}.

 The following construction originates in \bfcite{CE56}{XVII.\S 7} and \bfcite{Gr57}{Th.\ 2.4.1}. In its present form, it has been carried out by {\sc Haas} in the classical framework~\bfcit{Ha77}. 
 We do not claim any originality.

 I do not know whether the use of injectives in $\Al$ in the following construction can be avoided; in any case, it would be desirable to do so.

\eq

We set out to define the {\it proper Grothendieck spectral sequence functor}
\[
\Al_{(F,G)} \;\lrafl{30}{\EEGrp{F}{G}}\; \bo\bZiffp,\,\Cl\bc\; .
\]

{\it We define $\EEGrp{F}{G}$ on objects.} Suppose given $X\in\Ob\Al_{(F,G)}$. Choose an $(F,G)$-acyclic resolution $A_X\in\Ob\CC^{[0}(\Al)$ of $X$. Choose a CE\nobreakdash-resolution 
$J_X\in\Ob\CCCP(\Inj\Bl)$ of $FA_X$. Let $\EEGr{F}{G}(X) := \EEI(GJ_X) = \EE(\tI GJ_X)\in\Ob\bo\bZiff,\,\Cl\bc$ be the {\it Grothendieck spectral sequence} of $X$ with respect to $F$ and $G$. 
Accordingly, let 
\[
\EEGrp{F}{G}(X) \; :=\; \EEIp(GJ_X) \; =\; \EEp(\tI GJ_X) \;\in\;\Ob\bo\bZiffp,\,\Cl\bc
\]
be the {\it proper Grothendieck spectral sequence} of $X$ with respect to $F$ and $G$.

{\it We define $\EEGrp{F}{G}$ on morphisms.} Suppose given $X\in\Ob\Al_{(F,G)}$, and let $A_X$ and $J_X$ be as above. Choose an injective resolution $I_X\in\Ob\CC^{[0}(\Inj\Al)$ of $X$. 
Choose a quasiisomorphism $A_X\lraa{p_X} I_X$ that is mapped to $1_X$ by $\HH^0$ and to a quasiisomorphism by $F$; cf.\ Remark \ref{RemGr2}. Choose a CE\nobreakdash-resolution 
$K_X\in\Ob\CCCP(\Inj\Bl)$ of $FI_X$. Choose a morphism $J_X\lraa{q_{X}}K_X$ in $\CCCP(\Inj\Bl)$ that is mapped to $F p_X$ by $\HH^0\big((-)^{\ast,-}\big)$; cf.\ Remark~\ref{RemEx4}.

Note that $J_X\lraa{q_{X}}K_X$ can be written as a composite in $\CCCCE(\Inj\Bl)$ of a rowwise homotopism, followed by a double homotopism; cf.\ Proposition \ref{PropCEqis9}. Hence, so can
$G J_X\lrafl{25}{G q_{X}} GK_X$. Thus $\EEIp(G J_X)\mrafl{25}{\EEIp(G q_{X})} \EEIp(GK_X)$ is an isomorphism; cf.\ Remark \ref{RemSSDC2}, Lemma \ref{LemSSDC3}.

Suppose given $X\lraa{f} Y$ in $\Al_{(F,G)}$. Choose a morphism $I_X\lraa{f'} I_Y$ in $\CC^{[0}(\Al)$ that is mapped to $f$ by $\HH^0$. Choose a morphism $K_X\lraa{f''}K_Y$ in 
$\CCCP(\Inj\Bl)$ that is mapped to $Ff'$ by $\HH^0\big((-)^{\ast,-}\big)$; cf.\ Remark~\ref{RemEx4}. Let
\[
\EEGrp{F}{G}(X\lraa{f} Y) \; :=\; \left(\EEIp(G J_X)\mraisofl{25}{\EEIp(G q_{X})} \EEIp(GK_X) \mrafl{25}{\EEIp(Gf'')} \EEIp(GK_Y) \mlaisofl{25}{\EEIp(G q_Y)} \EEIp(G J_Y) \right)\; .
\]

\bq
 The procedure can be adumbrated as follows.
 \begin{center}
 \begin{picture}(400,650)
 \put( 100,   0){$X$}

 \put( 150,  10){\vector(1,0){230}}
 \put( 260,  25){$\scm f$}

 \put( 400,   0){$Y$}

 \put(   0, 200){$A_X$}
 \put(  40, 240){\vector(1,1){50}}
 \put(  20, 270){$\scm p_X$}
 \put( 100, 300){$I_X$}

 \put( 150, 310){\vector(1,0){230}}
 \put( 260, 325){$\scm f'$}

 \put( 300, 200){$A_Y$}
 \put( 340, 240){\vector(1,1){50}}
 \put( 320, 270){$\scm p_Y$}
 \put( 400, 300){$I_Y$}

 \put(   0, 500){$J_X$}
 \put(  40, 540){\vector(1,1){50}}
 \put(  20, 570){$\scm q_X$}
 \put( 100, 600){$K_X$}

 \put( 170, 610){\vector(1,0){210}}
 \put( 260, 625){$\scm f''$}

 \put( 300, 500){$J_Y$}
 \put( 340, 540){\vector(1,1){50}}
 \put( 320, 570){$\scm q_Y$}
 \put( 400, 600){$K_Y$}
 \end{picture}
 \end{center}

\eq

{\it We show that this defines a functor $\EEGrp{F}{G}:\Al_{(F,G)}\lra\bo\bZiffp,\,\Cl\bc$.} We need to show independence of the construction from the choices of $f'$ and $f''$, for then 
functoriality follows by appropriate choices.

Let $I_X\lraa{\w f'} I_Y$ and $K_X\lraa{\w f''}K_Y$ be alternative choices. The residue classes of $f'$ and $\w f'$ in $\KK^{[0}(\Al)$ coincide, whence so do the residue classes of 
$F f'$ and $F \w f'$ in $\KK^{[0}(\Bl)$. Therefore, the residue classes of $f''$ and $\w f''$ in $\KKKP(\Bl)$ coincide; cf.\ Proposition \ref{PropCEqis6}. Hence, so do the residue classes
of $G f''$ and $G\w f''$ in $\KKKP(\Cl)$. Thus $\EEIp(G f'') = \EEIp(G \w f'')$; cf.\ Lemma \ref{LemSSDC3}.

{\it We show that alternative choices of $A_X$, $I_X$ and $p_X$, and of $J_X$, $K_X$ and $q_X$, yield isomorphic proper Grothendieck spectral sequence functors.}

Let $\w A_X\lraa{\w p_X} \w I_X$ and $\w J_X\lraa{\w q_X} \w K_X$ be alternative choices, where $X$ runs through $\Ob\Al_{(F,G)}$.

Suppose given $X\lraa{f} Y$ in $\Al_{(F,G)}$. We resolve the commutative quadrangle
\[
\xymatrix{
X\ar[r]^f\ar@{=}[d] & Y\ar@{=}[d] \\
X\ar[r]^f           & Y           \\
}
\]
in $\Al$ to a commutative quadrangle
\[
\xymatrix{
I_X\ar[r]^{f'}\ar[d]_{u_X}    & I_Y\ar[d]^{u_Y} \\
\w I_X\ar[r]^{\w f'}          & \w I_Y    \\
}
\]
in $\KK^{[0}(\Al)$, in which $u_X$ and $u_Y$ are homotopisms; cf.\ Remark~\ref{RemEx4}. Then we resolve the commutative quadrangle
\[
\xymatrix{
FI_X\ar[r]^{F f'}\ar[d]_{F u_X}   & FI_Y\ar[d]^{F u_Y} \\
F\w I_X\ar[r]^{F\w f'}            & F\w I_Y    \\
}
\]
in $\KK^{[0}(\Bl)$ to a commutative quadrangle
\[
\xymatrix{
K_X\ar[r]^{f''}\ar[d]_{v_X}    & K_Y\ar[d]^{v_Y} \\
\w K_X\ar[r]^{\w f''}          & \w K_Y    \\
}
\]
in $\KKKP(\Bl)$; cf.\ Proposition \ref{PropCEqis6}. Therein, $v_X$ and $v_Y$ are each composed of a rowwise homotopism, followed by a double homotopism; cf.\ Proposition \ref{PropCEqis9}. 
So are $G v_X$ and $G v_Y$. An application of $\EEIp\big(G (-)\big)$ yields the sought isotransformation, viz.\
\[
\left(\EEIp(G J_X)\mraisofl{25}{\EEIp(G q_{X})} \EEIp(G K_X) \mraisofl{25}{\EEIp(G v_X)} \EEIp(G \w K_X) \mlaisofl{25}{\EEIp(G \w q_X)} \EEIp(G \w J_X)\right)
\]
at $X\in\Ob\Al_{(F,G)}$; cf.\ Remark \ref{RemSSDC2}, Lemma \ref{LemSSDC3}.


\bq
 Finally, we recall the starting point of the whole enterprise.

\eq

\begin{Remark}[\bfcite{CE56}{XVII.\S 7}, \bfcite{Gr57}{Th.\ 2.4.1}]
\label{RemGr2_5}
Suppose given $X\in\Ob\Al_{(F,G)}$ and $k,\,\ell\,\in\,\Z_{\ge 0}$. We have
\[
\ba{lcl}
\EEGrp{F}{G}(-k+1/\!-\!k-1\bby\!-\!k/\!-\!k-2)^{+k+\ell}(X)         & \iso & (\RR^k G)(\RR^\ell F)(X)             \\
\EEGrp{F}{G}(\infty/\!-\!\infty\bby\infty/\!-\!\infty)^{+k+\ell}(X) & \iso & \big(\RR^{k+\ell}(G\0 F)\big)(X)\; , \\
\ea
\]
naturally in $X$.
\end{Remark}

{\it Proof.} Keep the notation of the definition of $\EEGrp{F}{G}\,$.

We shall prove the first isomorphism. By Lemma~\ref{LemSSDC1}, we have
\[
\EEGrp{F}{G}(-k+1/\!-\!k-1\bby\!-\!k/\!-\!k-2)^{+k+\ell}(X)\;\iso\; \HH^k(\HH^\ell(G J_X^{-,\ast}))\; .
\]
Since $J_X$ is rowwise split, we have $\HH^\ell(G J_X^{-,\ast}) \iso G(\HH^\ell J_X^{-,\ast})$. Note that $\HH^\ell J_X^{-,\ast}$ is an injective resolution of $\HH^\ell FA_X$; 
cf.\ Remark \ref{RemCEqis1}.(1). By Remark~\ref{RemGr2}, $\HH^\ell F A_X\lraisofl{28}{\HH^\ell F p_X} \HH^\ell F I_X \iso (\RR^\ell F)(X)$. So
\[
\HH^k(\HH^\ell(G J_X^{-,\ast}))\;\iso\; \HH^k(G(\HH^\ell J_X^{-,\ast})) \;\iso\; (\RR^k G)(\HH^\ell F A_X)\;\iso\; (\RR^k G)(\RR^\ell F)(X)\; .
\]
We shall prove naturality of the first isomorphism. Suppose given $X\lraa{f} Y$ in $\Al_{(F,G)}$. Consider the following commutative diagram. Abbreviate 
$E := \EEp(-k+1/\!-\!k-1\bby\!-\!k/\!-\!k-2)^{+k+\ell}$.
{\footnotesize
\[
\xymatrix@C=20mm{
E(\tI GJ_X)\ar[r]^{E(\tI G q_X)}_\sim\ar[d]^\wr                                   & E(\tI GK_X)\ar[r]^{E(\tI Gf'')}\ar[d]^\wr                                   & E(\tI GK_Y)\ar[d]^\wr                  & E(\tI GJ_Y)\ar[l]_{E(\tI G q_Y)}^\sim\ar[d]^\wr                                   \\
\HH^k\HH^\ell GJ_X^{-,\ast}\ar[r]^{\HH^k\HH^\ell G q_X^{-,\ast}}_\sim\ar[d]^\wr   & \HH^k\HH^\ell GK_X^{-,\ast}\ar[r]^{\HH^k\HH^\ell Gf''^{-,\ast}}\ar[d]^\wr   & \HH^k\HH^\ell GK_Y^{-,\ast}\ar[d]^\wr  & \HH^k\HH^\ell GJ_Y^{-,\ast} \ar[l]_{\HH^k\HH^\ell G q_Y^{-,\ast}}^\sim\ar[d]^\wr  \\
\HH^k G \HH^\ell J_X^{-,\ast}\ar[r]^{\HH^k G\HH^\ell q_X^{-,\ast}}_\sim\ar[d]^\wr & \HH^k G\HH^\ell K_X^{-,\ast}\ar[r]^{\HH^k G\HH^\ell f''^{-,\ast}}\ar[d]^\wr & \HH^k G\HH^\ell K_Y^{-,\ast}\ar[d]^\wr & \HH^k G\HH^\ell J_Y^{-,\ast} \ar[l]_{\HH^k G\HH^\ell q_Y^{-,\ast}}^\sim\ar[d]^\wr \\
(\RR^k G) \HH^\ell FA_X\ar[r]^{(\RR^k G)\HH^\ell Fp_X}_\sim                       & (\RR^k G)\HH^\ell FI_X\ar[r]^{(\RR^k G)\HH^\ell F f'}\ar[d]^\wr             & (\RR^k G)\HH^\ell F I_Y\ar[d]^\wr      & (\RR^k G)\HH^\ell FA_Y \ar[l]_{(\RR^k G)\HH^\ell Fp_Y}^\sim                       \\
                                                                                  & (\RR^k G)(\RR^\ell F)(X)\ar[r]^{(\RR^k G)(\RR^\ell F)(f)}                   & (\RR^k G)(\RR^\ell F)(Y)               &                                                                                   \\
}
\]
}%

We shall prove the second isomorphism. By Lemma~\ref{LemCEqis4}, the quasiisomorphism $FA_X\lra \t J_X$ maps to a quasiisomorphism $GFA_X\lra \t GJ_X \iso G\t J_X$. By Lemma~\ref{LemCEqis4}, 
the quasiisomorphism $A_X\lraa{p_X} I_X$ maps to a quasiisomorphism $GFA_X\lrafl{25}{GFp_X} GFI_X$. So
\[
\ba{rcl}
\EEGrp{F}{G}(\infty/\!-\!\infty\bby\infty/\!-\!\infty)^{+k+\ell}(X) & \hspace*{-1mm}\iso\hspace*{-1mm} & \HH^{k+\ell}(\t GJ_X)\;\iso\; \HH^{k+\ell}(G\t J_X)\;\iso\; \HH^{k+\ell}(GFA_X) \\
                                                                    & \hspace*{-1mm}\iso\hspace*{-1mm} & \HH^{k+\ell}(GFI_X) \;\iso\; \big(\RR^{k+\ell}(G\0 F)\big)(X)\; . \\
\ea
\]
We shall prove naturality of the second isomorphism. Consider the following diagram. Abbreviate $\w E := \EEGrp{F}{G}(\infty/\!-\!\infty\bby\infty/\!-\!\infty)^{+k+\ell}$.
{\footnotesize
\[
\xymatrix@C=20mm{
\w E(\tI GJ_X)\ar[r]^{\w E(\tI G q_X)}_\sim\ar[d]^\wr           & \w E(\tI GK_X)\ar[r]^{\w E(\tI Gf'')}\ar[d]^\wr                                              & \w E(\tI GK_Y)\ar[d]^\wr            & \w E(\tI GJ_Y)\ar[l]_{\w E(\tI G q_Y)}^\sim\ar[d]^\wr            \\
\HH^{k+\ell}\t GJ_X\ar[r]^{\HH^{k+\ell}\t G q_X}_\sim\ar[d]^\wr & \HH^{k+\ell}\t GK_X\ar[r]^{\HH^{k+\ell}\t Gf''}\ar[d]^\wr                                    & \HH^{k+\ell}\t GK_Y\ar[d]^\wr       & \HH^{k+\ell}\t GJ_Y \ar[l]_{\HH^{k+\ell}\t G q_Y}^\sim\ar[d]^\wr \\
\HH^{k+\ell} G\t J_X\ar[r]^{\HH^{k+\ell} G \t q_X}_\sim         & \HH^{k+\ell} G\t K_X\ar[r]^{\HH^{k+\ell} G\t f''}                                            & \HH^{k+\ell} G\t K_Y                & \HH^{k+\ell}G\t J_Y \ar[l]_{\HH^{k+\ell} G\t q_Y}^\sim           \\
\HH^{k+\ell} GFA_X\ar[r]^{\HH^{k+\ell} GFp_X}_\sim\ar[u]_\wr    & \HH^{k+\ell} GFI_X\ar[r]^{\HH^{k+\ell} GFf'}\ar[d]^\wr\ar[u]                                 & \HH^{k+\ell} GF I_Y\ar[d]^\wr\ar[u] & \HH^{k+\ell} GFA_Y \ar[l]_{\HH^{k+\ell}GFp_Y}^\sim\ar[u]_\wr     \\
                                                                & \big(\RR^{k+\ell}(G\0 F)\big)(X)\ar[r]^{(\RR^{k+\ell}(G\0 F))(f)} & \big(\RR^{k+\ell}(G\0 F)\big)(Y)    &                                                                                             \\
}
\]
}%
\qed

\subsection{Haas transformations}
\label{SecGrothHaas}

\bq
 The following transformations have been constructed in the classical framework by {\sc Haas}~\bfcit{Ha77}. We do not claim any originality.

\eq

\subsubsection{Situation}
\label{SecHaasSit}


Consider the following diagram of abelian categories, left exact functors and transformations,
\[
\xymatrix{
\Al\ar[r]^F \ar[d]_U & \Bl\ar[r]^G\ar[d]_V & \Cl\ar[d]_W \\
\Al'\ar[r]^{F'}      & \Bl'\ar[r]^{G'}     & \Cl'\; , \!\!  \\
\ar@2"2,1"+<5mm,5mm>;"1,2"+<-5mm,-5mm>^\mu 
\ar@2"2,2"+<5mm,5mm>;"1,3"+<-5mm,-5mm>^\nu 
}
\] 
i.e.\ $F'\0 U\lraa{\mu} V\0 F$ and $G'\0 V\lraa{\nu} W\0 G$. Suppose that the conditions (1,\,2,\,3) hold.
\begin{itemize}
\item[(1)] The categories $\Al$, $\Bl$, $\Al'$ and $\Bl'$ have enough injectives. 
\item[(2)] The functors $U$ and $V$ carry injectives to injectives.
\item[(3)] The functor $F$ carries injective to $G$-acyclic objects. The functor $F'$ carries injective to $G'$-acyclic objects. 
\end{itemize}

We have $\Al_{(F,G)} = \Al$ since an injective resolution is an $(F,G)$-acyclic resolution. Likewise, we have $\Al'_{(F',G')} = \Al'$.

Note in particular the case $U = 1_\Al\,$, $V = 1_\Bl$ and $W = 1_\Cl\,$.

We set out to define the {\it Haas transformations}
\[
\EEGrp{F'}{G'}\big(U(-)\big) \;\lrafl{28}{\text{h}^\text{I}_\mu}\; \EEGrp{F}{G'\0 V}\big(-\big) \;\lrafl{28}{\text{h}^\text{II}_\nu}\; \EEGrp{F}{W\0 G}\big(-\big)\; ,
\]
where $\text{h}^\text{I}_\mu$ depends on $F$, $F'$, $G'$, $U$, $V$ and $\mu$, and where $\text{h}^\text{II}_\nu$ depends on $F$, $G$, $G'$, $V$, $W$ and $\nu$.

\subsubsection{Construction of the first Haas transformation}
\label{SecHaasFirst}

Given $T\in\Ob\Al$, we let $\EEGrp{F}{G}(T)$ be defined via an injective resolution $I_T$ of $T$ and via a CE-resolution $J_T$ of $FI_T$; cf.\ \S\ref{SecGrothDef}.

Given $T'\in\Ob\Al'$, we let $\EEGrp{F'}{G'}(T')$ be defined via an injective resolution $I'_{T'}$ of $T'$ and via a CE-resolution $J'_{T'}$ of $F'I'_{T'}$; cf.\ \S\ref{SecGrothDef}.

{\it We define $\text{\rm h}^\text{\rm I}_\mu$.} Let $X\in\Ob\Al$. By Remark \ref{RemEx3}, there is a unique morphism $I'_{UX} \lraa{h'X} UI_X$ in $\KK^{[0}(\Al')$ that maps to $1_{UX}$ under $\HH^0$. 
Let $J'_{UX}\lraa{h''X} VJ_X$ be the unique morphism in $\KKKP(\Bl')$ that maps to the composite morphism $\left(F'I'_{UX}\lrafl{25}{F'h'X} F'U I_X\lraa{\mu} VFI_X\right)\ru{6}$ in $\KK^{[0}(\Bl')$ under 
$\HH^0\big((-)^{\ast,-}\big)$; cf.\ Lemma \ref{LemCEnew1}. Let the {\it first Haas transformation} be defined by
\[
\left(\EEGrp{F'}{G'}\big(UX\big) \;\lrafl{28}{\text{h}^\text{I}_\mu X}\; \EEGrp{F}{G'\0 V}\big(X\big)\right)
\;\;\; :=\;\;\; \left(\EEI(G'J'_{UX})\mrafl{25}{\EEI(G'h''X)} \EEI(G'VJ_X)\right) \; .
\]

{\it We show that $\text{\rm h}^\text{\rm I}_\mu$ is a transformation.} Let $X\lraa{f} Y$ be a morphism in $\Al$. Let $I_X\lraa{f'} I_Y$ resolve $X\lraa{f} Y$. Let $J_X\lraa{f''} J_Y$ resolve 
$FI_X\lraa{f'} FI_Y$. Let $I'_{UX}\lraa{\w f'} I'_{UY}$ resolve $UX\lraa{Uf} UY$. Let $J'_{UX}\lraa{\w f''} J_{UY}$ resolve $F'I_{UX}\lraa{F'\w f'} F'I_{UY}\ru{5}$. The quadrangle
\[
\xymatrix{
UX\ar@{=}[r]\ar[d]_{Uf} & UX\ar[d]^{Uf} \\
UY\ar@{=}[r]            & UY            \\
}
\]
commutes in $\Al'$. Hence, by Remark \ref{RemEx3}, applied to $I'_{UX}$ and $UI_Y$, the resolved quadrangle
\[
\xymatrix{
I'_{UX}\ar[r]^{h'X}\ar[d]_{\w f'} & UI_X\ar[d]^{Uf'} \\
I'_{UY}\ar[r]_{h'Y}               & UI_Y
}
\]
commutes in $\KK^{[0}(\Al')$. Hence both quadrangles in
\[
\xymatrix@C=14mm{
F'I'_{UX}\ar[r]^{F'h'X}\ar[d]_{F'\w f'} & F'UI_X\ar[d]^{F'Uf'}\ar[r]^{\mu} & VFI_X\ar[d]^{VFf'} \\
F'I'_{UY}\ar[r]_{F'h'Y}                 & F'UI_Y              \ar[r]_{\mu} & VFI_Y \\
}
\]
commute in $\KK^{[0}(\Bl')$. By Lemma \ref{LemCEnew1}, applied to $J'_{UX}$ and $VJ_Y$, the outer quadrangle in the latter diagram can be resolved to the commutative quadrangle
\[
\xymatrix@C=34mm{
J'_{UX}\ar[r]^{h''X}\ar[d]_{\w f''} & VJ_X\ar[d]^{Vf''} \\
J'_{UY}\ar[r]_{h''Y}                & VJ_Y \\
}
\]
in $\KKKP(\Bl')$. Applying $\EEI\big(G'(-)\big)$ and employing the definitions of $\EEGrp{F'}{G'}\,$, $\EEGrp{F}{G'\0 V}$ and $\text{h}^\text{I}_\mu\,$, we obtain the sought commutative diagram
\[
\xymatrix@C=20mm{
\EEGrp{F'}{G'}(UX)\ar[r]^{\text{h}^\text{I}_\mu X}\ar[d]_{\EEGrp{F'}{G'}(Uf)} & \EEGrp{F}{G'\0 V}(X)\ar[d]^{\EEGrp{F}{G'\0 V}(f)} \\
\EEGrp{F'}{G'}(UY)\ar[r]_{\text{h}^\text{I}_\mu Y}                            & \EEGrp{F}{G'\0 V}(Y) \\
}
\]
in $\bo\bZiffp,\,\Cl'\bc$.

\subsubsection{Construction of the second Haas transformation}
\label{SecHaasSecond}

We maintain the notation of \S\ref{SecHaasFirst}.

Given $X\in\Ob\Al$, we let the {\it second Haas transformation} be defined by
\[
\left(\EEGrp{F}{G\0 V}\big(X\big) \;\lraa{\text{h}^\text{II}_\nu X}\; \EEGrp{F}{W\0 G}\big(X\big)\right)
\;\;\; :=\;\;\; \left(\EEIp(G'V J_X)\lrafl{25}{\EEIp(\nu)} \EEIp(WG J_X)\right) \; .
\]
It is a transformation since $\nu$ is.

\section{The first comparison} 
\label{SecFirstComp}

\subsection{The first comparison isomorphism}
\label{SecFirstCompIso}

Suppose given abelian categories $\Al$, $\Al'$ and $\Bl$ with enough injectives and an abelian \mb{category $\Cl$.}

Let $\Al\ti\Al'\lraa{F}\Bl$ be a biadditive functor. Let $\Bl\lraa{G}\Cl$ be an additive functor. 

Suppose given objects $X\in\Ob\Al$ and $X'\in\Ob\Al'$. Suppose the following properties to hold.

\begin{itemize}
\item[(a)\phantom{$'$}] The functor $F(-,X'):\Al\lra\Bl$ is left exact.
\item[(a$'$)] The functor $F(X,-):\Al'\lra\Bl$ is left exact.
\item[(b)\phantom{$'$}] The functor $G$ is left exact.
\item[(c)\phantom{$'$}] The object $X$ possesses a $\big(F(-,X'),\, G\big)$-acyclic resolution $A\in\Ob\CC^{[0}(\Al)$.
\item[(c$'$)] The object $X'$ possesses a $\big(F(X,-),\, G\big)$-acyclic resolution $A'\in\Ob\CC^{[0}(\Al')$.
\end{itemize}

Moreover, the resolutions appearing in (c) and (c$'$) are stipulated to have the following properties.

\begin{itemize}
\item[(d)\phantom{$'$}] For all $k\ge 0$, the quasiisomorphism $\Conc X\lra A$ is mapped to a quasiisomorphism $\Conc F(X,A'^k)\lra F(A,A'^k)$ under $F(-,A'^k)$.
\item[(d$'$)] For all $k\ge 0$, the quasiisomorphism $\Conc X'\lra A'$ is mapped to a quasiisomorphism $\Conc F(A^k,X')\lra F(A^k,A')$ under $F(A^k,-)$.
\end{itemize}

The conditions (d,\,d$'$) are e.g.\ satisfied if $F(-,A'^k)$ and $F(A^k,-)$ are exact for all $k\ge 0$. 

\begin{Theorem}[first comparison]
\label{ThFC1}
The proper Grothendieck spectral sequence for the functors $F(X,-)$ and $G$, evaluated at $X'$, is isomorphic to the proper Grothendieck spectral sequence for the 
functors $F(-,X')$ and $G$, evaluated at $X$; i.e.
\[
\EEGrp{F(X,-)}{G}(X') \;\iso\; \EEGrp{F(-,X')}{G}(X)
\]
in $\bo\bZiffp,\Cl\bc$.
\end{Theorem}

{\it Proof.} Let \mb{$J_A\,,\; J_{A'}\,,\; J_{A,A'}\;\in\;\Ob\CCCP(\Inj\Bl)$} be CE-resolutions of the complexes $F(A,X'),\, F(X,A'),\, \t F(A,A')\,\in\,\Ob\CC^{[0}(\Bl)$, respectively.

The quasiisomorphism $\Conc X\lra A$ induces a morphism $F(\Conc X,A')\lra F(A,A')$, yielding $F(X,A')\lra \t F(A,A')$, which is a quasiisomorphism since $\Conc F(X,A'^k)\lra F(A,A'^k)$ is 
a quasiisomorphism for all $k\ge 0$ by (d).

Choose a CE-resolution $J_{A'}\lra J_{A,A'}$ of $F(X,A')\lra \t F(A,A')$; cf.\ Remark~\ref{RemEx4}. Since the morphism $F(X,A')\lra \t F(A,A')$ is a quasiisomorphism, $J_{A'}\lra J_{A,A'}$ is a 
composite in $\CCCCE(\Inj\Bl)$ of a rowwise homotopism and a double homotopism; cf.\ Proposition \ref{PropCEqis9}. So is $G J_{A'}\lra G J_{A,A'}$. Hence, by Remark \ref{RemSSDC2} and by 
Lemma \ref{LemSSDC3}, we obtain an isomorphism of the proper spectral sequences of the first filtrations of the total complexes,
\[
\EEGrp{F(X,-)}{G}(X') \= \EEIp(GJ_{A'}) \;\lraiso\; \EEIp(GJ_{A,A'})\; .
\]
Likewise, we have an isomorphism
\[
\EEGrp{F(-,X')}{G}(X) \= \EEIp(GJ_A) \;\lraiso\; \EEIp(GJ_{A,A'})\; .
\]
We compose to an isomorphism $\EEGrp{F(X,-)}{G}(X') \,\lraiso\, \EEGrp{F(-,X')}{G}(X)$ as sought.\qed

\subsection{Naturality of the first comparison isomorphism}
\label{SecFirstCompIsoNat}

\bq
We narrow down the assumptions just as we have done for the introduction of the Haas transformations in \S\ref{SecHaasSit} in order to be able to express, in this narrower case, a naturality of the 
first comparison isomorphism from Theorem~\ref{ThFC1}.

\eq

Suppose given abelian categories $\Al$, $\Al'$ and $\Bl$ with enough injectives and an abelian \mb{category $\Cl$.}

Let $\Al\ti\Al'\lraa{F}\Bl$ be a biadditive functor. Let $\Bl\lraa{G}\Cl$ be an additive functor. 

Suppose that the following properties hold.

\begin{itemize}
\item[(a)\phantom{$'$}] The functor $F(-,X'):\Al\lra\Bl$ is left exact for all $X'\in\Ob\Al'$.
\item[(a$'$)] The functor $F(X,-):\Al'\lra\Bl$ is left exact for all $X\in\Ob\Al$.
\item[(b)\phantom{$'$}] The functor $G$ is left exact.
\item[(c)\phantom{$'$}] For all $X'\in\Ob\Al'$, the functor $F(-,X')$ carries injective objects to $G$-acyclic objects.
\item[(c$'$)] For all $X\in\Ob\Al$, the functor $F(X,-)$ carries injective objects to $G$-acyclic objects.
\item[(d)\phantom{$'$}] The functor $F(I,-)$ is exact for all $I\in\Ob\Inj\Al$.
\item[(d$'$)] The functor $F(-,I')$ is exact for all $I'\in\Ob\Inj\Al'$.
\end{itemize}

\begin{Proposition}
\label{PropFCNat1}
Suppose given $X\lraa{x} \w X$ in $\Al$ and $X'\in\Ob\Al'$. Note that we have a transformation $F(x,-) : F(X,-)\;\lra\; F(\w X,-)$.
The following quadrangle, whose vertical isomorphisms are given by the construction in the proof of Theorem {\rm\ref{ThFC1}}, commutes.
\[
\xymatrix@C=20mm{
\EEGrp{F(X,-)}{G}(X')\ru{-3} \ar[r]^{\text{\rm h}^\text{\rm I}_{F(x,-)} X'}\ar[d]_{\wr} & \EEGrp{F(\w X,-)}{G}(X')\ru{-3}\ar[d]^{\wr} \\
\EEGrp{F(-,X')}{G}(X)\ar[r]^{\EEGrp{F(-,X')}{G}(x)}                                     & \EEGrp{F(-,X')}{G}(\w X)                    \\
}
\]
\end{Proposition}

For the definition of the first Haas transformation $\text{h}^\text{I}_{F(x,-)}$, see \S\ref{SecHaasFirst}. 

An analogous assertion holds with interchanged roles of $\Al$ and $\Al'$.

{\it Proof of Proposition~{\rm\ref{PropFCNat1}}.} Let $I$ resp.\ $\w I$ be an injective resolution of $X$ resp.\ $\w X$ in $\Al$. Let $I\lraa{\h x} \w I$ be a resolution of $X\lraa{x} \w X\ru{5}$. 
Let $I'$ be an injective resolution of $X'$ in $\Al'$. 

Let $J_{I'}^{(X)}$ resp.\ $J_{I'}^{(\w X)}$ be a CE-resolution of $F(X,I')$ resp.\ $F(\w X,I')$.

Let $J_{I,I'}$ resp.\ $J_{\w I,I'}$ be a CE-resolution of $\t F(I,I')$ resp.\ $\t F(\w I,I')$.

Let $J_I$ resp.\ $J_{\w I}$ be a CE-resolution of $F(I,X')$ resp.\ $F(\w I,X')$. 

We have a commutative diagram
\[
\xymatrix@C=15mm{
F(X,I')\ar[r]^{F(x,I')}\ar[d]     & F(\w X,I')\ar[d] \\
\t F(I,I')\ar[r]^{\t F(\h x, I')} & \t F(\w I,I')    \\
F(I,X')\ar[r]^{F(\h x,X')}\ar[u]  & F(\w I,X')\ar[u] \\
} 
\]
in $\CC^{[0}(\Bl)$, hence in $\KK^{[0}(\Bl)$. By Proposition \ref{PropCEqis6}, it can be resolved to a commutative diagram
\[
\xymatrix{
J_{I'}^{(X)}\ar[r]\ar[d] & J_{I'}^{(\w X)}\ar[d] \\
J_{I,I'}\ar[r]           & J_{\w I,I'}           \\
J_I\ar[r]\ar[u]          & J_{\w I}\ar[u]        \\
} 
\]
in $\KKKP(\Bl)$. Application of $\EEIp\big(G(-)\big)$ yields the result; cf.\ Lemma \ref{LemSSDC3}.

\bq
 We refrain from investigating naturality of the first comparison isomorphism in $G$.

\eq

\section{The second comparison} 
\label{SecSecondComp}

\subsection{The second comparison isomorphism}
\label{SecSecondCompIso}

Suppose given abelian categories $\Al$ and $\Bl'$ with enough injectives, and abelian categories $\Bl$ and $\Cl$.

Let $\Al\lraa{F} \Bl'$ be an additive functor. Let $\Bl\ti\Bl'\lraa{G} \Cl$ be a biadditive functor. 

Suppose given objects $X\in\Ob\Al$ and $Y\in\Ob\Bl$. Let $B\in\Ob\CC^{[0}(\Bl)$ be a resolution of $Y$, i.e.\ suppose a quasiisomorphism $\Conc Y\lra B$ to exist. Suppose the 
following properties to hold.

\begin{itemize}
\item[(a)] The functor $F$ is left exact.
\item[(b)] The functor $G(Y,-)$ is left exact.
\item[(c)] The object $X$ possesses an $(F,\, G(Y,-))$-acyclic resolution $A\in\Ob\CC^{[0}(\Al)$.
\item[(d)] The functor $G(B^k,-)$ is exact for all $k\ge 0$.
\item[(e)] The functor $G(-,I')$ is exact for all $I'\in\Ob\Inj\Bl'$.
\end{itemize}

\begin{Remark}
\label{RemSC0_5}
Suppose given a morphism $D\lraa{f} D'$ in $\CCCP(\Cl)$. If $\HH^\ell(f^{-,\ast})$ is a quasiisomorphism for all $\ell\ge 0$, then $f$ induces an isomorphism
\[
\EEIp(D)\;\lrafl{27}{\EEIp(f)}\;\EEIp(D')
\]
of proper spectral sequences. 
\end{Remark}

{\it Proof.} By Lemma \ref{LemProper1}, it suffices to show that $\EEI(\al+1/\al-1\bby\al/\al-2)^{+k}(f)$ is an isomorphism for all $\al\in\Z$ and all $k\in\Z$. By Lemma \ref{LemSSDC1}, this amounts to 
isomorphisms $\HH^k\HH^\ell(f^{-,\ast})$ for all $k,\,\ell\,\ge\, 0$, i.e.\ to quasiisomorphisms $\HH^\ell(f^{-,\ast})$ for all $\ell\ge 0$.
\qed

Consider the double complex $G(B,F\!A)\in\Ob\CCCP(\Cl)$, where the indices of $B$ count rows and the indices of $A$ count columns. To the first filtration of its total complex, we can associate the 
proper spectral sequence $\EEIp(G(B,F\!A))\in\Ob\bo\bZiffp,\Cl\bc$.

\begin{Theorem}[second comparison]
\label{ThSC1}
The proper Grothendieck spectral sequence for the functors $F$ and $G(Y,-)$, evaluated at $X$, is isomorphic to $\EEIp(G(B,F\!A))$; i.e.\ 
\[
\EEGrp{F}{G(Y,-)}(X) \;\iso\; \EEIp(G(B,F\!A))
\]
in $\bo\bZiffp,\Cl\bc$.
\end{Theorem}

{\it Proof.} Let $J'\in\Ob\CCCP(\Inj\Bl')$ be a CE-resolution of $F\!A$. By definition, $\EEGrp{F}{G(Y,-)}(X) = \EEIp(G(Y,J'))$. By Remark \ref{RemSC0_5}, it suffices to find $D\in\Ob\CCCP(\Cl)$ 
and two morphisms of double complexes
\[
G(B,F\!A)\;\lraa{u}\; D\;\llaa{v}\; G(Y,J')
\]
such that $\HH^\ell(u^{-,\ast})$ and $\HH^\ell(v^{-,\ast})$ are quasiisomorphisms for all $\ell\ge 0$.

Given a complex $U\in\Ob\CC^{[0}(\Bl)$, recall that we denote by $\text{Conc}_2 U\in\Ob\CCCP(\Bl)$ the double complex whose row number $0$ is given by $U$, and whose other rows are zero. 

We have a diagram 
\[
G(B,\Conc_2 F\!A)\;\lra\; G(B,J') \;\lla\; G(\Conc Y,J')
\]
in $\CCCCP(\Cl)$. Let $\ell\ge 0$. Application of $\HH^\ell\big((-)^{-,=,\ast}\big)$ yields a diagram
$$
\HH^\ell\big(G(B,\Conc_2 F\!A)^{-,=,\ast}\big)\;\;\lra\;\; \HH^\ell\big(G(B,J')^{-,=,\ast}\big) \;\;\lla\;\; \HH^\ell\big(G(\Conc Y,J')^{-,=,\ast}\big) 
\leqno (\ast)
$$
in $\CCCP(\Cl)$. We have
\[
\HH^\ell\big(G(B,\Conc_2 F\!A)^{-,=,\ast}\big) \;\;\iso\;\; G\Big(B\,,\,\HH^\ell\big((\Conc_2 FA)^{-,\ast}\big)\Big) \;\=\; G\big(B,\Conc\HH^\ell(F\!A)\big)
\]
and
\[
\HH^\ell\big(G(B,J')^{-,=,\ast}\big) \;\;\iso\;\; G\big(B,\HH^\ell(J'^{-,\ast})\big)\; ,
\]
since the functor $G(B^k,-)$ is exact for all $k\ge 0$ by (d), or, since the CE-resolution $J$ is rowwise split. Since the CE-resolution $J'$ is rowwise split, we moreover have 
\[
\HH^\ell\big(G(\Conc Y,J')^{-,=,\ast}\big) \;\;\iso\;\; G\big(\Conc Y,\HH^\ell(J'^{-,\ast})\big)\; .
\]
So the diagram $(\ast)$ is isomorphic to the diagram
$$
G\big(B,\Conc\HH^\ell(F\!A)\big) \;\;\lra\;\;  G\big(B,\HH^\ell(J'^{-,\ast})\big) \;\;\lla\;\; G\big(\Conc Y,\,\HH^\ell(J'^{-,\ast})\big) \; ,
\leqno (\ast\ast)
$$
whose left hand side morphism is induced by the quasiisomorphism $\Conc\HH^\ell(F\!A)\lra \HH^\ell(J'^{-,\ast})$, and whose right hand side morphism is induced by the quasiisomorphism 
$\Conc Y\lra B$. 

By exactness of $G(B^k,-)$ for $k\ge 0$, the left hand side morphism of $(\ast\ast)$ is a rowwise quasiisomorphism. Since $\HH^\ell(J'^{k,\ast})$ is injective, the functor
$G(-,\HH^\ell(J'^{k,\ast}))$ is exact by (e), and therefore the right hand side morphism of $(\ast\ast)$ is a columnwise quasiisomorphism. Thus an application of $\t$ to $(\ast\ast)$ yields two
quasiisomorphisms; cf.\ \S\ref{SecTotal}. Hence, also an application of $\t$ to $(\ast)$ yields two quasiisomorphisms in the diagram
\[
\t\HH^\ell\big(G(B,\Conc_2 F\!A)^{-,=,\ast}\big)\;\;\lra\;\; \t\HH^\ell\big(G(B,J')^{-,=,\ast}\big) \;\;\lla\;\; \t\HH^\ell\big(G(\Conc Y,J')^{-,=,\ast}\big) \; .
\]
Note that $\t\0\HH^\ell\big((-)^{-,=,\ast}\big) = \HH^\ell\big((-)^{-,\ast}\big)\0\t_{1,2}$, where $\t_{1,2}$ denotes taking the total complex in the first and the second index of a triple
complex; cf.\ \S\ref{SecPlanewise}. Hence we have a diagram 
\[
\HH^\ell\Big(\big(\t_{1,2} G(B,\Conc_2 F\!A)\big)^{-,\ast}\Big)\;\;\lra\;\; \HH^\ell\Big(\big(\t_{1,2} G(B,J')\big)^{-,\ast}\Big) 
\;\;\lla\;\; \HH^\ell\Big(\big(\t_{1,2} G(\Conc Y,J')\big)^{-,\ast}\Big) 
\]
consisting of two quasiisomorphisms. This diagram in turn, is isomorphic to
\[
\HH^\ell\Big(G(B,F\!A)^{-,\ast}\Big)\;\;\lra\;\; \HH^\ell\Big((\t_{1,2} G(B,J'))^{-,\ast}\Big) 
\;\;\lla\;\; \HH^\ell\Big(\big(G(Y,J')\big)^{-,\ast}\Big) \; ,
\]
where the left hand side morphism is obtained by precomposition with the isomorphism $G(B,FA^k)\lraiso \t\Conc_1 G(B,FA^k) = (\t_{1,2} G(B,\Conc_2 F\!A))^{-,k}$, where $k\ge 0$; cf.\ \S\ref{SecTotal}.

Hence we may take
\[
\big(G(B,F\!A)\;\lraa{u}\; D\;\llaa{v} G(B,J') \big) \;\;\; :=\;\;\; \Big(G(B,F\!A)\;\lra\; \t_{1,2} G(B,J') \;\lla\; G(Y,J')\Big)\; .
\]
\qed

\subsection{Naturality of the second comparison isomorphism}
\label{SecNatuSecond}

\bq
Again, we narrow down the assumptions just as we have done for the introduction of the Haas transformations in \S\ref{SecHaasSit} to express a naturality of the second comparison isomorphism from 
Theorem~\ref{ThSC1}.

\eq

Suppose given abelian categories $\Al$ and $\Bl'$ with enough injectives, and abelian categories $\Bl$ and $\Cl$. Suppose given additive functors $\Al\lradoublefl{28}{-35}{F}{\w F}\Bl'$ 
and a transformation $F\lraa{\phi}\w F$. Let $\Bl\ti\Bl'\lraa{G}\Cl$ be a biadditive functor.

Suppose given a morphism $X\lraa{x}\w X$ in $\Al$ and an object $Y\in\Ob\Bl$. Let $B\in\Ob\CC^{[0}(\Bl)$ be a resolution of $Y$, i.e.\ suppose a quasiisomorphism $\Conc Y\lra B$ to exist. Suppose the 
following properties to hold.

\begin{itemize}
\item[(a)] The functors $F$ and $\w F$ are left exact and carry injective to $G(Y,-)$-acyclic objects.
\item[(b)] The functor $G(Y,-)$ is left exact.
\item[(c)] The functor $G(B^k,-)$ is exact for all $k\ge 0$.
\item[(d)] The functor $G(-,I')$ is exact for all $I'\in\Ob\Inj\Bl'$.
\end{itemize}

Let $A\lraa{a}\w A$ in $\CC^{[0}(\Inj\Al)$ be an injective resolution of $X\lraa{x}\w X$ in $\Al$. Note that we have a commutative quadrangle
\[
\xymatrix@C=14mm{
G(B,FA)\ar[r]^{G(B,\phi A)}\ar[d]_{G(B,Fa)} & G(B,\w FA)\ar[d]^{G(B,\w Fa)} \\
G(B,F\w A)\ar[r]^{G(B,\phi\w A)}            & G(B,\w F\w A)                 \\
}
\]
in $\CCCP(\Cl)$. 

Note that once chosen injective resolutions $A$ of $X$ and $\w A$ of $\w X$, the image of $G(B,Fa)$ in $\KKKP(\Cl)$ does not depend on the choice of the resolution $A\lraa{a}\w A$ of $X\lraa{x}\w X$, for 
$\smash{\CC^{[0}(\Al)\mrafl{25}{G(B,F(-))} \CCCP(\Cl)}\ru{5}$ maps an elementary split acyclic complex to an elementary horizontally split acyclic complex.

\begin{Lemma}
\label{LemNatSecIsoA}
The quadrangle
\[
\xymatrix@C=18mm{
\EEGrp{F}{G(Y,-)}(X)\ar[r]^{\EEGrp{F}{G(Y,-)}(x)}\ar[d]_\wr & \EEGrp{F}{G(Y,-)}(\w X)\ar[d]^\wr \\
\EEIp(G(B,FA))\ar[r]^{\EEIp(G(B,Fa))}                       & \EEIp(G(B,F\w A))                 \\
}
\]
commutes, where the vertical isomorphisms are those constructed in the proof of Theorem~{\rm\ref{ThSC1}}.
\end{Lemma}

{\it Proof.} Let $J'\lraa{\h a}\w J'$ be a CE-resolution of $FA\lraa{Fa} F\w A$. Consider the following commutative diagram in $\CCCP(\Cl)$.
\[
\xymatrix@C=14mm{
G(Y,J')\ar[r]^{G(Y,\h a)}\ar[d]            & G(Y,\w J')\ar[d]    \\
\t_{1,2} G(B,J')\ar[r]^{\t_{1,2}G(B,\h a)} & \t_{1,2} G(B,\w J') \\
G(B,FA)\ar[r]^{G(B,Fa)}\ar[u]              & G(B,F\w A)\ar[u]    \\
}
\]
An application of $\EEIp$ yields the result.\qed

\begin{Lemma}
\label{LemNatSecIsoB}
The quadrangle
\[
\xymatrix@C=22mm{
\EEGrp{F}{G(Y,-)}(X)\ar[r]^{\text{\rm h}^\text{\rm I}_\phi X}\ar[d]_\wr & \EEGrp{\w F}{G(Y,-)}(X)\ar[d]^\wr \\
\EEIp(G(B,FA))\ar[r]^{\EEIp(G(B,\phi A))}                       & \EEIp(G(B,\w F A))                \\
}
\]
commutes, where the vertical morphisms are those constructed in the proof of Theorem~{\rm\ref{ThSC1}}.
\end{Lemma}

For the definition of the first Haas transformation $\text{h}^\text{I}_{F(x,-)}$, see \S\ref{SecHaasFirst}. 

{\it Proof.} Let $J'\lraa{\h\phi}\breve{J}'$ be a CE-resolution of $FA\lraa{F\phi} \w FA$. Consider the following commutative diagram in $\CCCP(\Cl)$.
\[
\xymatrix@C=14mm{
G(Y,J')\ar[r]^{G(Y,\h\phi)}\ar[d]            & G(Y,\breve{J}')\ar[d]  \\
\t_{1,2} G(B,J')\ar[r]^{\t_{1,2}G(B,\h\phi)} & \t_{1,2}G(B\breve{J}') \\
G(B,FA)\ar[r]^{G(B,\phi A)}\ar[u]            & G(B,\w FA)\ar[u]       \\
}
\]
An application of $\EEIp$ yields the result.\qed

\bq
 We refrain from investigating naturality of the second comparison isomorphism in $Y$.

\eq

\section{Acyclic CE-resolutions} 
\label{SecAcycCE}

\bq
 We record {\sc Beyl'}s Theorem \bfcite{Be81}{Th.\ 3.4} (here Theorem \ref{ThAC1}) in order to document that it fits in our context. The argumentation is entirely due to 
 {\sc Beyl} \mb{\bfcite{Be81}{Sec.\ 3}}, so we do not claim any originality.

\eq

Let $\Al$, $\Bl$ and $\Cl$ be abelian categories. Suppose $\Al$ and $\Bl$ to have enough injectives. Let $\Al\lraa{F}\Bl\lraa{G}\Cl$ be left exact functors.

\subsection{Definition}

Let $T\in\Ob\CC^{[0}(\Bl)$. In this \S\ref{SecAcycCE}, a CE-resolution of $T$ will {\sf synonymously} (and not quite correctly) be called an {\it injective CE-resolution,} to emphasise the fact
that its object entries are injective.

We regard $\CC^{[0}(\Bl)$ as an exact category as in Remarks \ref{RemCEqis2_5} and \ref{RemCEqis3_5}.

\begin{Definition}
\label{DefD1}\rm
A double complex $B\in\CCCP(\Bl)$ is called a {\it $G$-acyclic CE-resolution} of $T$ if the following conditions are satisfied.

\begin{itemize}
\item[(1)] We have $\HH^0(B^{\ast,-}) \iso T$ and $\HH^k(B^{\ast,-}) \iso 0$ for all $k\ge 1$.
\item[(2)] The morphism of complexes $B^{k,\ast}\lra B^{k+1,\ast}$, consisting of vertical differentials of $B$, is a pure morphism for all $k\ge 0$.
\item[(3)] The object $\BB^\ell(B^{k,\ast})$ is $G$-acyclic for all $\,k,\,\ell\,\ge\, 0\,$.
\item[(4)] The object $\ZZ^\ell(B^{k,\ast})$ is $G$-acyclic for all $\,k,\,\ell\,\ge\, 0\,$.
\end{itemize}

A {\it $G$-acyclic CE-resolution} is a $G$-acyclic CE-resolution of some $T\in\Ob\CC^{[0}(\Bl)$.
\end{Definition}

From (3,\,4) and the short exact sequence $\ZZ^\ell(B^{k,\ast})\lra B^{k,\ell}\lra\BB^{\ell+1}(B^{k,\ast})$, we conclude that 
$B^{k,\ell}$ is $G$-acyclic for all $\,k,\,\ell\,\ge\, 0\,$.

From (3,\,4) and the short exact sequence $\BB^\ell(B^{k,\ast})\lra\ZZ^\ell(B^{k,\ast})\lra\HH^\ell(B^{k,\ast})$, we conclude that 
$\HH^\ell(B^{k,\ast})$ is $G$-acyclic for all $\,k,\,\ell\,\ge\, 0\,$.

\begin{Example}
\label{ExD3}\rm
An injective CE-resolution of $T$ is in particular a $G$-acyclic CE-resolution of $T$.
\end{Example}

Note that given $Y\in\Ob\CC(\Bl)$ and $\ell\in\Z$, we have $\ZZ^\ell GY\iso G\ZZ^\ell Y$, whence the universal property of the cokernel $\HH^\ell G Y$ of $GY^{\ell-1}\lra\ZZ^\ell GY$ induces 
a morphism $\HH^\ell GY\lra G\HH^\ell Y$. This furnishes a transformation $\HH^\ell(GX^{k,\ast}) \lraa{\theta X} G \HH^\ell(X^{k,\ast})$, natural in $X\in\Ob\CCCP(\Bl)$.

\begin{Remark}
\label{RemD2}
If $B$ is a $G$-acyclic CE-resolution, then $\HH^\ell(GB^{-,\ast}) \lraa{\theta B} G \HH^\ell(B^{-,\ast})$ is an isomorphism for all $\ell\ge 0$.
\end{Remark}

{\it Proof.} The sequences 
\[
\ba{rcccl}
G\BB^\ell(B^{k,\ast})     & \lra & G\ZZ^\ell(B^{k,\ast}) & \lra & G\HH^\ell(B^{k,\ast}) \\
G\ZZ^{\ell-1}(B^{k,\ast}) & \lra & G B^{k,\ell-1}        & \lra & G\BB^\ell(B^{k,\ast}) \\
\ea
\]
are short exact for $k,\, \ell\,\ge\, 0$ by $G$-acyclicity of $\BB^\ell(B^{k,\ast})$ resp.\ of $\ZZ^{\ell-1}(B^{k,\ast})$. In particular, the cokernel of $GB^{k,\ell-1}\lra G\ZZ^\ell(B^{k,\ast})$
is given by $G\HH^\ell(B^{k,\ast})$. \qed

\subsection{A theorem of Beyl}

Let $X\in\Ob\Al_{(F,G)}$. Let $A\in\Ob\CC^{[0}(\Al)$ be a $(F,G)$-acyclic resolution of $X$. Let $B\in\CCCP(\Bl)$ be a $G$-acyclic CE\nobreakdash-resolution of $F\!A$.

\begin{Theorem}[{\sc Beyl,} \bfcite{Be81}{Th.\ 3.4}]
\label{ThAC1}
We have an isomorphism of proper spectral sequences
\[
\EEGrp{F}{G}(X) \;\iso\; \EEIp(GB)
\]
in $\bo\bZiffp,\Cl\bc$.
\end{Theorem}

{\it Proof.} Since the proper Grothendieck spectral sequence is, up to isomorphism, independent of the choice of an injective CE-resolution, as pointed out in \S\ref{SecGrothDef}, our assertion is 
equivalent to the existence of an injective CE-resolution $J$ of $F\!A$ such that $\EEIp(GJ)\;\iso\; \EEIp(GB)$. So by Remark~\ref{RemSC0_5}, it suffices to show that there exists an injective 
CE-resolution $J$ of $F\!A$ and a morphism $B\lra J$ that induces a quasiisomorphism $\HH^\ell(GB^{-,\ast})\lra \HH^\ell(GJ^{-,\ast})$ for all $\ell\ge 0$. By Remark \ref{RemD2} and 
Example~\ref{ExD3}, it suffices to show that $G\HH^\ell(B^{-,\ast})\lra G\HH^\ell(J^{-,\ast})$ is a quasiisomorphism for all $\ell\ge 0$.

By the conditions (1,\,2) on $B$ and by $G$-acyclicity of $\HH^\ell(B^{k,\ast})$ for $k,\,\ell\,\ge\, 0$, the complex $\HH^\ell(B^{-,\ast})$ is a $G$-acyclic resolution of 
$\HH^\ell(F\!A)$; cf.\ Remark \ref{RemCEqis3}. 

By Remark \ref{RemEx1}, there exists $J\in\Ob\CCCP(\Inj\Bl)$ with vertical pure morphisms and split rows, and a morphism $B\lra J$ consisting rowwise of pure monomorphisms such that 
$\HH^k(B^{\ast,-})\lra\HH^k(J^{\ast,-})$ is an isomorphism of complexes for all $k\ge 0$. In particular, the composite $(\Conc_2 F\!A\lra B\lra J)$ turns $J$ into an injective
CE-resolution of $F\!A$.

Let \mb{$\ell\ge 0$}. Since $B$ is a $G$-acyclic and $J$ an injective CE-resolution of $F\!A$, both $\Conc\HH^\ell(F\!A)\lra \HH^\ell(B^{-,\ast})$ and $\Conc\HH^\ell(F\!A)\lra \HH^\ell(J^{-,\ast})$ are 
quasiisomorphisms. Hence $\HH^\ell(B^{-,\ast})\lra \HH^\ell(J^{-,\ast})$ is a quasiisomorphism, too. Now Lemma \ref{LemCEqis4} shows that $G\HH^\ell(B^{-,\ast})\lra G\HH^\ell(J^{-,\ast})$ is a 
quasiisomorphism as well.\qed

\section{Applications}
\label{SecAppl}

\bq
We will apply Theorems \ref{ThFC1} and \ref{ThSC1} in various algebraic situations. In particular, we will re-prove a theorem of Beyl; viz.\ Theorem~\ref{ThApp3} in \S\ref{SecApplLHS}.

\eq

In several instances below, we will make tacit use of the fact that a left exact functor between abelian categories respects injectivity of objects provided it has an exact left adjoint.

\subsection{A Hopf algebra lemma}
\label{SecHopfLemma}

\bq
We will establish Lemma~\ref{LemCrux} in \S\ref{SecLemma}, needed to prove an acyclicity that enters the proof of the comparison result Theorem~\ref{ThAppHopf3} in \S\ref{SecApplHSH} for Hopf algebra
cohomology, which in turn allows to derive comparison results for group cohomology and Lie algebra cohomology; cf.\ \S\S\,\ref{SecApplLHS}, \ref{SecApplHS}.

\eq

\subsubsection{Definition}
\label{SecHopfDef}

Let $R$ be a commutative ring. Write $\ts := \ts_R\,$. A {\it Hopf algebra over $R$} is an $R$-algebra $H$ together with $R$-algebra morphisms $H\lraa{\eps} R$ {\it (counit)} and 
$H\lraa{\De} H\ts H$ {\it (comultiplication)}, and an $R$-linear map $H\lraa{S} H$ {\it (antipode)} such that the following conditions (i--iv) hold. 

Write $x\De = \sum_i xu_i\ts xv_i$ for $x\in H$, where $u_i$ and $v_i$ are chosen maps from $H$ to $H$, and where $i$ runs over a suitable indexing set. 
Note that $\sum_i (r\cdot x+s\cdot y)u_i\ts (r\cdot x+s\cdot y)v_i = r\cdot\left(\sum_i x u_i\ts x v_i\right) + s\cdot\left(\sum_i y u_i\ts y v_i\right)$ for $x,\, y\,\in\, H$ and $r,\,s\,\in\, R$,
whereas $u_i$ and $v_i$ are not necessarily $R$-linear maps.

\bq
 The elegant Sweedler notation~\cite[\S 1.2]{Sw69} for the images under $\De (\De\ts 1)$ etc.\ led the author, being new to Hopf algebras, to confusion in a certain case. So we will express them in 
 these more naive terms.

\eq

Write $H\ts H\lraa{\nabla} H$, 
$x\ts y\lramaps x\cdot y$ and $R\lraa{\et} H$, $r\lramaps r\cdot 1_H$. Write $H\ts H\lraa{\ta} H\ts H$, $x\ts y\lramaps y\ts x$.

\begin{itemize}
\item[(i)] We have $\De (\eps\ts\id_H) = (x\lramaps 1_R\ts x)$, i.e.\ $\sum_i xu_i\eps\cdot xv_i = x$ for $x\in H$.
\item[(i$'$)] We have $\De (\id_H\ts\eps) = (x\lramaps x\ts 1_R)$, i.e.\ $\sum_i xu_i\cdot xv_i\eps = x$ for $x\in H$.
\item[(ii)] We have $\De (\id_H\ts\De) = \De (\De\ts\id_H)$, i.e.\ $\sum_{i,j} x u_i \ts x v_i u_j\ts x v_i v_j = \sum_{i,j} x u_i u_j\ts x u_i v_j\ts x v_i$ for $x\in H$.
\item[(iii)] We have $\De (S\ts \id_H) \nabla = \eps\et$, i.e.\ $\sum_i x u_i S\cdot x v_i = x\eps\cdot 1_H$ for $x\in H$.
\item[(iii$'$)] We have $\De (\id_H\ts S) \nabla = \eps\et$, i.e.\ $\sum_i x u_i\cdot x v_i S = x\eps\cdot 1_H$ for $x\in H$.
\item[(iv)] We have $S^2 = \id_H$.
\end{itemize}

In particular, imposing (iv), we stipulate a Hopf algebra to have an {\it involutive} antipode.

\subsubsection{Some basic properties}
\label{SecHopfBasic}

\bq
 In an attempt to be reasonably self-contained, we recall some basic facts on Hopf algebras needed for Lemma~\ref{LemCrux} below; cf.\ \bfcite{Sw69}{Ch.\ IV}, \bfcite{Ab77}{\S 2}, \bfcite{Mo93}{\S\S 1-3}.
 In doing so, we shall use direct arguments. 

\eq

Suppose given a Hopf algebra $H$ over $R$.

\begin{Remark}[\bfcite{Sw69}{Prop.\ 4.0.1}, \bfcite{Ab77}{Th.\ 2.1.4}, \bfcite{Mo93}{3.4.2}]
\label{RemOmni}\Absit

The following hold.
\begin{itemize}
\item[\rm (1)] We have $\sum_i (x\cdot y)u_i\ts (x\cdot y)v_i = \sum_{i,j} (xu_i\cdot yu_j)\ts (x v_i\cdot y v_j)$ for $x,\, y\,\in\, H$.
\item[\rm (2)] We have $1_H S = 1_H$.
\item[\rm (3)] We have $(x\cdot y)S = yS\cdot xS$ for $x,\, y\,\in\, H$.
\item[\rm (4)] We have $S\eps = \eps$.
\item[\rm (5)] We have $\De (S\ts S) \ta = S\De$, i.e.\ $\sum_i x u_i S\ts x v_i S = \sum_i x S v_i \ts x S u_i$ for $x\in H$.
\item[\rm (6)] We have $x\cdot y = \sum_i \smash{\left(\sum_j (x u_i) u_j\cdot y\cdot (x u_i) v_j S\right)}\cdot x v_i$ for $x,\, y\,\in\, H$.
\item[\rm (6$'$)] We have $y\cdot x = \sum_i x u_i \cdot \smash{\left(\sum_j (x v_i) u_j S\cdot y\cdot (x v_i) v_j\right)}$ for $x,\, y\,\in\, H$.
\item[\rm (7)] We have $\sum_i x v_i\cdot x u_i S = x\eps\cdot 1_H$ for $x\in H$.
\item[\rm (7$'$)] We have $\sum_i x v_i S\cdot x u_i = x\eps\cdot 1_H$ for $x\in H$.
\end{itemize}
\end{Remark}

{\it Proof.} Ad (1). Given $x,\, y\,\in\, H$, we obtain
\[
\sumn{i} (xy)u_i\ts (xy)v_i \= (xy)\De \= x\De\cdot y\De \= \sumn{i,j} (xu_i\cdot yu_j)\ts (x v_i\cdot y v_j) \; .
\]

Ad (2). Remarking that $1_H\De = 1_H\ts 1_H$, we obtain
\[
1_HS \= 1_H \De(S\ts\id_H)\nabla \;\aufgl{(iii)}\; 1_H\eps\cdot 1_H \= 1_H\; .
\]

Ad (3). Given $x,\, y\,\in\, H$, we obtain
\[
\barcl
(x\cdot y)S
& \aufgl{2$\,\ti\,$(i$'$)} & \sum_{i,k} (x u_i\cdot x v_i\eps\cdot y u_k\cdot  y v_k\eps) S \\
& \aufgl{(iii$'$)} & \sum_{i,j,k} (x u_i\cdot y u_k\cdot  y v_k\eps) S \cdot x v_i u_j\cdot x v_i v_j S \\
& \aufgl{(iii$'$)} & \sum_{i,j,k,\ell} (x u_i\cdot y u_k) S \cdot x v_i u_j\cdot y v_k u_\ell\cdot y v_k v_\ell S \cdot  x v_i v_j S \\
& \aufgl{2$\,\ti\,$(ii)} & \sum_{i,j,k,\ell} (x u_i u_j\cdot y u_k u_\ell)S \cdot x u_i v_j\cdot  y u_k v_\ell\cdot y v_k S \cdot x v_i S \\
& \aufgl{(1)} & \sum_{i,j,k} (x u_i \cdot y u_k) u_j S \cdot (x u_i\cdot y u_k)v_j \cdot y v_k S \cdot x v_i S \\
& \aufgl{(iii)} & \sum_{i,k} (x u_i \cdot y u_k) \eps \cdot y v_k S \cdot x v_i S \\
& = & \sum_{i,k} (y u_k \eps \cdot y v_k) S \cdot (x u_i\eps\cdot x v_i) S \\
& \aufgl{2$\,\ti\,$(i)} & yS \cdot xS \; . \\
\ea
\]

Ad (4). Note that $(y\eps\cdot z)\eps = y\eps\cdot z\eps = (y\cdot z)\eps$ for $y,\, z\,\in\, H$. Given $x\in H$, we obtain
\[
xS\eps
\;\aufgl{(i)}\; (\sumn{i} x u_i\eps\cdot xv_i)S\eps 
\= (\sumn{i} x u_i\eps\cdot xv_i S)\eps 
\= (\sumn{i} x u_i\cdot xv_i S)\eps
\;\aufgl{(iii$'$)}\; (x\eps\cdot 1_H)\eps
\= x\eps\; .
\]

Ad (5). Given $x\in H$, we obtain
\[
\barcl
x\De (S\ts S)\ta
& \aufgl{(i)} & \sum_i (x u_i\eps\cdot x v_i)\De (S\ts S)\ta \\
& = & \sum_i (x u_i\eps\cdot 1_H)\De\cdot xv_i\De (S\ts S)\ta \\
& \aufgl{(iii)} & \sum_{i,j} (x u_i u_j S \cdot x u_i v_j) \De\cdot x v_i \De (S\ts S)\ta \\
& = & \sum_{i,j} x u_i u_j S\De \cdot x u_i v_j \De \cdot x v_i \De (S\ts S)\ta \\
& \aufgl{(ii)} & \sum_{i,j} x u_i S\De \cdot x v_i u_j \De \cdot x v_i v_j \De (S\ts S)\ta \\
& = & \sum_{i,j,k,\ell} x u_i S\De \cdot (x v_i u_j u_k \ts x v_i u_j v_k) \cdot (x v_i v_j v_\ell S \ts x v_i v_j u_\ell S) \\
& = & \sum_{i,j,k,\ell} x u_i S\De \cdot (x v_i u_j u_k\cdot x v_i v_j v_\ell S \ts x v_i u_j v_k\cdot x v_i v_j u_\ell S) \\
& \aufgl{(ii)} & \sum_{i,j,k,\ell} x u_i S\De \cdot (x v_i u_j\cdot x v_i v_j v_k v_\ell S \ts x v_i v_j u_k\cdot x v_i v_j v_k u_\ell S) \\
& \aufgl{(ii)} & \sum_{i,j,k,\ell} x u_i S\De \cdot (x v_i u_j\cdot x v_i v_j v_k S \ts x v_i v_j u_k u_\ell\cdot x v_i v_j u_k v_\ell S) \\
& \aufgl{(iii$'$)} & \sum_{i,j,k} x u_i S\De \cdot (x v_i u_j\cdot x v_i v_j v_k S \ts x v_i v_j u_k \eps\cdot 1_H) \\
& = & \sum_{i,j,k} x u_i S\De \cdot (x v_i u_j\cdot (x v_i v_j v_k \cdot x v_i v_j u_k \eps)S\ts 1_H) \\
& \aufgl{(i)} & \sum_{i,j} x u_i S\De \cdot (x v_i u_j\cdot x v_i v_j S\ts 1_H) \\
& \aufgl{(iii$'$)} & \sum_i x u_i S\De \cdot (x v_i\eps\cdot 1_H\ts 1_H) \\
& = & \sum_i (x u_i\cdot x v_i\eps)S\De \\
& \aufgl{(i$'$)} & x S\De\; . \\
\ea
\]

Ad (6). Given $x,\, y\,\in\, H$, we obtain
\[
x\cdot y
\;\aufgl{(i$'$)}\; \sumn{i} x u_i\cdot y \cdot x v_i\eps 
\;\aufgl{(iii)}\; \sumn{i,j} x u_i\cdot y \cdot x v_i u_j S\cdot x v_i v_j 
\;\aufgl{(ii)}\; \sumn{i,j} x u_i u_j\cdot y \cdot x u_i v_j S\cdot x v_i \; . 
\]

Ad (6$'$). Given $x\in H$, we obtain
\[
y\cdot x
\;\aufgl{(i)}\; \sumn{i} x u_i\eps\cdot y \cdot x v_i 
\;\aufgl{(iii$'$)}\; \sumn{i,j} x u_i u_j\cdot x u_i v_j S\cdot y \cdot x v_i 
\;\aufgl{(ii)}\; \sumn{i,j} x u_i\cdot x v_i u_j S\cdot y \cdot x v_i v_j \; . 
\]

Ad (7). Given $x\in H$, we have 
\[
\sumn{i} x v_i\cdot x u_i S
\aufgl{(iv)}    \sumn{i} xS^2 v_i\cdot xS^2 u_i S 
\aufgl{(5)}     \sumn{i} xS u_i S\cdot xS v_i S^2 
\aufgl{(iv)}    \sumn{i} xS u_i S\cdot xS v_i 
\aufgl{(iii)}   xS\eps\cdot 1_H 
\aufgl{(4)}     x\eps\cdot 1_H \; . 
\]

Ad (7$'$). Given $x\in H$, we have
\[
\sumn{i} x v_i S\cdot x u_i 
\aufgl{(iv)}     \sumn{i} xS^2 v_i S\cdot xS^2 u_i  
\aufgl{(5)}      \sumn{i} xS u_i S^2\cdot xS v_i S 
\aufgl{(iv)}     \sumn{i} xS u_i \cdot xS v_i S
\aufgl{(iii$'$)} xS\eps\cdot 1_H 
\aufgl{(4)}      x\eps\cdot 1_H \; . 
\]
\qed

In the present \S\ref{SecHopfLemma}, we shall refer to the assertions Remark~\ref{RemOmni}.(1--7$'$) just by (1--7$'$).

\subsubsection{Normality}
\label{SecNormality}

Suppose given a Hopf algebra $H$ over $R$, and an $R$-subalgebra $K\tm H$. Suppose $H$ and $K$ to be flat as modules over $R$.

Note that $K\ts K\lra H\ts H$ is injective. We will identify $K\ts K$ with its image.

The $R$-subalgebra $K\tm H$ is called a {\it Hopf-subalgebra} if $K\De\tm K\ts K$ and $KS\tm K$. In this case, we may and will suppose the maps $u_i$ and $v_i$ to restrict to 
maps from $K$ to $K$. 

Suppose $K\tm H$ to be a Hopf-subalgebra. It is called {\it normal,} if for all $a\in K$ and all $x\in H$, we have 
\[
\sumn{i} x u_i\cdot a \cdot x v_i S\;\in\; K\Icm\text{and}\Icm\sumn{i} x u_i S\cdot a \cdot x v_i\;\in\; K\; .
\]

An ideal $I\tm H$ is called a {\it Hopf ideal} if $I\De\tm I\ts H + H\ts I$ (where we have identified $I\ts H$ and $H\ts I$ with their images in $H\ts H$), $I\eps = 0$ and $IS\tm I$. In this case, 
the quotient $H/I$ carries a Hopf algebra structure via 
\[
\ba{rclrcl}
H/I & \lraa{\eps}  & R \; ,          & x + I & \lramaps & x\eps \\ 
H/I & \lraa{\De}   & H/I\ts H/I\; ,  & x + I & \lramaps & \sum_i (x u_i + I)\ts (x v_i + I) \\ 
H/I & \lraa{S}     & H/I\; ,         & x + I & \lramaps & xS + I \; . \\
\ea
\]

Suppose $K\tm H$ to be a normal Hopf subalgebra. Write $K^+ := \Kern(K\lraa{\eps} R)$. By (6,\,6$'$,\,3,\,4) and by writing 
\[
k\De \= \big({\ncm\sum_i} (ku_i - ku_i\eps)\ts kv_i\big) + 1\ts k
\]
for $k\in K^+$, the ideal $H K^+ = K^+ H$ is a Hopf ideal in $H$.

\subsubsection{Some remarks and a lemma}
\label{SecLemma}

Suppose given a Hopf algebra $H$ over $R$ and a normal Hopf-subalgebra $K\tm H$. Suppose $H$ and $K$ to be flat as modules over $R$. 

Write $\b H := H/HK^+$. Given $x\in H$, write $\b x := x + HK^+\in \b H$ for its residue class. 

Let $N'$, $N$, $M$, $M'$ and $Q$ be $H$-modules. Let $P$ be an $\b H$-module, which we also consider as an $H$-module via $H\lra\b H$, $x\lramaps\b x$.

We write $\liu{K}{(N,M)} = \liu{K}{(N|_K,M|_K)}$ for the $R$-module of $K$-linear maps from $N$ to $M$.

\begin{Remark}
\label{RemMod}
Given $f\in\liu{R}{(N,M)}$ and $x\in H$, we define $x\cdot f\in\liu{R}{(N,M)}$ by
\[
[n](x\cdot f) \; :=\; \sumn{i} x u_i\cdot [x v_i S\cdot n] f
\]
for $n\in N$. This defines a left $H$-module structure on $\liu{R}{(N,M)}$.
\end{Remark}

\bq
 Formally, squared brackets mean the same as parentheses. Informally, squared brackets are to accentuate the arguments of certain maps.

\eq

{\it Proof.} We {\it claim} that $x'\cdot(x\cdot f) = (x'\cdot x)\cdot f$ for $x,\, x'\,\in\, H$. Suppose given $n\in N$. We obtain
\[
\barcl
[n](x'\cdot (x\cdot f))
& = & \sum_i x' u_i\cdot [x' v_i S\cdot n](x\cdot f) \\
& = & \sum_{i,j} x' u_i\cdot x u_j\cdot [x v_j S\cdot x' v_i S\cdot n]f \\
& \aufgl{(3)} & \sum_{i,j} (x' u_i\cdot x u_j)\cdot [(x' v_i \cdot x v_j)S\cdot n]f \\
& \aufgl{(1)} & \sum_i (x'\cdot x) u_i\cdot [(x'\cdot x)v_i S\cdot n]f \\
& = & [n]((x'\cdot x)\cdot f) \; . \\
\ea
\]

We {\it claim} that $1_H\cdot f = f$. Suppose given $n\in N$. We obtain
\[
[n](1_H\cdot f)
\=  \sumn{i} 1_H u_i\cdot [1_H v_i S\cdot n]f 
\= 1_H \cdot [1_H S\cdot n] f 
\aufgl{(2)}  [n]f\; , 
\]
remarking that $1_H\De = 1_H\ts 1_H\,$. 
\qed

\bq
 I owe to {\sc G.\ Hi\ss} the hint to improve a previous weaker version of Corollary~\ref{RemCrux} below by means of the following Remark~\ref{RemFix}.

\eq

Denote by 
\[
M^K \; := \; \{ m\in M\; :\; \text{$a\cdot m = a\eps\cdot m$ for all $a\in K$}\} 
\]
the fixed point module of $M$ under $K$.

\begin{Remark}
\label{RemFix}
Letting $\b x\cdot m := x\cdot m$ for $x\in H$ and $m\in M^K$, we define an $\b H$-module structure on $M^K$.
\end{Remark}

{\it Proof.} The value of the product $\b x\cdot m$ does not depend on the chosen representative $x$ of $\b x$ since, given $y\in H$, $a\in K^+$ and $m\in M^K$, we have 
\[
y\cdot a\cdot m \= y\cdot a\eps\cdot m \= 0 \; .
\]
It remains to be shown that given $x\in H$ and $m\in M^K$, the element $x\cdot m$ lies in $M^K$. In fact, given $a\in K$, we obtain
\[
\barcl
a\cdot x\cdot m
& \aufgl{(6$'$)} & \sum_i x u_i \cdot \left(\sum_j (xv_i) u_j S\cdot a\cdot (x v_i) v_j\right)\cdot m \vspace*{1mm}\\
& = & \sum_i x u_i \cdot \left(\sum_j (xv_i) u_j S\cdot a\cdot (x v_i) v_j\right)\!\eps\cdot m \\
& = & \sum_{i,j} x u_i \cdot x v_i u_j S\eps\cdot a\eps\cdot x v_i v_j\eps \cdot m \\
&\aufgl{(4)} & \sum_{i,j} x u_i \cdot x v_i u_j\eps\cdot a\eps\cdot x v_i v_j\eps \cdot m \\
&\aufgl{(ii)} & \sum_{i,j} x u_i u_j \cdot x u_i v_j\eps\cdot a\eps\cdot x v_i\eps \cdot m \\
&\aufgl{(i$'$)} & \sum_{i} x u_i\cdot a\eps\cdot x v_i\eps \cdot m \\
&\aufgl{(i$'$)} & a\eps \cdot x \cdot m \; .\\
\ea
\]
\qed

\begin{Remark}
\label{RemHomAsFix}
We have $\left(\!\liu{R}{(N,M)}\right)^K = \liu{K}{(N,M)}$, as subsets of $\liu{R}{(N,M)}$.
\end{Remark}

{\it Proof.} The module $\left(\liu{R}{(N,M)}\right)^K$ consists of the $R$-linear maps $N\lraa{f} M$ that satisfy
\[
\sumn{i} x u_i\cdot [x v_i S\cdot n]f \= x\eps\cdot [n]f \; . 
\]
for $x\in H$ and $n\in N$. The module $\liu{K}{(N,M)}$ consists of the $R$-linear maps $N\lraa{f} M$ that satisfy
\[
[x\cdot n]f \= x\cdot [n]f
\]
for $x\in H$ and $n\in N$. By (iii$'$), we have $\left(\!\liu{R}{(N,M)}\right)^K \om \liu{K}{(N,M)}$.

It remains to show that $\left(\!\liu{R}{(N,M)}\right)^K \tm \liu{K}{(N,M)}$. Given $f\in \left(\!\liu{R}{(N,M)}\right)^K$, $x\in H$ and $n\in N$, we obtain
\[
\barcl
x\cdot [n]f
& \aufgl{(i$'$)} & \sum_i x u_i\cdot x v_i\eps\cdot [n]f \\
& = & \sum_i x u_i\cdot [x v_i\eps\cdot n]f \\
& \aufgl{(iii)} & \sum_{i,j} x u_i\cdot [x v_i u_j S\cdot x v_i v_j\cdot n]f \\
& \aufgl{(ii)} & \sum_{i,j} x u_i u_j\cdot [x u_i v_j S\cdot x v_i\cdot n]f \\
& = & \sum_i x u_i \eps\cdot [x v_i\cdot n]f \\
& \aufgl{(i)} & [x\cdot n]f \; .\\
\ea
\]
\qed

\begin{Corollary}
\label{RemCrux}
Given $f\in\liu{K}{(N,M)}$ and $x\in H$, we define $\b x\cdot f\in\liu{K}{(N,M)}$ by
\[
[n](\b x\cdot f) \; :=\; \sumn{i} x u_i\cdot [x v_i S\cdot n] f
\]
for $n\in N$. This defines a left $\b H$-module structure on $\liu{K}{(N,M)}$.
\end{Corollary}

{\it Proof.} By Remark~\ref{RemMod}, we may apply Remark~\ref{RemFix} to $\liu{R}{(N,M)}$. By Remark \ref{RemHomAsFix}, the assertion follows.\qed

\begin{Remark}
\label{RemCruxFun}
Given $f\in\liu{K}{(N,M)}$, $x\in H$, and $H$-linear maps $N'\lraa{\nu} N$, $M\lraa{\mu} M'$, we obtain 
\[
\nu (\b x\cdot f) \mu \= \b x\cdot (\nu f\mu)\; .
\]
\end{Remark}

{\it Proof.} Given $n'\in N'$, we obtain
\[
[n']\big(\nu (\b x\cdot f) \mu\big)
\= \big(\sumn{i} x u_i\cdot [x v_i S\cdot n'\nu] f\big) \mu 
\= \sumn{i} x u_i\cdot [x v_i S\cdot n'](\nu f \mu) 
\= [n'](\b x\cdot (\nu f\mu))\; . \;\; \qed
\]

\bq
 The following Lemma~\ref{LemCrux} has been suggested by the referee, and has been achieved with the help of {\sc G.~Carnovale.} It is reminiscent of \bfcite{Sch92}{Cor.\ 4.3}, but easier. It resembles a
 bit a Fourier inversion.

\eq

Note that the right $\b H$-module structure on $\b H$ induces a left $\b H$-module structure on $\liu{R}{(\b H,M)}$.

\begin{Lemma}
\label{LemCrux}
We have the following mutually inverse isomorphisms of $\b H$-modules.
\[
\barcl
\liu{K}{(H,M)}                                   & \lraisoa{\Phi} & \liu{R}{(\b H,M)} \\
f                                                & \lramaps       & (\b x\ramaps\sum_i x u_i\cdot [x v_i S] f) \vspace*{2mm} \\ 
\liu{K}{(H,M)}                                   & \llaisoa{\Psi} & \liu{R}{(\b H,M)} \\
(x\ramaps\sum_j x v_j\cdot [\,\ol{x u_j S}\,] g) & \llamaps       & g \\ 
\ea
\]
\end{Lemma}

{\it Proof.} We {\it claim} that $\Phi$ is a welldefined map. We have to show that $f\Phi$ is welldefined, i.e.\ that its value at $\b x$ does not depend on the representing element $x$.
Suppose given $y\in H$ and $a\in K^+$. We obtain
\[
\barcl
\sum_i (ya) u_i\cdot [(ya) v_i S] f
& \aufgl{(1)} & \sum_{i,j} y u_i \cdot a u_j\cdot [(y v_i\cdot a v_j) S] f \\
& \aufgl{(3)} & \sum_{i,j} y u_i \cdot a u_j\cdot [a v_j S\cdot y v_i S] f \\
& = & \sum_{i,j} y u_i \cdot a u_j\cdot a v_j S\cdot [y v_i S] f \\
& \aufgl{(iii$'$)} & \sum_i y u_i \cdot a\eps \cdot [y v_i S] f \\
& = & 0 \; . \\
\ea
\]

We {\it claim} that $\Phi$ is $\b H$-linear. Suppose given $y\in H$ and $x\in H$. We obtain
\[
\barcl
[\b x]((\b y f)\Phi)
& = & \sum_i x u_i\cdot [x v_i S] (\b y f) \\
& = & \sum_{i,j} x u_i\cdot y u_j\cdot [y v_j S\cdot x v_i S]f \\
& \aufgl{(3)} & \sum_{i,j} x u_i\cdot y u_j\cdot [(x v_i\cdot y v_j)S]f \\
& \aufgl{(1)} & \sum_i (x\cdot y) u_i\cdot [(x\cdot y) v_i S]f \\
& = & [\b x](\b y (f\Phi)) \; . \\
\ea
\]

We {\it claim} that $\Psi$ is a welldefined map. We have to show that $g\Psi$ is $K$-linear. Suppose given $a\in K$ and $x\in H$. Note that $a u_i\in K$ for all $i$, whence also $a u_i S\in K$, and 
therefore $a u_i S \con_{HK^+} a u_i S\eps\cdot 1_H$. We obtain
\[
\barcl
[a\cdot x](g\Psi)
& = & \sum_j (a\cdot x) v_j\cdot [\,\ol{(a\cdot x) u_j S}\,] g \\
& \aufgl{(1)} & \sum_{i,j} a v_i\cdot x v_j\cdot [\,\ol{(a u_i\cdot x u_j) S}\,] g \\
& \aufgl{(3)} & \sum_{i,j} a v_i\cdot x v_j\cdot [\,\ol{x u_j S}\cdot\ol{a u_i S}\,] g \\
& = &           \sum_{i,j} a v_i\cdot x v_j\cdot [\,\ol{x u_j S}\cdot\ol{a u_i S\eps}\,] g \\
& \aufgl{(4)} & \sum_{i,j} a u_i\eps\cdot a v_i\cdot x v_j\cdot [\,\ol{x u_j S}\,] g \\
& \aufgl{(i)} & \sum_j a\cdot x v_j\cdot [\,\ol{x u_j S}\,] g \\
& = & a\cdot [x](g\Psi)\; . \\
\ea
\]

We {\it claim} that $\Phi\Psi = \id_{\!\liu{K}{(H,M)}}$. Suppose given $x\in H$. We obtain
\[
\barcl
[x](f\Phi\Psi)
& = & \sum_j x v_j\cdot [\,\ol{x u_j S}\,](f\Phi) \\
& = & \sum_{i,j} x v_j\cdot x u_j S u_i\cdot [x u_j S v_i S] f \\
& \aufgl{(5)} & \sum_{i,j} x v_j\cdot x u_j v_i S\cdot [x u_j u_i S^2] f \\
& \aufgl{(iv)} & \sum_{i,j} x v_j\cdot x u_j v_i S\cdot [x u_j u_i] f \\
& \aufgl{(ii)} & \sum_{i,j} x v_j v_i\cdot x v_j u_i S\cdot [x u_j] f \\
& \aufgl{(7)} & \sum_j x v_j \eps\cdot [x u_j] f \\
& \aufgl{(i)} & [x] f \; . \\
\ea
\]

We {\it claim} that $\Psi\Phi = \id_{\!\liu{R}{(\b H,M)}}$. Suppose given $x\in H$. We obtain
\[
\barcl
[\b x](g\Psi\Phi)
& = & \sum_i x u_i\cdot [x v_i S] (g\Psi) \\
& = & \sum_{i,j} x u_i\cdot x v_i S v_j\cdot [\,\ol{x v_i S u_j S}\,] g \\
& \aufgl{(5)} & \sum_{i,j} x u_i\cdot x v_i u_j S\cdot [\,\ol{x v_i v_j S^2}\,] g \\
& \aufgl{(iv)} & \sum_{i,j} x u_i\cdot x v_i u_j S\cdot [\,\ol{x v_i v_j}\,] g \\
& \aufgl{(ii)} & \sum_{i,j} x u_i u_j\cdot x u_i v_j S\cdot [\,\ol{x v_i}\,] g \\
& \aufgl{(iii$'$)} & \sum_i x u_i\eps\cdot [\,\ol{x v_i}\,] g \\
& \aufgl{(i)} & [\b x] g \; .\\
\ea
\]

Finally, it follows by $\b H$-linearity of $\Phi$ and by $\Psi = \Phi^{-1}$ that $\Psi$ is $\b H$-linear.\qed

The tensor product $N\ts M$ is an $H$-module via $\De$. Note that $R$ is an $H$-module via $\eps$. Note that $R\ts M \iso M\iso M\ts R$ as $H$\nobreakdash-modules by (i,\,i$'$).

\begin{Remark}[cf.\ \bfcite{Ben91}{Lemma 3.5.1}]
\label{RemBens}
We have mutually inverse isomorphisms of $R$-modules
\[
\barcl
\liu{\b H}{(P,\liu{K}{(Q,M)})}  & \lraisoa{\al} & \liu{H}{(P\ts Q, M)}     \\
f                               & \lramaps      & (p\ts q \ramaps [q](pf)) \\
\liu{\b H}{(P,\liu{K}{(Q,M)})}  & \llaisoa{\be} & \liu{H}{(P\ts Q, M)}     \\
(p\ramaps (q\ramaps [p\ts q]g)) & \llamaps      & g \; ,                   \\
\ea
\]
natural in $P\in\Ob\b H\!\Modl$, $Q\in \Ob H\!\Modl$ and $M\in\Ob H\!\Modl$.
\end{Remark}

{\it Proof.} We {\it claim} that $\al$ is welldefined. We have to show that $f\al$ is $H$-linear. Suppose given $x\in H$. We obtain
\[
\barcl
x \cdot (p\ts q)
& = & \sum_i \ol{x u_i}\cdot p \ts x v_i \cdot q \\
& \lramapsa{f\al} & \sum_i [x v_i\cdot q]((\ol{x u_i}\cdot p) f) \\
& = & \sum_i [x v_i\cdot q](\ol{x u_i}\cdot (p f)) \\
& = & \sum_{i,j} x u_i u_j\cdot [x u_i v_j S\cdot x v_i\cdot q](p f) \\
& \aufgl{(ii)} & \sum_{i,j} x u_i\cdot [x v_i u_j S\cdot x v_i v_j \cdot q](p f) \\
& \aufgl{(iii)} & \sum_i x u_i\cdot [x v_i \eps \cdot q](p f) \\
& \aufgl{(i$'$)} & x \cdot [q](p f) \\
& = & x\cdot [p\ts q](f\al) \; . \\
\ea
\]

We {\it claim} that $\be$ is welldefined. First, we have to show that $[p](g\be)$ is $K$-linear. Suppose given $a\in K$. We obtain
\[
a\cdot q\;\lramapsfl{30}{\;[p](g\be)}\; [p\ts a\cdot q]g \;\aufgl{(i)}\;\sumn{i} [\ol{a u_i\eps} \cdot p \ts a v_i\cdot q]g \= \sumn{i} [\ol{a u_i} \cdot p \ts a v_i\cdot q]g \= a\cdot [p\ts q]g \; .
\]
Second, we have to show that $g\be$ is $\b H$-linear. Suppose given $x\in H$. We obtain
\[
\barcl
\b x\cdot p
& \lramapsa{g\be} & (q\ramaps [\b x\cdot p\ts q]g) \\
& \aufgl{(i)} & (q\ramaps \sum_i [\ol{x u_i\cdot x v_i\eps}\cdot p\ts q)]g) \\
& \aufgl{(iii$'$)} & (q\ramaps \sum_{i,j} [\ol{x u_i}\cdot p\ts x v_i u_j\cdot x v_i v_j S\cdot q)]g) \\
& \aufgl{(ii)} & (q\ramaps \sum_{i,j} [\ol{x u_i u_j}\cdot p\ts x u_i v_j\cdot x v_i S\cdot q)]g) \\
& = & (q\ramaps \sum_i x u_i\cdot [p\ts x v_i S\cdot q]g) \\
& = & \b x\cdot (q\ramaps [p\ts q]g) \; . \\
\ea
\]

Finally, $\al$ and $\be$ are mutually inverse.
\qed

\begin{Corollary}
\label{CorBens}
We have $\liu{\b H}{(P,M^K)} \iso \liu{\b H}{(P,\liu{K}{(R,M)})} \iso \liu{H}{(P, M)}$ as $R$\nobreakdash-modules, natural in $P$ and $M$.
\end{Corollary}

{\it Proof.} Note that $M \iso \liu{R}{(R,M)}$ as $H$-modules, whence $M^K\iso \liu{K}{(R,M)}$ as $\b H$-modules by Remarks~\ref{RemFix},~\ref{RemHomAsFix}. Now the assertion follows from 
Remark~\ref{RemBens}, letting $Q = R$. \qed

\subsection{Comparing Hochschild-Serre-Hopf with Grothendieck}
\label{SecApplHSH}

Let $R$ be a commutative ring. Suppose given a Hopf algebra $H$ over $R$ (with involutive antipode) and a normal Hopf-subalgebra $K\tm H$; cf.\ \S\ref{SecNormality}. 
Write $\b H := H/HK^+$. Suppose $H$, $K$ and $\b H$ to be projective as modules over $R$. Suppose $H$ to be projective as a module over $K$. 

Let $B\in\Ob\CC(H\!\Modl)$ be a projective resolution of $R$ over $H$. Let $\b B\in\Ob\CC(\b H\!\Modl)$ be a projective 
resolution of $R$ over $\b H$. Note that since $\b H$ is projective over $R$, $\b B|_R\in\Ob\CC(R\!\Modl)$ is a projective resolution of $R$ over $R$. Let $M$ be an $H$-module. 

By Corollary~\ref{RemCrux} and by Remark~\ref{RemCruxFun}, we have a biadditive functor
\[
\ba{rclcl}
(H\!\Modl)^\0 & \ti & H\!\Modl & \lraa{U} & \b H\!\Modl \\
(X            & ,   & X')      & \lramaps & U(X,X') := \liu{K}{(X,X')}\; . \\
\ea
\]
Write 
\[
\ba{rclcl}
(\b H\!\Modl)^\0 & \ti & \b H\!\Modl & \lraa{V} & R\!\Modl \\
(Y               & ,   & Y')         & \lramaps & V(Y,Y') := \liu{\b H}{(Y,Y')} \\
\ea
\]
for the usual Hom-functor.

In particular, we shall consider the functors 
\[
\ba{lclcl}
H\!\Modl & \lrafl{28}{U(R,-)} & \b H\!\Modl                   & \lrafl{28}{V(R,-)} & R\!\Modl                  \\
X        & \lramapsfl{}{}     & U(R,X)\iso X^K\hspace*{-15mm} &                    &                           \\
         &                    & Y                             & \lramaps           & V(R,Y)\iso Y^{\bar H}\; . \\
\ea
\]

On the other hand, we shall consider the double complex 
\[
D(M) \= D^{-,=}(M) \;:=\; V\big(\b B_{-}\,,\; U(B_{=}\,,\,M)\big) \= \liu{\b H}{\big(\b B_{-}\,, \liu{K}{(B_{=}\,,\,M)}\big)}\; .
\]
Note that $D(M)$ is isomorphic in $\CCCP(R\!\Modl)$ to $\liu{H}{\big(\b B_{-}\ts_R B_{=}\, , \; M\big)}$, naturally in $M$; cf.\ Remark~\ref{RemBens}.

\begin{Lemma}
\label{LemAppHopf2}
The $\b H$-module $U(H,M)$ is $V(R,-)$-acyclic.
\end{Lemma}

{\it Proof.} By Lemma \ref{LemCrux}, this amounts to showing that $\liu{R}{(\b H, M)}$ is $V(R,-)$-acyclic, which in turn amounts to showing that 
$V\big(\b B\,, \liu{R}{(\b H, M)}\big) = \liu{\b H}{\big(\b B\,, \liu{R}{(\b H, M)}\big)}$ has vanishing cohomology in degrees $\ge 1$. 
Now, 
\[
\liu{\b H}{\big(\b B\,, \liu{R}{(\b H, M)}\big)} \;\iso\; \liu{R}{(\b H\ts_{\b H}\b B\,,\, M)} \;\iso\; \liu{R}{(\b B\,,\, M)}\; ,
\]
whose cohomology in degree $i\ge 1$ is $\Ext_R^i(R,M) \iso 0$.\qed

\begin{Lemma}
\label{LemAppHopf2a}
Given a projective $H$-module $P$, the $\b H$-module $U(P,M)$ is $V(R,-)$-acyclic.
\end{Lemma}

{\it Proof.} It suffices to show that $U(\coprod_\Gamma H,\, M) \iso \prod_\Gamma U(H,\, M)$ is $V(R,-)$-acyclic for any indexing set $\Gamma$. By Lemma~\ref{LemAppHopf2}, it remains to 
be shown that $\RR^i V(R, \prod_\Gamma Y)$ is isomorphic to $\prod_\Gamma\RR^i V(R,Y)$ for a given $\b H$\nobreakdash-module $Y$ and for $i\ge 1$. Having chosen an injective resolution $J$ of $Y$, we may 
choose the injective resolution $\prod_\Gamma J$ of $\prod_\Gamma Y$. Then
\[
{\ncm
\RR^i V(R,\prod_\Gamma Y) \;\iso\; \HH^i V(R,\prod_\Gamma J) \;\iso\; \HH^i \prod_\Gamma V(R, J) \;\iso\; \prod_\Gamma \HH^i V(R, J) \;\iso\; \prod_\Gamma \RR^i V(R, Y)\; .
}
\]
\qed

\begin{Theorem}
\label{ThAppHopf3}
The proper spectral sequences 
\[
\EEIp(D(M)) \Icm\text{and}\Icm\EEGrp{U(R,-)}{\,V(R,-)}(M) 
\]
are isomorphic (in $\bo\bZiffp,\,R\!\Modl\bc$), naturally in $M\in\Ob H\!\Modl$.
\end{Theorem}

{\it Proof.} To apply Theorem~\ref{ThFC1} with, in the notation of \S\ref{SecFirstCompIso}, 
\[
\Big(\Al\ti\Al'\lrafl{25}{F}\Bl\lrafl{25}{G}\Cl\Big) \= \Big((H\!\Modl)^\0\ti H\!\Modl\lrafl{25}{U}\b H\!\Modl\;\lrafl{25}{V(R,-)}\;R\!\Modl\Big)\; ,
\]
and with $X = R$ and $X' = M$, we verify the conditions (a--d$'$) of loc.\ cit. in this case. 

Ad (c). We {\it claim} that $B$ is a $\big(U(-,M),\,V(R,-)\big)$-acyclic resolution of $R$. We have to show that $U(B_i,M)$ is $V(R,-)$-acyclic for $i\ge 0$; cf.~\S\ref{SecGrothDef}. Since $B_i$
is projective over $H$, this follows by Lemma \ref{LemAppHopf2a}. This proves the {\it claim.}

Ad (c$'$). Let $I$ be an injective resolution of $M$ over $H$. We {\it claim} that $I$ is a $\big(U(R,-),\,V(R,-)\big)$-acyclic resolution of $M$. We have to show that $U(R,I^i)$ is $V(R,-)$-acyclic
for $i\ge 0$. In fact, by Corollary~\ref{CorBens}, $U(R,I^i)$ is an injective $\b H$-module. This proves the {\it claim.}

Ad (d,\,d$'$). We {\it claim} that $U(B_i,-)$ and $U(-,I^i)$ are exact for $i\ge 0$; cf.\ \S\ref{SecFirstCompIso}. The former follows from $H$ being projective over $K$. The latter is 
a consequence of $I^i|_K$ being injective in $K\!\Modl$ by exactness of $K\!\Modl\;\lrafl{25}{H\ts_K -}\; H\!\Modl$. This proves the {\it claim.}

So an application of Theorem~\ref{ThFC1} yields
\[
\EEGrp{U(R,-)}{V(R,-)}(M) \;\iso\; \EEGrp{U(-,M)}{V(R,-)}(R) \; .
\]

To apply Theorem \ref{ThSC1} with, in the notation of \S\ref{SecSecondCompIso}, 
\[
\Big( \Al\lrafl{25}{F} \Bl'\; , \;\; \Bl\ti\Bl'\lrafl{25}{G} \Cl \Big) \= \Big( (H\!\Modl)^\0\;\;\lrafl{25}{U(-,M)}\;\b H\!\Modl\; , \;\; (\b H\!\Modl)^\0\ti\b H\!\Modl\lrafl{25}{V} \Cl \Big)\; ,
\]
and with $X = R$ and $Y = R$, we verify the conditions (a--e) of loc.\ cit. in this case. 

Ad (c). We have already remarked that $B$ is a $\big(U(-,M),\,V(R,-)\big)$-acyclic resolution of $R$.

Ad (d). As a resolution of $R$ over $\b H$, we choose $\b B$.

So an application of Theorem~\ref{ThSC1} yields
\[
\EEGrp{U(-,M)}{V(R,-)}(R) \;\iso\; \EEIp\Big( V\!\big(\b B_{-}\,,\; U(B_{=}\,,\,M)\big)\Big) \; .\vspace*{-3mm}
\]

Naturality in $M\in\Ob H\!\Modl$ remains to be shown. Suppose given $M\lraa{m}\w M$ in $H\!\Modl$. Note that the requirements of \S\ref{SecFirstCompIsoNat} are met. By Proposition~\ref{PropFCNat1},
with roles of $\Al$ and $\Al'$ interchanged, we have the following commutative quadrangle.
\[
\xymatrix@C=43mm{
\EEGrp{U(R,-)}{V(R,-)}(M)\ar[r]^{\EEGrp{U(R,-)}{V(R,-)}(m)}             & \EEGrp{U(R,-)}{V(R,-)}(\w M) \\
\EEGrp{U(-,M)}{V(R,-)}(R)\ar[r]^{\text{h}^\text{I}_{U(-,m)}R}\ar[u]^\wr & \EEGrp{U(-,\w M)}{V(R,-)}(R)\ar[u]_\wr \\
}
\]
Note that the requirements of \S\ref{SecNatuSecond} are met. By Lemma~\ref{LemNatSecIsoB}, we have the following commutative quadrangle.
\[
\xymatrix@C=30mm{
\EEGrp{U(-,M)}{V(R,-)}(R)\ar[r]^{\text{h}^\text{I}_{U(-,m)}R}\ar[d]_\wr                           & \EEGrp{U(-,\w M)}{V(R,-)}(R)\ar[d]^\wr                       \\
\EEIp\Big( V\!\big(\b B_{-}\,,\; U(B_{=}\,,\,M)\big)\Big)\ar[r]^{\EEIp(V(\b B_{-},\;U(B_{=},m)))} & \EEIp\Big( V\!\big(\b B_{-}\,,\; U(B_{=}\,,\,\w M)\big)\Big) \\
}
\]
\qed

\subsection{Comparing Lyndon-Hochschild-Serre with Grothendieck}
\label{SecApplLHS}

Let $R$ be a commutative ring. Let $G$ be a group and let $N\nrm G$ be a normal subgroup. Write $\b G := G/N$. Let $M$ be an $RG$-module. Write 
$\BBar_{G;R}\in\Ob\CC(RG\Modl)$ for the bar resolution of $R$ over $RG$, having $(\BBar_{G;R})_i = RG^{\ts (i+1)}$ for $i\ge 0$, the tensor product being taken over $R$.

Note that $RG$ is a Hopf algebra over $R$ via 
\[
\ba{rclrcl}
RG & \lraa{\De}  & RG\ts RG \; , & g & \lramaps & g\ts g \\
RG & \lraa{S}    & RG       \; , & g & \lramaps & g^{-1} \\ 
RG & \lraa{\eps} & R        \; , & g & \lramaps & 1\; ,  \\
\ea 
\]
where $g\in G$; cf.\ \S\ref{SecHopfDef}. Moreover, $RN$ is a normal Hopf subalgebra of $RG$ such that $RG/(RG)(RN)^+\iso R\b G$; cf.\ \S\ref{SecNormality}.

Note that $RG$, $RN$ and $R\b G$ are projective over $R$, and that $RG$ is projective over $RN$.

We have functors $RG\Modl\lrafl{28}{(-)^N} R\b G\Modl\lrafl{28}{(-)^{\b G}} R\!\Modl$, taking respective fixed points.

\begin{Theorem}[{\sc Beyl,} \bfcite{Be81}{Th.\ 3.5}]
\label{ThApp3}
The proper spectral sequences 
\[
\EEGrp{(-)^N}{\,(-)^{\b G}}(M) \Icm\text{and}\Icm\EEIp\Big(\liu{RG}{\big((\BBar_{\b G;R})_{-}\ts_R (\BBar_{G;R})_{=}\, , \; M\big)}\Big)
\]
are isomorphic (in $\bo\bZiffp,\,R\!\Modl\bc$), naturally in $M\in\Ob RG\Modl$.
\end{Theorem}

\bq
{\sc Beyl} uses his Theorem~\ref{ThAC1} to prove Theorem~\ref{ThApp3}. We shall re-derive it from Theorem~\ref{ThAppHopf3}, which in turn relies on the Theorems \ref{ThFC1} and \ref{ThSC1}.

\eq

{\it Proof.} This follows by Theorem~\ref{ThAppHopf3}. \qed

\subsection{Comparing Hochschild-Serre with Grothendieck}
\label{SecApplHS}

Let $R$ be a commutative ring. Let $\gfk$ be a Lie algebra over $R$ that is free as an $R$-module. Let $\nfk\nrm\gfk$ be an ideal such that $\nfk$ and $\b\gfk := \gfk/\nfk$ are free as 
$R$-modules. Let $M$ be a $\gfk$-module, i.e.\ a $\Ul(\gfk)$-module. Write $\BBar_{\gfk;R}\in\Ob\CC(\Ul(\gfk)\Modl)$ for the Chevalley-Eilenberg 
resolution of $R$ over $\Ul(\gfk)$, having $(\BBar_{\gfk;R})_i = \Ul(\gfk)\ts_R \wdg^i \gfk$ for $i\ge 0$; cf.\ \mb{\bfcite{CE56}{XIII.\S7}} or \mb{\bfcite{We95}{Th.\ 7.7.2}}.

Note that $\Ul(\gfk)$ is a Hopf algebra over $R$ via 
\[
\ba{rclrcl}
\Ul(\gfk) & \lraa{\De}  & \Ul(\gfk)\ts\Ul(\gfk) \; , & g & \lramaps & g\ts 1 + 1\ts g \\ 
\Ul(\gfk) & \lraa{S}    & \Ul(\gfk)             \; , & g & \lramaps & -g              \\
\Ul(\gfk) & \lraa{\eps} & R                     \; , & g & \lramaps & 0 \; ,          \\
\ea
\]
where $g\in\gfk$; cf.\ \S\ref{SecHopfDef}.

Note that $\Ul(\gfk)$, $\Ul(\nfk)$ and $\Ul(\b\gfk)$ are projective over $R$, and that $\Ul(\gfk)$ is projective over $\Ul(\nfk)$; cf.~\bfcite{We95}{Cor.\ 7.3.9}.

We have functors $\Ul(\gfk)\Modl\lrafl{28}{(-)^\nfk} \Ul(\b\gfk)\Modl\lrafl{28}{(-)^{\b\gfk}} R\!\Modl$, taking respective annihilated submodules; cf.~\bfcite{We95}{p.\ 221}.

\begin{Theorem}
\label{PropApp3_0_1}
The proper spectral sequences 
\[
\EEGrp{(-)^\nfk}{\,(-)^{\b\gfk}}(M) \Icm\text{and}\Icm \EEIp\Big(\liu{\Ul(\gfk)}{\big((\BBar_{\b\gfk;R})_{-}\ts_R (\BBar_{\gfk;R})_{=}\, , \; M\big)}\Big)
\]
are isomorphic (in $\bo\bZiffp,\,R\!\Modl\bc$), naturally in $M\in\Ob\Ul(\gfk)\Modl$.
\end{Theorem}

\bq
 Cf.\ {\sc Barnes,} \bfcite{Ba85}{Sec.\ IV.4, Ch.\ VII}.

\eq

{\it Proof.} This follows by Theorem~\ref{ThAppHopf3}. \qed

\subsection{Comparing two spectral sequences for a change of rings}
\label{SecAppChR}

\bq
 The following application is taken from \bfcite{CE56}{XVI.\S6}.

\eq

Let $R$ be a commutative ring. Let $A\lraa{\phi} B$ be a morphism of $R$-algebras. Consider the functors $A\!\Modl\lrafl{28}{\liu{A}{(B,-)}} B\!\Modl$ and 
$(B\!\Modl)^\0\ti B\!\Modl\lrafl{28}{\liu{B}{(-,=)}} R\!\Modl$.

Let $X$ be an $A$-module, let $Y$ be a $B$-module. 

\bq
 We shall compare two spectral sequences with $\EE_2$-terms $\Ext^i_B(Y,\Ext_A^j(B,X))$, converging to $\Ext^{i+j}_A(Y,X)$. If one views $X\indII_{\! A}^{\! B} := \liu{A}{(B,X)}$ as a way to induce
 from $A\!\Modl$ to $B\!\Modl$, this measures the failure of the Eckmann-Shapiro-type formula $\Ext_B^i(Y,X\indII_{\! A}^{\! B}) \auf{?}{\iso} \Ext_A^i(Y,X)$, which holds if $B$ is projective over $A$.

\eq

Let $I\in\Ob\CC^{[0}(A\!\Modl)$ be an injective resolution of $X$. Let $P\in\Ob\CC^{[0}(B\!\Modl)$ be a projective resolution of $Y$. 

\begin{Proposition}
\label{PropApp3_1}
The proper spectral sequences 
\[
\EEGrp{\liu{A}{(B,-)}}{\,\liu{B}{(Y,-)}}(X)\Icm\text{and}\Icm\EEIp\!\Big(\liu{B}{\big(P_{-},\,\liu{A}{(B,I^{=})}\big)}\Big) 
\]
are isomorphic (in $\bo\bZiffp,\,R\!\Modl\bc$).
\end{Proposition}

{\it Proof.} To apply Theorem \ref{ThSC1}, if suffices to remark that for each injective $A$-module $I'$, the $B$-module $\liu{A}{(B,I')}$ is injective, and thus $\liu{B}{(Y,-)}$-acyclic.
\qed

\begin{Remark} 
\label{RemApp3_2}\rm
The functor $\liu{A}{(B,-)}$ can be replaced by $\liu{A}{(M,-)}$, where $M$ is an $A$-$B$-bimodule that is flat over $B$.
\end{Remark}

\subsection{Comparing two spectral sequences for $\Ext$ and $\ts$}
\label{SecAppExtTs}

Let $R$ be a commutative ring. Let $S$ be a ring. Let $A$ be an $R$-algebra. Let $M$ be an \mb{$R$-$S$-bimodule}. Let $X$ and $X'$ be $A$-modules. Assume that $X$ is flat over $R$. Assume that 
$\Ext_R^i(M,X') \iso 0$ for $i\ge 1$. 

\begin{Example}
\label{ExApp3_2_5}\rm
Let $T$ be a discrete valuation ring, with maximal ideal generated by $t$. Let $R = T/t^\ell$ for some $\ell\ge 1$. Let $S = T/t^k$, where $1\le k\le \ell$. Let $G$ be a finite group, and let
$A = RG$. Let $M = S$. Let $X$ and $X'$ be $RG$-modules that are both finitely generated and free over $R$. 
\end{Example}

Consider the functors
\[
(A\!\Modl)^\0\ti A\!\Modl \;\lrafl{27}{\liu{A}{(-,=)}}\; R\!\Modl \;\lrafl{27}{\liu{R}{(M,-)}}\; S\Modl
\]

\begin{Proposition}
\label{PropApp3_3}
\hspace*{3mm} The proper Grothendieck spectral sequences 
\[
\EEGrp{\liu{A}{(X,-)}}{\,\liu{R}{(M,-)}}(X') \Icm\text{and}\Icm \EEGrp{\liu{A}{(-,X')}}{\,\liu{R}{(M,-)}}(X) 
\]
are isomorphic (in $\bo\bZiffp,\,S\Modl\bc$).
\end{Proposition}

Both have $\EE_2$-terms $\Ext^i_R\!\big(M,\,\Ext_A^j(X,X')\big)$ and converge to $\Ext_A^{i+j}(X\ts_R M,X')$. In particular, in the situation of Example \ref{ExApp3_2_5}, both have $\EE_2$-terms 
$\Ext_R^i\!\big(S,\,\Ext_{RG}^j(X,X')\big)$ and converge to $\Ext_{RG}^{i+j}(X/t^k,X')$.

{\it Proof of Proposition~{\rm\ref{PropApp3_3}}.} To apply Theorem \ref{ThFC1}, we comment on the conditions in \S\ref{SecFirstCompIso}.\vspace*{-2mm}
\begin{itemize}
\item[(c)\phantom{$'$}] Given a projective $A$-module $P$, we want to show that the $R$-module $\liu{A}{(P,X')}$ is \mb{$\liu{R}{(M,-)}$-acyclic}. We may assume that $P = A$, which is to be viewed as an 
$A$-$R$-bimodule. Now, we have $\Ext_R^i\!\big(M,\liu{A}{(A,X')}\big) \iso \Ext_R^i(M,X') \iso 0$ for $i\ge 1$ by assumption.
\item[(c$'$)] Given an injective $A$-module $I'$, the $R$-module $\liu{A}{(X,I')}$ is injective since $X$ is flat over $R$ by assumption.\qed
\end{itemize}

\subsection{Comparing two spectral sequences for $\text{\it $\El\!$xt}$ of sheaves}
\label{SecAppExtSh}

Let $T\lraa{f} S$ be a flat morphism of ringed spaces, i.e.\ suppose that
\[
\Ol_T\ts_{f^{-1}\Ol_S} - \; :\; f^{-1}\Ol_S\Modl\;\lra\;\Ol_T\Modl
\]
is exact. Consequently, $f^\ast:\Ol_S\Modl\lra\Ol_T\Modl$ is exact.

Given $\Ol_{\! S}$-modules $\Fl$ and $\Fl'$, we abbreviate $\liu{\Ol_{\! S}}{(\Fl,\Fl')} := \Hom_{\Ol_{\! S}}(\Fl,\Fl')\in \Ob R\!\Modl$ and 
$\liu{\Ol_{\! S}}{(\!(\Fl,\Fl')\!)} := \text{\it $\Hl$om}_{\Ol_{\! S}}(\Fl,\Fl')\in\Ob\Ol_{\! S}\Modl$.

Let $\Fl$ be an $\Ol_{\! S}$-module that has a locally free resolution $\Bl\in\Ob\CC(\Ol_{\! S}\Modl)$; cf.\ \mb{\bfcite{Har77}{Prop.~III.6.5}}. Let $\Gl\in\Ob\Ol_T\Modl$.
Let $\Al\in\Ob\CC^{[0}(\Ol_{\! T}\Modl)$ be an injective resolution of $\Gl$.

Consider the functors $\Ol_T\Modl\lraa{f_\ast}\Ol_{\! S}\Modl$ and $(\Ol_{\! S}\Modl)^\0\ti\Ol_{\! S}\Modl\mrafl{29}{\liu{\Ol_{\! S}}{(\!(-,=)\!)}}\Ol_{\! S}\Modl$.

\begin{Proposition}
\label{PropApp4} 
The proper spectral sequences 
\[
\EEGrp{f_\ast}{\,\liu{\Ol_{\! S}}{(\!( \Fl, -)\!)}}(\Gl) \Icm\text{and}\Icm \EEIp\big( \liu{\Ol_{\! S}}{(\!( \Bl_{-}, f_\ast\Al^{=})\!)} \big)
\]
are isomorphic (in $\bo\bZiffp,\,\Ol_{\! S}\Modl\bc$).
\end{Proposition}

In particular, both spectral sequences have $\EE_2$-terms $\text{\it $\El\!$xt}_{\Ol_{\! S}}^{\,i}\!\big(\Fl,\, (\RR^j f_\ast)(\Gl)\big)$ and converge to 
$(\RR^{i+j} \text{I}\!\Gamma_{\!\Fl})(\Gl)$, where $\text{I}\!\Gamma_{\!\Fl}(-) := \liu{\Ol_{\! S}}{(\!(\Fl,f_\ast(-))\!)} \iso f_\ast \liu{\Ol_T}{(\!(f^\ast \Fl,-)\!)}$. For 
example, if $S = \{\ast\}$ is a one-point-space and if we write $R := \Ol_S(S)$, then we can identify $\Ol_S\Modl = R\!\Modl$. If, in this case, $\Fl = R/rR$ for some $r\in R$, then 
$\text{I}\!\Gamma_{\! R/rR}(\Gl) \iso \Gamma(T,\Gl)[r] := \{ g\in\Gl(T)\; :\; rg = 0\}$.

{\it Proof of Proposition~{\rm\ref{PropApp4}}.} To apply Theorem \ref{ThSC1}, we comment on the conditions in \S\ref{SecSecondCompIso}.
\begin{itemize}
\item[(c)] Since $f_\ast$ maps injective $\Ol_T$-modules to injective $\Ol_{\! S}$-modules by flatness of $T\lraa{f} S$, the complex $\Al$ is an 
$\big(f_\ast,\, \liu{\Ol_{\! S}}{(\!( \Fl, -)\!)}\big)$-acyclic resolution of $\Gl$.
\item[(e)] If $\Il$ is an injective $\Ol_{\! S}$-module and $U\tm S$ is an open subset, then $\Il|_U$ is an injective $\Ol_U$\nobreakdash-module; cf.\ \bfcite{Har77}{Lem.~III.6.1}. Hence 
$\liu{\Ol_{\! S}}{(\!(-,\Il)\!)}$ turns a short exact sequence of $\Ol_{\! S}$-modules into a sequence that is short exact as a sequence of abelian presheaves, and hence a fortiori short exact as 
a sequence of $\Ol_{\! S}$-modules. In other words, the functor $\liu{\Ol_{\! S}}{(\!(-,\Il)\!)}$ is exact. \qed
\end{itemize}

\begin{footnotesize}

\parskip1.2ex

\vspace*{5mm}

\begin{flushright}
Matthias K\"unzer\\
Lehrstuhl D f\"ur Mathematik\\
RWTH Aachen\\
Templergraben 64\\
D-52062 Aachen \\
kuenzer@math.rwth-aachen.de \\
www.math.rwth-aachen.de/$\sim$kuenzer\\
\end{flushright}

\end{footnotesize}

\end{document}